\documentclass{article}

\usepackage{amsthm,amsmath,amssymb}
\usepackage{mathrsfs} 
\usepackage{graphicx}
\usepackage[utf8]{inputenc}
\usepackage[T1]{fontenc}
\usepackage[margin=1in]{geometry} 
\usepackage{esint} 

\usepackage[usenames,dvipsnames,svgnames,table]{xcolor}
\usepackage[colorlinks=true, pdfstartview=FitV,linkcolor=ForestGreen,citecolor=ForestGreen, urlcolor=blue]{hyperref}

\usepackage{hyperref}

\newtheorem{thm}{Theorem}[section]
\newtheorem{prop}[thm]{Proposition}

\newtheorem{lemma}[thm]{Lemma}
\newtheorem*{claim}{Claim}
\newtheorem{defn}[thm]{Definition}
\theoremstyle{remark}
\newtheorem{remark}[thm]{Remark}
\numberwithin{equation}{section}

\newcommand{\R}{\mathbb R}
\newcommand{\N}{\mathbb N}

\newcommand{\M}{\mathbf M}

\newcommand{\mB}{\mathcal{B}}
\newcommand{\eps}{\varepsilon}
\newcommand{\Sing}{\mathbf{S}}
\newcommand{\Pka}{\mathscr{P}^{m}}

\newcommand{\Pfinal}{\mathscr{P}^{\finalexp}}
\newcommand{\finalexp}{\ms}
\newcommand{\gammas}{\gamma_\ast}
\newcommand{\ms}{m_\ast}

\newcommand{\fin}{f_{\mathrm{in}}}

\newcommand{\dd} {\, \mathrm{d}}
\newcommand{\dt}{\, \mathrm{d} t}
\newcommand{\ds}{\, \mathrm{d} s}
\newcommand{\dw}{\, \mathrm{d} w}
\newcommand{\dv}{\, \mathrm{d} v}
\newcommand{\dz}{\, \mathrm{d} z}

\newcommand{\un}{\mathbf{1}}

\newcommand{\cI}{\mathcal{I}}

\newcommand{\ned}{N_{\mathrm{ED}}}

\newcommand{\Csob}{C_{\mathrm{Sob}}}

\newcommand{\cln}{c^{\, \mathrm{ln}}_1}
\newcommand{\clnb}{c^{\, \mathrm{ln}}_2}
\newcommand{\clnt}{c^{\, \mathrm{ln}}_3}

\newcommand{\Cinterpol}{C_{\mathrm{int}}}

\newcommand{\Cnonlin}{C_{\mathrm{nonlin}}}

\newcommand{\etad}{\eta_{\mathrm{D}}}
\newcommand{\etadg}{\eta_{\mathrm{DG}}}
\newcommand{\etadgp}{\eta_{\mathrm{DG+}}}

\newcommand{\fkp}{f^\kappa_+}

\newcommand{\fnikp}{f^{n_i,\kappa}_+}
\newcommand{\fnkp}{f^{n,\kappa}_+}

\newcommand{\kap}{\kappa}
\newcommand{\kbar}{\overline{\kap}}
\newcommand{\Gkp}{G^\kap_+}

\newcommand{\ka}{k}
\newcommand{\aout}{{a^{\mathrm{out}}_\delta}}
\newcommand{\anout}{{a^{\mathrm{out}}_n}}
\newcommand{\out}{{\, \mathrm{out}}}
\newcommand{\inn}{{\, \mathrm{in}}}
\newcommand{\Aout}{A^{\mathrm{out}}}

\newcommand{\ain}{{a^{\mathrm{in}}_\delta}}
\newcommand{\anin}{a_n^{\mathrm{in}}}

\newcommand{\lnp}{\ln_{+}}

\DeclareMathOperator*{\esssup}{ess-sup}

\title{\bfseries Local regularity  for the space-homogenous Landau
  equation with very soft potentials}

\author{François Golse, Cyril Imbert, Sehyun Ji and Alexis F. Vasseur}

\date{\today}

\begin{document}

\maketitle

\begin{abstract}
This paper deals with the space-homogenous Landau equation with very soft potentials, including the Coulomb case. This nonlinear equation is of parabolic type with diffusion matrix given by the convolution product of the solution with the matrix 
$a_{ij} (z)=|z|^\gamma (|z|^2 \delta_{ij} - z_iz_j)$ for $\gamma \in [-3,-2)$. We derive local truncated entropy estimates and use them to establish two facts. Firstly, we prove that the set of singular points (in time and velocity) for the weak solutions 
constructed as in [C. Villani, Arch. Rational Mech. Anal. 143 (1998), 273-307] has zero $\Pfinal$ parabolic Hausdorff measure with $\ms:=  \frac72 |2+\gamma|$. Secondly, we prove that if such a weak solution is axisymmetric, then it is smooth 
away from the symmetry axis. In particular, radially symmetric weak solutions are smooth away from the origin. 
\end{abstract}

\tableofcontents


\section{Introduction}


\subsection{Presentation of the problem and main results}\label{ss:MainRes}


We are concerned with the regularity of weak solutions $f\equiv f(t,v)\ge 0$ of the space-homogenous Landau equation with very soft potentials (see \cite{V98} on p. 280)
\begin{equation}\label{e:landau}
\partial_t f(t,v)  = Q (f,f)(t,v)\,,\qquad (t,v)\in(0,+\infty) \times \R^3\,.
\end{equation}
The Landau collision integral $Q(f,f)$ is defined as
\begin{align}\label{coll_oper} 
Q (f,f)(t,v):=\nabla_v\cdot \left( \int_{\R^3} a(v-w) (f(t,w) \nabla_vf(t,v) - f(t,v) \nabla_wf(t,w) ) \dw \right)\,, 
\end{align}
where $a(z) := |z|^{2+\gamma} \Pi (z)$ for $\gamma \in [-3,-2)$ in the case of very soft potentials, and 
\begin{equation}\label{Def-Pi}
\Pi(z):= I - \frac{z}{|z|}\otimes \frac{z}{|z|}\,,\qquad z\not=0
\end{equation}
is the orthogonal projection on $(\R z)^\perp$. The Landau equation is supplemented with an initial condition
\begin{equation}
\label{e:ic}
f(0,v) = \fin (v)\ge 0\,,\qquad v \in \R^3.
\end{equation}
The Landau equation \cite{Landau1936} in the Coulomb case, corresponding to $\gamma=-3$, is the fundamental collisional kinetic theory used in plasma physics. The Landau
collision integral replaces the Boltzmann collision integral, which diverges in the case of a Coulomb interaction between point particles: see \S 41 in \cite{LL10} and \cite{V98} on p. 280. In the Coulomb case, the Landau equation, henceforth referred
to as the Landau-Coulomb equation, becomes
\begin{equation}\label{e:landauNC}
\partial_t f(t,v)  = \mathrm{Tr}\left[(f(t,\cdot)\star a(v))\nabla^2_vf(t,v)\right]+8\pi f(t,v)^2\,,\qquad (t,v)\in(0,+\infty) \times \R^3\,.
\end{equation}

In \cite{V98}, Villani has constructed some kind of global weak solutions of the Landau equation with very soft potentials, called $H$-solutions (see section 6 in \cite{V98}, and section 7 for the Coulomb case) for all measurable initial data $f_{in}\ge 0$ 
a.e. on $\R^3$ satisfying
\begin{equation}\label{ini-bnd}
\int_{\R^3}f_{in}(v)(1+|v|^2+|\ln f_{in}(v)|)\dv<\infty\,.
\end{equation}
These solutions satisfy the conservation laws
\begin{equation}\label{ConsLaws}
\int_{\R^3}\left(\begin{matrix}1\\ v\\|v|^2\end{matrix}\right)f(t,v)\dv=\int_{\R^3}\left(\begin{matrix}1\\ v\\|v|^2\end{matrix}\right)f_{in}(v)\dv\,,\qquad t\ge 0\,,
\end{equation}
together with the following variant of Boltzmann's H-Theorem: for all $t\ge 0$
\begin{equation}\label{HThm}
\int_{\R^3}f\ln f(t,v)\dv\!+\!\int_0^t\int_{\R^3\times\R^3}\left|\Pi(v\!-\!w)(\nabla_v\!-\!\nabla_w)\sqrt{f(s,v)f(s,w)|v\!-\!w|^{2+\gamma}}\right|^2\dv\dw\ds\!\le\!\int_{\R^3}f_{in}\ln f_{in}(v)\dv.
\end{equation}
Desvillettes has proved in \cite{d15} that any $H$-solution $f$ of the Landau equation with $\gamma\in[-3,-2)$ satisfies
\begin{equation}\label{DesviProd}
\int_0^\infty\int_{\R^3}(1+|v|^2)^{\gamma/2}|\nabla_v\sqrt{f(t,v)}|^2\dv<\infty\,.
\end{equation}
One observes moreover that,  if $f(t,v)$ is a $H$-solution of the Landau equation with $\gamma \in[-3,-2)$, the rescaled function $f_\eps(t,v):=\eps^{\gamma+5}f(\eps^2t,\eps v)$ is also a $H$-solution of the Landau equation. In particular, 
in the case of the Landau-Coulomb equation ($\gamma=-3$), the rescaled function $f_\eps(t,v)=\eps^2f(\eps^2t,\eps v)$ is also a $H$-solution of the Landau-Coulomb equation. These properties suggest that $H$-solutions of the Landau-Coulomb 
equation could be analogous to (the square Euclidean norm of) Leray solutions of the Navier-Stokes equations for incompressible fluids in three space dimensions.

Whether $H$-solutions of the Landau equation with very soft potentials are smooth for all positive times, or develop singularities in finite time is an outstanding open problem in the analysis of kinetic equations: see section 1.3 (2) in chapter 2E of 
\cite{MR1942465}, and section 7 of \cite{luis-survey}. 

However, a most remarkable recent result by Guillen and Silvestre \cite{GuilSilv} claims the global propagation of regularity for solutions of the Landau equation for all $\gamma\in[-3,1]$ with $C^1$ initial data of finite Fisher information 
satisfying a Gaussian decay assumption: see Theorem 1.2 in \cite{GuilSilv}.

Conditional regularity results have been obtained on this problem. Silvestre has proved (in Theorem 3.8 of \cite{zbMATH06668299}) that $H$-solutions of the Landau-Coulomb equation in $L^\infty((0,+\infty);L^p(\R^3;(1+|v|)^k\dv)))$ with $p>3/2$ and 
$k>8$ are bounded and therefore regular (with the special case of radially decreasing solutions already treated in Theorem 1.1 of \cite{MR3599518}). Desvillettes, He and Jiang have obtained a new Lyapunov functional for the Landau-Coulomb 
equation in \cite{desvillettes2020new}, and proved the global existence and uniqueness of regular solutions for initial data that are sufficiently close to the centered, reduced Gaussian distribution --- or to a more general Maxwellian distribution function: 
see \eqref{DefMaxw} below for the definition of Maxwellian equilibria, and Theorem 1.2 (ii) in \cite{desvillettes2020new}. They also prove that $H$-solutions of the Landau-Coulomb equation cannot become singular after some finite time $T^*>0$ 
computed explicitly in terms of the initial data (see Theorem 1.2 (ii) in \cite{desvillettes2020new}).

Global $H$-solutions obtained by Villani's construction (see equation (49) in \cite{V98}), henceforth referred to as ``Villani solutions'', satisfy an additional ``truncated entropy'' inequality, which is a key ingredient in the regularity theory of such solutions. 
With M.P. Gualdani, we have studied the partial regularity in time of Villani solutions to the Landau-Coulomb equation. The set of times about which a Villani solution of the Landau-Coulomb equation is unbounded is proved in \cite{ggiv}  to be of Hausdorff 
dimension at most $1/2$. In \cite{Leray34}, the set of singular times of Leray's solutions of the Navier-Stokes equations for incompressible fluids in three space dimensions is proved to be also of Hausdorff dimension at most $1/2$. This is a further formal 
analogy between the Landau-Coulomb and the Navier-Stokes equations.

\smallskip
In the present paper, we pursue this analogy and first establish a partial regularity result \textit{\`a la} Caffarelli-Kohn-Nirenberg \cite{ckn} for the Landau equation. For $t_0>0$ and $v_0\in\R^3$, we call parabolic cylinder of size $0<r<\sqrt{t_0}$ 
the set $Q_r(t_0,v_0):=(t_0-r^2,t_0]\times B_r(v_0)$, where $B_r(v_0)$ is the Euclidean open ball of radius $r$ centered at $v_0$. The parabolic Hausdorff measure of dimension $m\in[0,5]$ of a Borel set $X\subset[0,+\infty)\times\R^3$ is
\[
\Pka (X)  := \lim_{\delta \to 0^+}  \inf \left\{ \sum_{i=1}^\infty r_i^m\text{ s.t. }X \subset \bigcup_{i=1}^\infty Q_{r_i} (t_i,v_i)\text{ with }r_i < \delta \right\}
\]
(see section 2.10 in chapter II of \cite{MR0257325}). Since $\text{diam}(Q_r(t_0,v_0))\le\sqrt{5}r$, the $m$-dimensional Hausdorff measure satisfies $\mathscr{H}^m(X)\le C\Pka(X)$ for some $C>0$ (see \cite{ckn} on p. 772). 

Next we define the notion of ``singular point'' for a $H$-solution of the Landau equation.

\bigskip
\begin{defn} 
A point $(t_0,v_0)\in(0,+\infty)\times\R^3$ is a \emph{singular point} of a $H$-solution $f$ of the Landau equation \eqref{e:landau} on $(0,+\infty)\times\R^3$ if $f\notin L^\infty(Q_r(t_0,v_0))$ for all $0\!<\!r\!<\!\sqrt{t_0}$.The set of singular points of 
$f$ is called $\Sing[f]$.
  
\end{defn}

Our first main result is a bound on the size of the singular set for Villani solutions of the Landau equation. The proof of this result is based on the De Giorgi method \cite{DeG} (see the proof of Lemma IV in \cite{DeG}).

\bigskip
\noindent
\textbf{Theorem A.} 
{\it 
Let $f$ be a global Villani solution of the Landau equation \eqref{e:landau} for $\gamma\in[-3,-2)$, with initial data \eqref{e:ic} satisfying \eqref{ini-bnd}. Then $\Pfinal(\Sing[f])=0$ with $\ms:=\tfrac72|2+\gamma|$.
}

\smallskip
In particular, radial Villani solutions $f$ of \eqref{e:landau} for $\gamma\in[-\tfrac{18}7,-2)$ with initial data \eqref{e:ic} satisfying \eqref{ini-bnd} have singular set $\Sing[f]\subset(0,+\infty)\times\{0\}$. Indeed, arguing as in Remark 3 on p. 808 of \cite{ckn}, 
any $(t_0,v_0)\in\Sing[f]$ with $v_0\not=0$ gives rise to a $2$-dimensional sphere in $\Sing[f]$, contradicting $\mathscr{H}^{\ms}(\Sing[f])\le C\Pfinal(\Sing[f])=0$ since $\gamma\in[-\tfrac{18}7,-2)$ implies $\ms\le 2$. Yet, this simples argument fails to 
settle the Coulomb case $\gamma=-3$.

\smallskip
Our second main result answers precisely the question of local regularity for Villani solutions of the Landau equation under appropriate symmetry assumptions.

\bigskip
\noindent
\textbf{Theorem B.} 
{\it 
Let $f$ be a global Villani solution of the Landau equation \eqref{e:landau} for $\gamma\in[-3,-2)$, with initial data \eqref{e:ic} satisfying \eqref{ini-bnd}. Assume that $f$ is axisymmetric, viz. there exists a measurable function $\mathcal F$ defined a.e.
on the set $(0,+\infty)\times(0,+\infty)\times\R$, together with $(\bar v,\omega)\in\R^3\times\mathbb{S}^2$ such that
$$
f(t,v)=\mathcal F(t,|(v-\bar v)\times\omega|,(v-\bar v)\cdot\omega)\quad\text{ for a.e. }(t,v)\in(0,+\infty)\times\R^3.
$$
Then the singular set of this axisymmetric solution satisfies $\Sing[f]\subset(0,+\infty)\times(\bar v+\R\omega)$.
}

\smallskip
\noindent
\textbf{Corollary C.} 
{\it 
Let $f$ be a global Villani solution of the Landau equation \eqref{e:landau} for $\gamma\in[-3,-2)$, with initial data \eqref{e:ic} satisfying \eqref{ini-bnd}. Assume that $f$ is radial, viz. there exists $\bar v\in\R^3$ and a measurable function $\mathfrak F$ 
defined a.e. on $(0,+\infty)\times(0,+\infty)$ such that
$$
f(t,v)=\mathfrak F(t,|v-\bar v|)\quad\text{ for a.e. }(t,v)\in(0,+\infty)\times\R^3.
$$
Then the singular set of this radial solution satisfies $\Sing[f]\subset(0,+\infty)\times\{\bar v\}$.
}

\smallskip
Corollary C is a direct consequence of Theorem B, but does not follow directly from partial regularity, i.e. Theorem A, for $\gamma\!\in\![-3,-\tfrac{18}7)$. The proof of Theorem B (and therefore of Corollary C) involves significantly different ideas, 
including a local bound on $(-\Delta_v)^{-1}f(t,v)$, a critical quantity which dominates in turn the diffusion matrix $f(t,\cdot)\star a$ of \eqref{e:landauNC}. 

To the best of our knowledge, Theorem B and Corollary C are the first regularity results on global weak solutions of the Landau-Coulomb equation for the largest class of initial data satisfying only the physically natural constraints \eqref{ini-bnd} 
of finite mass, energy and entropy, without any smallness condition on the quantity \eqref{ini-bnd}, and without any additional decay assumption on $f_{in}$. In particular, we do not know whether the breakthrough result by Guillen and Silvestre in
\cite{GuilSilv} can be applied to obtain global regularity for the class of solutions considered in Theorem B and Corollary C.


\subsection{Previous results and related works}


Besides the works already mentioned in the preceding section, there is a sizeable literature on the Landau equation \eqref{e:landau}. In this section, we only mention a few references directly related to our study.

\paragraph{General results on the Landau equation.} Although the Landau-Coulomb equation, i.e. \eqref{e:landau} with $\gamma=-3$, is the only physically relevant case, more general potentials are considered in the literature for purely 
mathematical purposes. The usual terminology is that $0< \gamma \le 1$ corresponds to hard potentials, $\gamma =0$ to the case of Maxwell molecules, and $-3\le\gamma <0$ to soft potentials. Soft potentials are further decomposed into 
moderately soft ($\gamma\in(-2,0)$) and very soft ($\gamma\in[-3,-2)$) potentials. The interested reader is referred to \cite{MR1942465} for a general introduction to the Landau equation, and to \cite{dvI,dv} for a complete discussion of 
the hard potential case. In the cases of Maxwell molecules and of moderately soft potentials, weak solutions of \eqref{e:landau} with initial data \eqref{e:ic} satisfying \eqref{ini-bnd} are known to become locally bounded for all positive times: 
see Theorem 1.1 in \cite{ gg2}. The construction of $H$-solutions in \cite{V98} already mentioned in the preceding section covers all the cases (i.e. $\gamma\in[-3,1]$). 

As a consequence of his lower bound on entropy production \eqref{DesviProd} (which is Theorem 3 in \cite{d15}) already mentioned above, Desvillettes has proved the propagation of moments for $H$-solutions (Proposition 4 in \cite{d15}) 
of the Landau-Coulomb equation. He also proves that $H$-solutions of the Landau equation are weak solutions in the traditional sense of this term (Corollary 1.1 of \cite{d15}).

Classical parabolic theory says that if a weak solution $f$ of the Landau equation \eqref{e:landau} is bounded in some parabolic cylinder $Q_r(t_0,v_0)$ with $t_0>r^2>0$ and $v_0\in\R^3$, then $f$ is smooth in a (smaller) parabolic cylinder, say 
on $Q_{r/2}(t_0,v_0)$. For the sake of being complete, we have sketched in Appendix A.1 a proof of this result. 

\paragraph{The Landau-Coulomb equation.} Fournier has proved in \cite{MR2718931} that $H$-solutions of the Landau-Coulomb equation in $L^1((0,T);L^\infty(\R^3))$ are unique. The large time behavior of $H$-solutions of the Landau-Coulomb
equation is studied in \cite{CDH17}. In \cite{MR1946444,CM17,cdh}, existence and uniqueness are studied in the close-to-equilibrium, non space-homogenous case. Other results about local existence and uniqueness or smoothing effects can 
be found in \cite{zbMATH06417110,gimv,HendersonSnelson2017,zbMATH06987670}. Concerning the possibility of finite time blow-up for the Landau-Coulomb equation, Bedrossian, Gualdani and Snelson \cite{bedrossian2021nonexistence} have 
proved that self-similar blow-up solutions do not exist as soon as the self-similar profile enjoys a mild decay assumption in the velocity variable. The reader is also referred to the recent survey article \cite{luis-survey} by Silvestre which contains 
in particular a conjecture about global boundedness of solutions of the Landau equation in the Coulomb case for bounded initial data.

\paragraph{The ``isotropic'' Landau-Coulomb model} In view of the considerable difficulty of the regularity  problem for the Landau-Coulomb equation, Krieger and Strain \cite{MR2901061} have considered a model equation, referred to as
``isotropic'', since it corresponds to \eqref{e:landau} where the matrix field $a(z)=|z|^{-1}\Pi(z)$ is replaced with $|z|^{-1}I$. The best result on this model is the global existence of radially decreasing regular solutions in \cite{MR3599518}, 
based on a barrier function argument. This quite remarkable result is perhaps the best argument in favor of global regularity for radial solutions of the Landau-Coulomb equation. See also Theorem 1.1 (2) of \cite{gg-almost-landau} for the 
isotropic analogue of \eqref{e:landau} in the non-Coulomb, very soft case, with $a((z)=|z|^{2+\gamma}\Pi(z)$ replaced with $|z|^{2+\gamma}I$: classical solutions exist for all positive times provided that $\gamma\in(\bar\gamma,-2]$ with 
$-3<\bar\gamma<-2.45$. Notice however that Maxwellian distribution functions, i.e. 
\begin{equation}\label{DefMaxw}
\mathcal M_{\rho,u,\theta}f(v):=\frac{\rho}{(2\pi\theta)^3}e^{-|v-u|^2/2\theta}
\end{equation}
for some $\rho,\theta>0$ and $u\in\R^3$, which are the equilibrium (i.e. time-independent) solutions of \eqref{e:landau}, of considerable importance in statistical mechanics, do not satisfy the ``isotropic'' Landau-Coulomb model. This suggests that, 
for all its similarities with the Landau-Coulomb equation, the isotropic Landau-Coulomb model is of questionable physical relevance.

\paragraph{The De Giorgi method and partial regularity.} Caffarelli, Kohn and Nirenberg \cite{ckn} have proved that the singular set of (suitable) weak solutions of the Navier-Stokes equations for incompressible fluids in three space dimensions
has zero $\mathscr{P}^1$ Hausdorff measure. The third author gave an alternative proof of this result in \cite{MR2374209} relying on arguments from De Giorgi's paper \cite{DeG}. Another application of the De Giorgi method to fluid mechanics 
is to be found in \cite{CV-qg}. The reader interested in an introduction to De Giorgi's ideas in \cite{DeG} with further applications is referred to the survey articles \cite{CV-survey,zbMATH06608022}. The truncated entropy used with the De Giorgi 
method in the present paper and in our earlier work with M. Gualdani \cite{ggiv} finds its origin in \cite{GouVa}.

\paragraph{The De Giorgi method and kinetic equations.} The De Giorgi or Moser iteration methods, initially devised for elliptic or parabolic equations, have been successfully applied to space inhomogeneous kinetic equations. The class of 
ultraparabolic equations, studied by the Italian school since the 1990s, contains kinetic equations of Fokker-Planck type. For instance, Moser's iteration method was  implemented by Pascucci and Polidoro \cite{PascucciPolidoro2004} in 
this context, resulting in the local boundedness of nonnegative subsolutions of ultraparabolic equations. De Giorgi's ideas have also proved useful  when applied to the Landau equation \cite{gimv} or to the Boltzmann equation without cutoff 
\cite{Imb_Silv2017}. The nonconstructive proof of H\"older's continuity in \cite{gimv} was recently revisited by Gu\'erand and the second author \cite{guerand2021logtransform} and Gu\'erand and Mouhot \cite{guerand2021quantitative}. 
To conclude, we mention that the De Giorgi method is also implemented in \cite{alonso2020giorgi,luis-stanley-close-to} to study the Boltzmann equation without cutoff in the close-to-equilibrium regime. 

\subsection{Main steps in the proofs of Theorems A and B}


Both the proofs of Theorems A and B are based on the De Giorgi method \cite{DeG,CV-survey} adapted to parabolic equations, and on the scale invariance properties of the Landau equation. In this section, we only consider the Coulomb
case $\gamma=-3$, which is the most interesting on physical grounds, and (probably) also the most difficult mathematically. Note that for other values of $\gamma $, the definitions and arguments below should be modified since the rescaled solution $f_\eps$ scales differently.

If $f\equiv f(t,v)$ is a solution of the Landau-Coulomb equation \eqref{e:landauNC}, and if $(t_0,v_0)\in(0,+\infty)\times\R^3$, then
\begin{equation}\label{f-lm}
f_{\lambda,\mu}(t,v):=\lambda f(t_0+\lambda t,v_0+\mu v)
\end{equation}
is also a solution of \eqref{e:landauNC} for all $\lambda,\mu>0$. We shall apply De Giorgi iterations to
\begin{equation}\label{DGFunc0}
\esssup_{t-r^2<s\le t}\int_{B_r(v)}(f(s,w)-\kap)_+\dw+\int_{Q_r(t,v)}\left|\nabla_w\sqrt{f(s,w)}\right|^2\un_{f(s,w)\ge\kap}dsdw
\end{equation}
for a sequence of decreasing $r$ and increasing $\kappa$. It is therefore essential that both terms in the expression above are scaled identically under the transformation $f\mapsto f_{\lambda,\mu}$ in \eqref{f-lm}. This singles out the case 
$\lambda=\mu^2$ and sugests considering the scaling transformation
\begin{equation}\label{f-eps}
f_\eps(t,v):=\eps^2f(t_0+\eps^2 t,v_0+\eps v)\,.
\end{equation}

As mentionned above, Villani solutions of the Landau equation satisfy a truncated entropy inequality, which is key to the De Giorgi nonlinearization procedure for the functional \eqref{DGFunc0}. First we define the function
\begin{equation}\label{Def-h+}
h_+(z):=z\lnp(z)-(z-1)_+\,,\quad\text{ where }\lnp(z):=\max(\ln z,0)\,,\qquad z>0\,.
\end{equation}
Pick $(t_0,v_0)\in(0,+\infty)\times\R^3$ and $\delta,\eps,\kap,r>0$, assuming that
$$
\eps\in(0,\min(\tfrac12,\sqrt{t_0}))\,,\quad\kap_\eps:=\eps^2\kap\in[1,2]\cap\mathbb Q\,,\quad r_\eps:=r/\eps\in(0,2]\,,\quad\delta_\eps:=\delta/\eps\in(0,1]\,.
$$
Then, for all Villani solutions $f$ of the Landau-Coulomb equation \eqref{e:landauNC} with initial data $f_{in}$ such that
\begin{equation}\label{ini-bnd-bis}
0<m_0\le\int_{\R^3}f_{in}(v)\dv\le M_0\,,\quad\int_{\R^3}|v|^2f_{in}(v)\dv\le E_0\,,\quad\int_{\R^3}f_{in}(v)|\ln f_{in}(v)|\dv\le H_0\,,
\end{equation}
where $m_0<M_0$ and $E_0,H_0$ are arbitrary positive constants, the scaled solution \eqref{f-eps} satisfies
\begin{equation}\label{TruncEntrScal}
\begin{aligned}
\esssup_{-r_\eps^2<t\le 0}\int_{B_{r_\eps}(0)}\kap_\eps h_+\left(\tfrac{f_\eps(t,v)}{\kap_\eps}\right)\dv+\int_{Q_{r_\eps}(0,0)}\tfrac{\left|\nabla_v(f_{\eps}(t,v)-\kap_\eps)_+\right|^2}{f_\eps(t,v)}\dt \dv
\\
\le C_0\int_{Q_{r_\eps+\delta_\eps}(0,0)}\left(\kap_\eps+\tfrac1{\delta_\eps^2}+\tfrac1{\delta_\eps^2}f_\eps(t,\cdot)\star\tfrac{1}{|\,\cdot\,|}(v)\right)
	f_\eps\left(\ln_+\left(\tfrac{f_\eps}{\kap_\eps}\right)+\ln_+\left(\tfrac{f_\eps}{\kap_\eps}\right)^2\right)(t,v)\dt \dv
\end{aligned}
\end{equation}
for some $C_0\equiv C'_0[|v_0|,m_0,M_0,E_0,H_0]\ge 1$ that is independent of $\eps$: see Lemma 3.5 below.

Following \cite{ckn}, we define a notion of ``suitable weak solution'' of the Landau-Coulomb equation \eqref{e:landauNC}; specifically, suitable weak solutions of \eqref{e:landauNC} are $H$-solutions satisfying an inequality that is a variant of 
Boltzmann's H-Theorem implying \eqref{TruncEntrScal}. Then we prove that any initial data satisfying \eqref{ini-bnd} launches at least one global suitable weak solution of \eqref{e:landauNC} with initial condition \eqref{e:ic}. This global suitable 
weak solution is obtained as a limit point of the approximation sequence defined in equation (49) of \cite{V98}. Any statement which holds true for suitable weak solutions applies therefore to Villani solutions (in the terminology of the previous 
section). The existence proof of suitable weak solutions occupies sections \ref{ss:approx-scheme}-\ref{sec:lei} and Appendix \ref{a:trunc-entropineq}, while the derivation of \eqref{TruncEntrScal} occupies section \ref{ss:SuitLocEntr}.

Applying the De Giorgi method to the left-hand side of \eqref{TruncEntrScal}, we arrive at the following local regularity criterion for any suitable weak solution $f$ of the Landau-Coulomb equation with initial data satisfying \eqref{ini-bnd-bis}. 
Assume that the scaled solution \eqref{f-eps} satisfies
$$
\esssup_{(t,v)\in Q_1(0,0)}f_\eps(t,\cdot)\star\frac{\un_{B_1(0)^c}}{|\,\cdot\,|}(v)\le Z_\eps\,,\qquad\text{ where }Z_\eps\ge 1\,.
$$
Then, there exists $\eta_{DG}\equiv\eta_{DG}[m_0,M_0,E_0,H_0,|v_0|]\in(0,1)$ such that
\begin{equation}\label{DG-0}
\begin{aligned}
\esssup_{-4<t\le 0}\int_{B_2(0)}(f_\eps(t,v)-1)_+\dv+\int_{Q_2(0,0)}|\nabla_v\sqrt{f_\eps}|^2\un_{f_\eps\ge 1}\dt \dv\le\tfrac{\eta_{DG}}{Z_\eps^{3/2}}\implies f_\eps\le 2\text{ a.e. on }Q_{1/2}(0,0)\,.
\end{aligned}
\end{equation}
See Lemma \ref{l:dg1} below.

The quantity in this local regularity criterion involves a sup in the time variable, and is therefore not an additive set function on parabolic cylinders. However, this regularity criterion can be improved as follows. For each $H$-solution $f$ of the 
Landau-Coulomb equation \eqref{e:landauNC}, set
$$
F(t,v,w):=\sqrt{f(t,v)f(t,w)/|v-w|}\Pi(v-w)(\nabla_v\ln f(t,v)-\nabla_w\ln f(t,w))\,.
$$
Under the same assumptions as above, and with $\eps_j=\lambda^j$ for some $\lambda\in(0,\tfrac14)$, there exists a threshold $\eta_{DG+}\equiv\eta_{DG+}[m_0,M_0,E_0,H_0,|v_0|]>0$ and $\lambda\in(0,\tfrac14)$ small enough such that 
\begin{equation}\label{DG+}
\begin{aligned}
\limsup_{j\to\infty}\eps_j^{-7/2}\int_{Q_{2\eps_j}(t_0,v_0)}\left(|\nabla_v\sqrt{f(t,v)}|^2\right.+&\left.\int_{\R^3}|F(t,v,w)|^2\dw\right)\dt \dv\le 2\eta_{DG+}&
\\
\implies&\limsup_{j\to\infty}\eps_j^{-3/2}\|f_{\eps_j}\|_{L^\infty_tL^1_v(Q_2(0,0))}\le\tfrac{\eta_{DG}}{3M_0^{3/2}}
\\
\implies&f\Big|_{Q_{\eps_{j_0}/2}(t_0,v_0)}\!\le 2\eps^{-2}_{j_0}\text{ a.e. for some }j_0>0\,.
\end{aligned}
\end{equation}
The two implications above are proved in Lemmas \ref{l:eventual-mass} and \ref{l:improved-dg} respectively. The first implication requires a control of the local mass by the dissipation rate. Setting $F_\eps(t,v,w):=\eps^{7/2}F(\eps^2t,\eps v,\eps w)$, 
for all $\lambda\in(0,\tfrac14)$
\begin{equation}\label{LocMass}
\begin{aligned}
\tfrac1{C_1\lambda^3}\|f_\eps\|_{L^\infty_tL^1_v(Q_{2\lambda}(0,0))}\le\left(1+\tfrac1{\lambda^4}\|F_\eps\|_{L^2(Q_2(0,0)\times\R^3)}\right)\|f_\eps\|_{L^\infty_tL^1_v(Q_2(0,0))}+\tfrac1{\lambda^8\eps}\|F_\eps\|^2_{L^2(Q_2(0,0)\times\R^3)}
\\
+\tfrac1{\lambda^8}\left(1+\|F_\eps\|^2_{L^2(Q_2(0,0)\times\R^3)}\right)\left\|\nabla_v\sqrt{f_\eps}\right\|^2_{L^2(Q_2(0,0))}&\,.
\end{aligned}
\end{equation}
The proof of this inequality relies on the weak formulation of the Landau-Coulomb equation with a test function involving the fundamental solution of the backward heat operator $\partial_t+\Delta_v$: see Lemma \ref{l:mass}.

\smallskip
Using the improved De Giorgi local regularity criterion \eqref{DG+} with a Vitali type covering lemma as in section 6 of \cite{ckn} leads to Theorem A: see Lemma \ref{l:vitali}

\smallskip
The proof of Theorem B is independent of the improved regularity criterion \eqref{DG+}. Returning to the truncated entropy inequality, we verify the first De Giorgi local regularity criterion \eqref{DG-0} by proving that the right-hand 
side of \eqref{TruncEntrScal} with $r_\eps=2$ and $\kap_\eps=\delta_\eps=\tfrac12$ vanishes as $\eps\to 0^+$. This is done in two steps:

\smallskip
\noindent
(a) assuming that $f$ is axisymmetric with axis of symmetry $\bar v+\R\omega$ (for $|\omega|=1$), prove that 
$$
\lim_{\eps\to 0^+}\int_{Q_{r_\eps+\delta_\eps}(0,0)}f_\eps(t,v)^2\dt \dv=0
$$
provided that the velocity $v_0$ in \eqref{f-eps} satisfies $\rho_0:=|(v_0-\bar v)\times\omega|>0$, and

\smallskip
\noindent
(b) prove that, for $0<\eps<\tfrac{\rho_0}6$
$$
\esssup_{(t,v)\in Q_3(0,0)}f_\eps(t,\cdot)\star\frac1{|\,\cdot\,|}(v)\le\frac{C_*}{\rho_0}\int_{\R^3}f(1+\ln_+f)(t,w)dw+C_*\rho_0^2\le Z_*[M_0,E_0,H_0,\rho_0]<\infty\,.
$$

\smallskip
The axial symmetry assumption on the suitable weak solution $f$ is an essential ingredient of the proofs of both statements (a) and (b), presented in section \ref{s:axi-sym} as Lemmas \ref{l:better-int} and \ref{l:lri} respectively. Both the proofs of 
Theorem A and Theorem B are based on the De Giorgi local regularity criterion \eqref{DG-0}

\subsection{Outline of this article}


In section~\ref{s:intro}, we collect preliminary results concerning the matrix diffusion coefficient and  the entropy dissipation bound, together with various useful estimates, obtained for example by interpolation. Section~\ref{s:approx} begins 
with the construction of suitable weak solutions. The heart of this section lies in the derivation of local entropy estimates (see \eqref{TruncEntrScal} above) tailored for the proof of the De Giorgi local regularity criterion in the next section.
It also contains the local mass estimates. The core of the present paper is section~\ref{s:degiorgi} where the statements and proofs of the De Giorgi local regularity criterion and of the eventual smallness of the local mass are to be found,
together with the proof of partial regularity (Theorem A). The regularity of axisymmetric solutions away from the axis of symmetry (Theorem B) belongs to Section~\ref{s:axi-sym}. We have gathered in the appendix a short proof of the fact 
that bounded solutions are smooth (Appendix \ref{a:smooth}), together with the lengthy and technical derivation of the truncated entropy inequalities (Appendix \ref{a:trunc-entropineq}).

\paragraph{Notation.} For all $r\in\R$, we set $r_+:=\max(r,0)$ and $r_-:=\max(-r,0)$. In particular, $r_\pm \ge 0$ and $r = r_+ - r_-$ while $|r|=r_+ + r_-$. For all $a,b\in\R$, we set $a \vee b:=\max(a,b)$ and $a \wedge b:=\min(a,b)$. 
The functions $h_+$ and $\lnp$ have been already introduced in \eqref{Def-h+}.

Given a set $A\subset\R^d$, the notation $\un_A$ designates the indicator function of the set $A$: in other words, $\un_A (\xi)= 1$ if $\xi \in A$, and $\un_A(\xi)=0$ otherwise. For simplicity, we set $\un_{z\ge\kappa}:=\un_{[\kappa,+\infty)}(z)$,
with similar definitions for $\un_{z>\kappa}$, or $\un_{z<\kappa}$ and $\un_{z\le\kappa}$.

We recall that $B_r(x)$ designates the open ball of the Euclidean space $\R^d$ with radius $r$ centered at $x\in\R^d$, and that $Q_r(t,x):=(t-r^2,t]\times B_r(x)$ is the parabolic cylinder of size $r$. To avoid unnecessarily heavy notations,
we shall sometimes replace $B_r(0)$ with $B_r$ and $Q_r(0,0)$ by $Q_r$. For a time interval $I$ and an open ball $\mB$ of $\R^3$, we set $L^q_t L^p_v (I \times \mB):=L^q(I,L^p (\mB))$. 

Weak convergence to $u$ of a sequence $\{u_n\}$ in the Banach space $E$ is denoted by $u_n\rightharpoonup u$ in $E$ as $n\to\infty$.

\paragraph{Constants.} A constant is said to be absolute if it is independent of all parameters and variables. The constants $c_\ast, c_{\ast \ast}, \bar c_1, \bar c_2, C_\ast \dots$ appearing below are absolute. Constants appearing in
functional inequalities such as $\Cinterpol$, $\Cnonlin\dots$ may depend on the spatial domains on which these functional inequalities hold.

The constants $C_0$, $C_1$, $C_2, \etadg, \etad, \etadgp$ appearing in the paper depend on the Euclidean norm $|v_0|$ of $v_0$ and on the bounds $m_0,M_0,E_0$ and $H_0$ on the macroscopic quantities attached to the initial data 
$\fin$. 

\paragraph{Acknowledgements.} The authors are indebted to A.V. Bobylev for fruitful discussions concerning the Landau equation in the radially symmetric case. We also thank M.P. Gualdani and N. Guillen for several discussions during 
the preparation of this paper.


\section{Diffusion, dissipation  and interpolation}
\label{s:intro}


This section contains preliminary results that will be used in the remainder of this work.

We have gathered together important properties of the diffusion matrix, starting with strict ellipticity. We continue with a local version of the entropy dissipation estimate due to L.~Desvillettes \cite{d15}. We conclude this section with various local 
estimates based on variants of the Sobolev embedding inequality.

\subsection{Derivatives of the matrix field {\it a}}

We take the collision kernel of the form 
\begin{equation}\label{Def-a}
a (z ) = \frac{\ka(|z|)}{|z|} \Pi (z)\,,\qquad z\not=0\,,
\end{equation}
where 
\begin{equation}\label{Ass-k}
\begin{aligned}
\ka\in C^1((0,+\infty))\text{ is positive, nondecreasing, such that }\ka(r)/r\text{ is nonincreasing in }(0,+\infty)&\,,
\\
\text{while }\ka'(|z|)/|z|^2\text{ belongs to }L^{1+\iota}(B_1)\text{ for some }\iota>0&\,.
\end{aligned}
\end{equation} 
Our main results correspond to $\ka (r) = r^{\gamma+3}$ with $\gamma \in [-3,-2)$. The following computations will be used repeatedly in the paper:
\begin{align}
\label{e:div-a}
(\nabla_z \cdot a) (z) &=   \ka (|z|) \nabla_z \left( \frac2{|z|}\right)\,,
\\
\label{e:hess-a}
(\nabla_z^2 \colon a) (z) &= -  \ka (0) 8 \pi \delta_0 - 2\frac{ \ka'(|z|)}{|z|^2}\,,
\end{align}
where $\delta_0$ denotes the Dirac mass at the origin in $\R^3$. We denote by $\sqrt{a(z)}$ the square root of the positive semidefinite real symmetric matrix $a(z)$, which is given by the formula
\begin{equation}\label{sqrta}
\sqrt{a(z)}:=\sqrt{\frac{k(|z|)}{|z|}}\Pi(z)\,,\qquad z\not=0\,.
\end{equation}

\subsection{Long- and short-range interactions}

The function $a$ is singular at the origin and needs to be mollified in order to construct solutions. Following \cite{ggiv}, we consider, for all $\delta \in (0,1)$,
\begin{equation}
\label{e:aout}
\aout (z) := X \left(\frac{|z|}\delta \right) a(z) = X (\delta^{-1}|z|)\frac{\ka (|z|)}{|z|} \Pi (z)\,,\qquad\text{ and }\quad\ain:=a-\aout\,,
\end{equation}
with $X \in C^\infty (\R)$ equal to $1$ in $[1,+\infty)$, supported in $[\tfrac12,+\infty)$ and such that $0 \le X' \le 3$. From \eqref{e:div-a} with $\ka(r)$ replaced with $\ka(r) X(\delta^{-1}r)$, we get
\begin{align*}
\nabla_z \cdot \aout (z) &= \ka(|z|) X (\delta^{-1} |z|) \nabla_z (2/ |z|)\,,
\\
\nabla_z^2 : \aout (z) &= - 2\bigg(\ka' (|z|)X(\delta^{-1}|z|) + \delta^{-1} \ka (|z|) X'(\delta^{-1}|z|) \bigg) |z|^{-2}\,.
\end{align*}
In particular, for $|z| \ge \delta/2$, we have
\begin{equation}
\label{e:div-aout}
|\nabla_z \cdot \aout (z) | \le 8 \frac{ \ka (\delta/2)}{\delta^2} \quad \text{ and } \quad |\nabla_z^2 : \aout (z)| \le 40 \frac{\ka (\delta/2)}{\delta^3} . 
\end{equation}
Indeed,  $\tfrac{\ka'(|z|)}{|z|}\le\tfrac{\ka (|z|)}{|z|^2}$ since $r\mapsto\tfrac{k(r)}r$ is nonincreasing and, for $|z|\ge \delta/2$, we have $\tfrac{\ka (|z|)}{|z|} \le2\tfrac{\ka (\delta/2)}{\delta}$. 

\subsection{Diffusion from long-range interactions}

We recall here \cite[Lemma~3.3]{ggiv} which provides us with a lower bound on the diffusion matrix corresponding to long-range interactions,
\[ 
\Aout (v) := \int_{\R^3} \aout (v-w) f(w) \dw .
\]
Since the statement used here is slightly different from that of \cite{ggiv}, we include a sketch of the proof for the reader's convenience.

\begin{lemma}[Diffusion matrix]\label{l:diff-mat}
Pick $M_0>m_0>0$ and $E_0,H_0>0$, and let $f \colon \R^3 \to[0,+\infty)$ be a measurable function such that
\begin{equation}\label{e:hydro-bounds}
\int_{\R^3} f(v) \dv \in [m_0,M_0],  \quad \int_{\R^3} f(v) |v|^2 \dv \le E_0 \quad \text{ and } \quad  \int_{\R^3} f \ln f (v) \dv \le H_0.
\end{equation}
There exists $\bar\delta \in (0,1)$ and $c_0,R_0 >0$ depending only on  $m_0,M_0,E_0,H_0$ such that for all $\delta \in (0,\bar\delta]$, we have
\[ 
 A^\out (v) \ge c_0 \frac{ \ka (|v|+R_0)}{(1+|v|)^3} I\,,
 \]
 where $I$ is the $3 \times 3$ identity matrix. 
\end{lemma}

\begin{proof}
Let $e \in \R^3$ be such that $|e|=1$. We know from \cite[Lemma~3.3]{zbMATH06668299} that there exists $R_0,\ell, \mu>0$ depending only on $m_0,M_0,E_0, H_0$ such that
\[
|\{ z \in B_{R_0} (-v) : f(v+z) \ge \ell \}| \ge \mu \quad \text{ and } | \{ z \in B_{R_0} (-v) : (e \cdot z/|z|)^2 \ge 1 -\eps_0 \}| \le \mu /2 \,,
\]
with $\eps_0  = \mu (\mu + 4 \pi R_0 (|v|+R_0)^2)^{-1}$ (see \cite[Lemma~3.3]{ggiv} for a justification of this explicit expression of $\eps_0$).

Recalling the definition of $A^\out$ and $\aout$, see \eqref{e:aout}, we write next
\begin{align*}
A^\out (v) e \cdot e & = \int_{\R^3} (\aout (z) e \cdot e) f(v+z) \dz = \int_{\R^3} \frac{\ka (|z|) X (|z|/\delta)}{|z|} (1 - (e \cdot z/|z|)^2)  f(v+z) \dz .
\intertext{Letting $\mathcal{A} = \{ z \in B_{R_0} (-v) : f(v+z) \ge \ell \} \setminus  \{ z \in B_{R_0} (-v) : (e \cdot z/|z|)^2 \ge 1 -\eps_0 \}$ and using the fact that $X=1$ in $[1,+\infty)$ while $\ka$ is nondecreasing, we get,}
A^\out (v) e \cdot e  & \ge  \eps_0 \ell\frac{\ka (|v|+R_0)}{ (|v|+R_0)} | \mathcal{A} \setminus B_\delta |. 
\\
\intertext{Recalling the expression for $\eps_0$ given above, we see that}
A^\out (v) e \cdot e     & \ge \frac{\ell \mu \ka (|v|+R_0)}{(\mu + 4 \pi R (|v|+R_0)^2)(|v|+R_0)} ( \mu /2 - 4 \pi \delta^3/3).
\end{align*}
We conclude by choosing $\bar\delta$ such that $4 \pi \bar\delta^3/3 \le \mu /3$. 
\end{proof}

\subsection{Local entropy dissipation estimate}

This subsection is devoted to the derivation of a bound from below for the localized truncated entropy dissipation term that will appear in local entropy estimates, see next section. We adapt here the reasoning from our previous work 
\cite{ggiv} in collaboration with M.P. Gualdani. It is simpler than the proof from \cite{d15} and was suggested to us by an anonymous referee.

\begin{prop}[Local entropy dissipation estimate]\label{p:lede}
Let $F,f,\Psi: \R^3 \to [0,+\infty)$ be such that $\Psi \in C^1$. Then for $\delta \in (0,\bar\delta)$ with $\bar\delta$ given by Lemma~\ref{l:diff-mat},
\begin{multline*}
\tfrac12  \iint_{\R^3 \times \R^3}f(v)f(w) \,  \aout (v-w) :  \left( \Psi (v) \frac{\nabla_v F(v)}{f(v)} - \Psi (w) \frac{\nabla_w F (w)}{f(w)} \right)^{\otimes 2} \dv \dw 
\\
\ge   c_0 \int_{\R^3}\frac{|\nabla_v F(v)|^2}{f(v)} \frac{k (|v|+R_0)}{(1+|v|)^3} \Psi^2 (v) \dv - 40 \frac{\ka (\delta)}{\delta^3} \left\{ \int_{\R^3} F(v) (\Psi (v) + \delta |\nabla_v \Psi (v)|) \dv \right\}^2\,,
\end{multline*}
where $c_0$ is given by Lemma~\ref{l:diff-mat}.
\end{prop}

\begin{proof}
We start from the left side of the desired estimate and write,
\begin{align}
\nonumber  
\tfrac12 \iint_{\R^3 \times \R^3}f(v)f(w) \aout (v-w) : &\left( \Psi (v) \frac{\nabla_v F(v)}{f(v)} - \Psi (w) \frac{\nabla_w F (w)}{f(w)} \right)^{\otimes 2} \dv \dw 
\\
\nonumber  
=& \iint_{\R^3 \times \R^3}f(v)f(w) \, \aout (v-w) : \left( \Psi (v) \frac{\nabla_v F(v)}{f(v)}  \right)^{\otimes 2} \dv \dw 
\\
\nonumber  & + \iint_{\R^3 \times \R^3} \aout (v-w)  \bigg( \Psi (v) \nabla_v F(v) \bigg) \cdot  \bigg( \Psi (w) \nabla_w F(w) \bigg) \dv \dw
\\
\intertext{integrating by parts twice in the second term}
\nonumber  =& \int_{\R^3} \Aout (t,v) : \left( \Psi (v) \frac{\nabla_v F(v)}{\sqrt{f(v)}}  \right)^{\otimes 2} \dv  
\\
\label{e:diff-1} & + \iint_{\R^3 \times \R^3}  \nabla_v \cdot \bigg( \nabla_w \cdot \bigg( \aout (v-w) \Psi (v) \Psi (w) \bigg) \bigg)F(v) F(w)\dv \dw
\end{align}
with $\Aout=\aout \ast f$. As far as the second term in the right-hand side of \eqref{e:diff-1} is concerned, we first compute
\begin{align*}
\nabla_v \cdot \bigg( \nabla_w \cdot \bigg( \aout (v-w) \Psi (v) \Psi (w) \bigg) \bigg)\!\!= & (-\nabla_z^2  : \aout) (v-w) \Psi (v) \Psi (w)   + \left(\aout (v-w) \nabla_v  \Psi (v)\right) \cdot \nabla_w \Psi (w) 
\\
& + (\nabla_z \cdot \aout)(v-w)\cdot(\Psi(v) \nabla_w \Psi (w) - \Psi(w) \nabla_v \Psi (v)).
\end{align*}
Use now \eqref{e:div-aout}  and the fact that $\aout (z)=0$ for $|z| \le \delta/2$ to get
\begin{align*}
\left|  \nabla_v \cdot \bigg( \nabla_w \cdot \bigg( \aout (v-w) \Psi (v) \Psi (w) \bigg) \bigg) \right| \le& \frac{40 \ka (\delta/2)}{\delta^{3}} \Psi (v) \Psi (w)   +\frac{2\ka (\delta/2)}{\delta} |\nabla_v  \Psi (v)|| \nabla_w \Psi (w)| 
\\
& + \frac{8\ka(\delta/2)}{ \delta^2} ( \Psi(v) |\nabla_w \Psi (w)| + \Psi(w) |\nabla_v \Psi (v)|) 
\\
\le& 40 \frac{\ka (\delta)}{\delta^3} ( \Psi(v) + \delta |\nabla_v \Psi (v)|)\; ( \Psi(w) + \delta |\nabla_w \Psi (w)|).
\end{align*}
Combining \eqref{e:diff-1} with Lemma~\ref{l:diff-mat} and the previous estimate yields the desired conclusion.
\end{proof}

\subsection{Interpolation, short-range interactions, and nonlinearization}

Throughout this section, $I \subset[0,+\infty)$ denotes a time interval and $\mB \subset \R^3$ an open ball.

We recall that $H$-solutions of the Landau equation satisfy the conservation laws of mass, energy and the inequality implied by Boltzmann's H-Theorem (see Definition 1 (c), formula (39) and Definition 2 in section 3 of \cite{V98}). 
Any $H$-solution $f$ of the Landau equation satisfies therefore $ f\in L^\infty_t L^1_v (I \times \mB)$, and $\nabla_v\sqrt{f}\in L^2(I \times \mB)$ by Theorem 3 of \cite{d15}. The next lemma makes precise the bounds in Lebesgue spaces 
that can be obtained by interpolation from these two pieces of information. 

\begin{lemma}[Interpolation]\label{l:interpolation}
Given $p, q \ge 1$ such that $\frac1p + \frac2{3q} \ge 1$, assume that $f \in L^\infty_t L^1_v (I \times \mB)$ is such that $\sqrt{f} \in  L^2(I \times \mB)$. Then
\[
\|f\|_{L^q_t L^p_v(I \times \mB)} \le \Cinterpol \left(\|f\|_{L^\infty_t L^1_v (I \times \mB)} + \left\|\nabla_v \sqrt{f}\right\|_{L^2 (I \times \mB)}^2\right)\,,
\]
where $\Cinterpol$ depends only on an upper bound of the radius of $\mB$ and the length of $I$, and on $p,q$.
\end{lemma}

\begin{proof}
By H\"older's inequality, it is enough to treat the case $\frac1p + \frac2{3q}=1$ since $I$ and $\mB$ are of finite measure. Observe that $p \in [1,3]$, and that
 \[
 \frac1p = \frac{1}1 \left( 1 - \frac1q \right)+ \frac13 \left( \frac1q \right) . 
 \]
Therefore, by H\"older's inequality,
\[
\|f\|_{L^q_t L^{p}_v (I \times \mB)}^q = \int_I \|f(t,\cdot)\|_{L^{p}_v (\mB)}^q \dt \le \int_I \|f(t,\cdot)\|_{L^1_v (\mB)}^{q-1} \|f(t,\cdot)\|_{L^{3}_v (\mB)} \dt
\]
so that, using again H\"older's inequality followed by Young's inequality,
\begin{align*}
\|f\|_{L^q_t L^{p}_v (I \times \mB)} & \le \|f\|_{L^\infty_t L^1_v (I \times \mB)}^{1-\frac1q} \|f\|_{L^1_t L^{3}_v (I \times \mB)}^{\frac1q} \le (1-\tfrac1q)\|f\|_{L^\infty_t L^1_v (I \times \mB)}+ \tfrac1q \|f\|_{L^1_t L^{3}_v (I \times \mB)}\,.
\end{align*}
Next we use the Sobolev inequality for $\sqrt{f}$:
\begin{align*}
\|f\|_{L^1_t L^{3}_v(I\times \mB)}^{\frac12} &= \|\sqrt{f}\|_{L^2_t L^{6}_v(I \times \mB)} \le \Csob ( \|\sqrt{f}\|_{L^2(I \times \mB)} + \|\nabla_v \sqrt{f}\|_{L^2 (I \times \mB)}) 
\\
& = \Csob ( \|f\|_{L^1(I \times \mB)}^{\frac12} + \|\nabla_v \sqrt{f}\|_{L^2 (I \times \mB)}) 
\\
& \le \Csob ( |I|^{1/2}\|f\|_{L^\infty_t L^1_v (I \times \mB)}^{\frac12} + \|\nabla_v \sqrt{f}\|_{L^2 (I \times \mB)}).
\end{align*}
Combining the two inequalities leads to the announced estimate with $\Cinterpol=(1+2\Csob^2)(|I|+1)$.
\end{proof}

\begin{lemma}[Control of short-range interactions]\label{l:source-nl}
Let $f\,,\,\, F:\,I\times\mathbb R^3\to[0,+\infty)$ be measurable functions defined a.e., and let $\rho\in(0,3)$. Then 
\[
\int_I\int_\mB\left(f(t,\cdot)\ast\frac{\un_{B_\delta}}{|\cdot|^{\rho}}(v)\right)F(t,v)\dv\dt\le C_{HLS,\rho}\|f\|_{L^q_tL^p_v(I\times(\mB+B_\delta))}\|F\|_{L^\theta_tL^\beta_v(I\times\mB)}
\]
with $1<p,q,\beta,\theta<\infty$ provided that
\[
\frac1q+\frac1\theta=1\quad\text{ and }\quad\frac1p+\frac1\beta+\frac{\rho}3=2\,,
\]
and where $C_{HLS,\rho}$ is an absolute constant for each $\rho\in(0,3)$. Besides, if $\rho<2$, one has
\[
\left\|f\ast_v\frac{\un_{B_\delta}}{|\cdot|^{\rho}}\right\|_{L^1_tL^\infty_v(I\times\mB)}\le\left(\frac{8\pi}{3(2-\rho)}\right)^{2/3}\delta^{2-\rho}\|f\|_{L^1_tL^3_v(I\times(\mB+B_\delta))}\,.
\]
\end{lemma}

\begin{proof}
The first inequality follows from the Hardy-Littlewood-Sobolev inequality (see Theorem 4.3 in \cite{LiebLoss}) in the $v$-variable (with $f$ replaced with $f\un_{\mB+B_\delta}$ and $F$ replaced with $F\un_\mB$), and from the H\"older inequality in 
the $t$-variable. The second inequality follows from the H\"older inequality applied to $f\un_{\mB+B_\delta}$ after observing that
$$
\left\|\frac{\un_{B_\delta}}{|\cdot|^{\rho}}\right\|^{3/2}_{L^1_tL^{3/2}_v(I\times\mB)}=4\pi\int_0^\delta r^{2-\frac32\rho}\dd r=\tfrac{8\pi}{3(2-\rho)}\delta^{3(1-\frac{\rho}2)}\,.
$$
\end{proof}

\begin{lemma}[Nonlinearization]\label{l:nonlin}
Given $p, q \ge 1$ such that  $\frac{1}p + \frac2{3q}=1$ and  $\kap,\kbar$ such that  $1\le \kap < \kbar < \kap+1$, for each $g:\,I\times\mathcal B\to[0,+\infty)$ we have
\[
  \| (g-\kbar)_+ \|_{L^q_t L^p_v(I \times \mB)}
  \le\Cnonlin \frac{\kbar}{(\sqrt{\kbar}-\sqrt{\kap})^{4}}\left( \| (g-\kap)_+\|_{L^\infty_t L^1_v (I \times \mathcal{B})} + \left\|\nabla_v(\sqrt{g}-\sqrt{\kap})_+ \right\|_{L^2 (I \times  \mathcal{B})}^2 \right)
\]
 where $\Cnonlin$ depends only on an upper bound for the radius of $\mathcal{B}$. 
\end{lemma}

\begin{proof}
We split the proof in two steps. 
  
\paragraph{Step~1.}
Consider $\Gkp := \Gamma ((\sqrt{g}-\sqrt{\kap})_+)$ with  $\Gamma (r) = \min (r,r^2)$. Observe first that
\begin{equation}\label{UpBndVN}
(\Gkp)^2=\min\left(\left(\sqrt{g}-\sqrt{\kap}\right)_+,\left(\sqrt{g}-\sqrt{\kap}\right)^2_+\right)^2\le\left(\sqrt{g}-\sqrt{\kap}\right)^2_+\le\left(\sqrt{g}-\sqrt{\kap}\right)_+\left(\sqrt{g}+\sqrt{\kap}\right)=(g-\kappa)_+\,.
\end{equation}
On the other hand
\begin{align*}
(g-\kbar)_+   &\le (\sqrt{g}-\sqrt{\kbar}+\sqrt{\kbar})^2_+ \un_{\{ g \ge \kbar\}} \le 2 (\sqrt{g}-\sqrt{\kbar})_+^2 + 2 \kbar \un_{\{g \ge \kbar\}}
\\
\intertext{use here that $y \un_{\{g \ge \kbar\}} = \Gamma (y \un_{\{g \ge \kbar\}})$ if $y \ge 1$ and $y \un_{\{g \ge \kbar\}} \le \un_{\{g \ge \kbar\}}$ if $y \le 1$,}
&\le 2 \left(\Gamma \left(\left(\sqrt{g}-\sqrt{\kbar}\right)_+\right)\right)^2 + 2( \kbar+1) \; \un_{\{g \ge \kbar\}} = \bigg( 2 (\Gkp)^2 +  2(\kbar+1) \bigg) \un_{\{ \sqrt{g}- \sqrt{\kap} \ge \sqrt{\kbar}-\sqrt{\kap} \}} 
\\
\intertext{use the fact that $\Gamma$ is nondecreasing and $0 < \sqrt{\kbar}-\sqrt{\kap}=\frac{\kbar-\kap}{\sqrt{\kbar}+\sqrt{\kap}} < 1$ since $0<\kbar-\kap<1$ and $\kbar>1$,}
& \le 2 \bigg( (\Gkp)^2 +  (\kbar+1)\bigg) \un_{\{\Gkp \ge (\sqrt{\kbar}-\sqrt{\kap})^2 \}} \le 2 \left[ 1 + (\kbar+1)(\sqrt{\kbar}-\sqrt{\kap})^{-4} \right] (\Gkp)^2 .
\end{align*}
Using again that $0 < \sqrt{\kbar}-\sqrt{\kap} < 1$ and $\kbar > 1$, we conclude that
\begin{equation}\label{e:fF}  
(g-\kbar)_+ \le \frac{6\kbar}{ (\sqrt{\kbar}-\sqrt{\kap})^{4}} (\Gkp)^2 .
\end{equation}

\paragraph{Step~2.} Lemma~\ref{l:interpolation} implies that
\begin{align*}
\|(g-\kbar)_+\|_{L^q_t L^p_v(I \times \mB)} & \le \frac{6\kbar}{(\sqrt{\kbar}-\sqrt{\kap})^{4}} \| (\Gkp)^2\|_{L^q_t L^p_v(I \times \mB)} 
\\
& \le \Cinterpol  \frac{6\kbar}{(\sqrt{\kbar}-\sqrt{\kap})^{4}}  (\|(\Gkp)^2\|_{L^\infty_t L^1_v (I \times \mB)} + \|\nabla \Gkp\|_{L^2 (I \times \mB)}^2) 
\\
& \le \Cinterpol  \frac{6\kbar}{(\sqrt{\kbar}-\sqrt{\kap})^{4}}  (\|\Gkp\|^2_{L^\infty_t L^2_v (I \times \mB)} + 4 \|\nabla (\sqrt{g}-\sqrt{\kap})_+\|_{L^2 (I \times \mB)}^2)\,,
\end{align*}
where we used that $\Gamma$ is a $2$-Lipschitz function in $[0,+\infty)$ to get the last line. With \eqref{UpBndVN}, this implies the desired inequality.
\end{proof}


\section{Suitable solutions: local entropy and mass estimates}
\label{s:approx}


For the reader's convenience, we first recall precisely Villani's notion of $H$-solution.

\begin{defn}[H-solutions \cite{V98}]
Given $T>0$ and $\fin \colon \R^3 \to [0,+\infty)$ measurable satisfying \eqref{ini-bnd}, a function $f: [0,T) \times \R^3 \to [0,+\infty)$ is said to be a \emph{H-solution} of \eqref{e:landau} with initial condition \eqref{e:ic} if the function $f$ belongs to
$C([0,T); \mathcal{D}'(\R^3)) \cap L^1((0,T);L^1(\R^3,\tfrac{\dv}{1+|v|}))$, if in addition
\[
(t,v,w)\!\mapsto\! F(t,v,w)\! :=\!  \sqrt{a (v\!-\!w)}\sqrt{f(t,v) f (t,w) }  \left( \nabla_v \ln f (t,v) \!-\! \nabla_w \ln f (t,w) \right)   \text{ belongs to } L^2([0,T)\times \R^3 \times \R^3)
\]
with $\sqrt{a(z)}$ defined in \eqref{sqrta} and if, for a.e. $t \in [0,T)$,
\begin{align*}
\int_{\R^3} \left(\begin{array}{c} 1 \\ v \\ |v|^2 \end{array}\right) f(t,v) \dv &=\int_{\R^3}\left(\begin{array}{c} 1 \\ v \\ |v|^2 \end{array}\right) \fin(v) \dv\,,
\\
\int_{\R^3} (f \ln f) (t,v) \dv &\le \int_{\R^3} (\fin \ln \fin) (v) \dv\,, 
\end{align*}
while, for all $\phi \in C^1_c ([0,T) \times \R^3)$,
\begin{multline*}
\int_{\R^3} \fin (v) \phi (0,v) \dv + \int_0^T \int_{\R^3} f (t,v) \partial_t \phi (t,v) \dt \dv
\\
=\tfrac12\int_0^T \iint_{\R^3 \times \R^3} F(t,v,w) \cdot \sqrt{a(v-w)}\sqrt{f(t,v)f(t,w)} (\nabla_v \phi (t,v) -\nabla_w \phi (t,w)) \dt \dv \dw\,.
\end{multline*}
\end{defn}

We do not know whether our main results, Theorems A and B, apply to all $H$-solutions of \eqref{e:landau}. However, they apply to all $H$-solutions {\it constructed as in Villani's global existence proof} of $H$-solutions for the Landau equation in 
\cite{V98}. 

We therefore introduce a notion of \emph{suitable weak solution} of the Landau equation, which are $H$-solutions satisfying the extra properties used in the proofs of Theorem A and B, and of Corollary C.

\begin{defn}[Suitable weak solutions]\label{d:suitable}
Given $m_0,M_0,E_0,H_0,T>0$ and $\fin\colon \R^3 \to [0,+\infty)$ measurable satisfying \eqref{ini-bnd-bis}, a H-solution $f$ of the Landau equation in $[0,T)\times\R^3$ with initial condition \eqref{e:ic} is said to be a \emph{suitable weak solution} if
\begin{enumerate}
\item
there exist a constant $\cI_0:=\cI_0[m_0,M_0,E_0, H_0]>0$ and $T>0$, such that, for all $\kappa\in[1,+\infty)\cap\mathbb Q$,
\[ 
\left\|F\right\|_{L^2 ((0,T) \times \R^3 \times \R^3)} + \left\|F^\kappa_+\right\|_{L^2 ((0,T) \times \R^3 \times \R^3)} \le  \cI_0   
\]
where
\(
F^\kappa_+ (t,v,w) = \sqrt{a(v-w)} \sqrt{f(t,v) f(t,w)}\left(\nabla_v \ln_+ (f/\kappa) (t,v) - \nabla_w \ln_+ (f/\kappa) (t,w) \right);
\)
\item  
for each $\Psi\in C^\infty_c((0,T) \times \R^3)$ and $\kappa\in[1,+\infty)\cap\mathbb Q$, there exists a negligible set $\mathcal{N}\subset (0,T)$ such that 
\begin{multline}\label{e:lee}
\int_{\R^3}    h_+^\kappa(f (t_2,v))  \varphi (t,v) \dv \le  \int_{\R^3}  h_+^\kappa (f (t_1,v))  \varphi (t,v)  \dv + \iint_{[t_1,t_2] \times \R^3} h_+^\kappa(f)(t,v) \partial_t \varphi(t,v)  \dt \dv
\\
- \iint_{[t_1,t_2] \times \R^3} ( A\nabla_v f(t,v) -f(\nabla_v \cdot A)(t,v)) \cdot \bigg( \varphi\frac{\nabla_v (f-\kappa)_+}{f}(t,v) + \lnp \left(\frac{f}\kappa\right)\nabla_v \varphi(t,v) \bigg) \dt \dv 
\end{multline}
for all $t_1,t_2 \in (0,T) \setminus \mathcal{N}$ with  $t_1 <t_2$, where $\varphi := \Psi^2$ while
\begin{equation}\label{Def-hkappa}
h_+^\kappa (r) := \kappa h_+ (r/\kappa) = r \ln_+ (r/\kappa) -(r-\kappa)_+\,,
\end{equation} 
and 
\begin{equation}\label{Def-A}
A(t,v) := \int_{\R^3} a(v-w) f(t,w) \dw\,.
\end{equation}
\end{enumerate}
\end{defn}

\smallskip
We first state a global existence result of suitable weak solutions of the Landau equation.

\begin{prop}[Existence of suitable weak solutions]\label{p:existence}
Let $\ka$ satisfy \eqref{Ass-k}, and let $\fin \colon \R^3 \to [0,+\infty)$ measurable and satisfying \eqref{ini-bnd}. Then for all $T>0$ there exists a suitable weak solution $f$ of the Landau equation \eqref{e:landau} in $[0,T) \times \R^3$ with 
initial condition \eqref{e:ic} and collision kernel \eqref{Def-a}.
\end{prop}

This section is devoted to the proof of Proposition \ref{p:existence} and to the derivation of local entropy and local mass estimates. 

\subsection{The approximation scheme}\label{ss:approx-scheme}

The main result of this article is about the H-solutions  constructed by the approximation scheme from \cite{V98}. It is explained in \cite[p.~297]{V98} that there exist global smooth solutions of the following approximate equation,
\begin{equation}\label{e:landau-reg}
\partial_t f^n (t,v) =  \nabla_v\cdot \left( \int_{\R^3} a_n (v-w) (f^n(t,w) \nabla_vf^n(t,v) - f^n(t,v) \nabla_wf^n(t,w) ) \dw \right) + \frac1n \Delta_v f^n (t,v)\,,
\end{equation}
where
\begin{equation}\label{e:def-an}
a_n (z)  := X(n|z|) a(z) = X (n|z|) \frac{\ka (|z|)}{|z|} \Pi (z)\,,
\end{equation}
with $X \in C^\infty (\R)$ equal to $1$  in $[1,+\infty)$, supported in $[\tfrac12,+\infty)$ and such that $0 \le X' \le 3$.

The approximate initial data $\fin^n$ is given by
\begin{equation}\label{e:finn}
\fin^n (v):= \tilde{\fin^n} (v) + \frac1n \frac{e^{-|v|^2/2}}{(2\pi)^{3/2}} \quad \text{ with } \quad \tilde{\fin^n} := \zeta_n \ast (\xi_n \fin)
\end{equation}
where $\xi_n$ is a truncation and $\zeta_n$ is a mollifier. More precisely, we consider a function $\xi\in C^\infty (\R^3,[0,1])$ such that $\un_{B_1}\le\xi\le\un_{B_2}$ and, for each integer $n\ge1$, we consider the cutoff function
$\xi_n (v) := \xi (v/n)$ and the mollifier $\zeta_n (v) := n^3 \xi (nv) / \|\xi\|_{L^1 (\R^3)}$.

We recall the notation $\ain$ and $\aout$ defined in \eqref{e:aout}, and we set
\begin{equation}\label{e:anin}
\anin(z):=X(n|z|)\ain(z)\,,\qquad\anout(z):=a_n(z)-\anin(z)=X(n|z|)X\left(\frac{|z|}\delta\right)a(z)\,.
\end{equation}

\begin{prop}[Approximate solutions]\label{p:approximate}
Let $\fin^n$ be given by \eqref{e:finn} for a measurable function $\fin\ge 0$ a.e. satisfying \eqref{ini-bnd-bis}. There exists a $H$-solution  $f$ of \eqref{e:landau} with initial data $\fin$ and a subsequence  $\{ f^{n_i}\}_i$ of the smooth 
global nonnegative solutions of  \eqref{e:landau-reg} with initial data $\fin^n$ such that, for all $T>0$ and   $R>0$, there exists a constant $\cI (T,R)>0$ depending only on $T,R$ and the hydrodynamic bounds $m_0,M_0,E_0$ and $H_0$ 
and $\cI_0>0$ depending only on  $m_0,M_0,E_0$ and $H_0$ such that
\begin{align}
\label{e:gradfn} 
\left\|\nabla_v \sqrt{f^n }\right\|_{L^2 ((0,T) \times B_R)}\le \cI (T,R)\,,
\\
\label{e:fn-equi}
\|f^n\|_{L^{5/3} ((0,T) \times B_R)} \le \cI (T,R)\,,
\\
\label{e:fn-mass}
\int_{\R^3} f^n (t,v) \dv  = \int_{\R^3} \fin^n (v) \dv\,,
\\
\label{e:dissip-all}
\left\|F^n\right\|_{L^2 ((0,T) \times \R^3 \times \R^3)} \le \cI_0 \,, 
\\
\nonumber
\left\|F^{ \kappa, n}_+\right\|_{L^2 ((0,T) \times \R^3 \times \R^3)} \le \cI\,,
\end{align}
and
\begin{align}
\label{e:ae-conv} 
f^{n_i} \to f  \text{ a.e. } & \text{ in } (0,+\infty) \times \R^3\,,
\\
\label{e:gradfn-conv}
\nabla_v \sqrt{f^{n_i}} \rightharpoonup \nabla_v \sqrt{f} & \text{ in } L^2((0,T) \times B_R)\,,
\\
\label{e:dissip-conv}
F^{n_i} \rightharpoonup F &\text{ in } L^2 ((0,T)\times \R^3 \times \R^3)\,,
\\
\label{e:dissipk-conv}
F^{\kappa, \inn,  n_i}_+ \rightharpoonup F^{\kappa, \inn}_+ & \text{ in } L^2 ((0,T)\times \R^3 \times \R^3)\,,
\end{align}
as $n_i\to\infty$. In the limits above,
\begin{equation}\label{e:Fn}
\left\{
\begin{aligned}
 F^n (t,v,w) &:=  \sqrt{a_n (v-w)}\sqrt{ f^n(t,v) f^n (t,w) }  \left( \frac{\nabla_v f^n (t,v)}{f^n(t,v)} - \frac{\nabla_w f^n (t,w)}{f^n(t,w)} \right)\,, 
 \\
 F (t,v,w) &:=  \sqrt{a (v-w)}\sqrt{ f(t,v) f (t,w) }  \left( \frac{\nabla_v f (t,v)}{f(t,v)} - \frac{\nabla_w f (t,w)}{f(t,w)} \right)\,, 
\end{aligned} \right.
\end{equation}
and
\begin{equation}\label{e:fkappain}
\left\{
\begin{aligned}
F^{\kappa, n}_+ (t,v,w) &:=  \sqrt{a (v-w)}\sqrt{ f^n(t,v) f^n (t,w) }\left( \frac{\nabla_v \fnkp (t,v)}{f^n(t,v)} - \frac{\nabla_w \fnkp (t,w)}{f^n(t,w)} \right)\,, 
\\
F^{\kappa,\inn, n}_+ (t,v,w) &:=  \sqrt{\anin (v-w)}\sqrt{ f^n(t,v) f^n (t,w) }\left( \frac{\nabla_v \fnkp (t,v)}{f^n(t,v)} - \frac{\nabla_w \fnkp (t,w)}{f^n(t,w)} \right)\,, 
\\
F^{\kappa,\inn}_+ (t,v,w) &:=  \sqrt{\ain (v-w)}\sqrt{ f(t,v) f (t,w) }  \left( \frac{\nabla_v \fkp (t,v)}{f(t,v)} - \frac{\nabla_w \fkp (t,w)}{f(t,w)} \right)\,,
\end{aligned}\right.
\end{equation}
with the notation
\begin{equation}\label{f-ka+}
\fkp(t,v):=(f(t,v)-\kappa)_+\,,\qquad\fnkp(t,v):=(f^n(t,v)-\kappa)_+\,.
\end{equation}
\end{prop}

\begin{proof}
We first pick $\ned \in \N$ depending only on $m_0,M_0,E_0,H_0$ such that for all $n \ge \ned$, we have
\begin{equation}
\label{e:hydro}
\int_{\R^3} \fin^n (v) \dv \in [m_0/2,2M_0], \quad \int_{\R^3} \fin^n (v) |v|^2 \dv \le 2 E_0, \quad \int_{\R^3} \fin^n \ln \fin^n (v) \dv \le 2 H_0\,.
\end{equation}

\paragraph{Entropy estimates.} Moreover, as explained in equations (17)-(18) of \cite{ggiv}, the conservation of mass/energy and Boltzmann's H-Theorem imply  for all $t>0$ and all $\kappa \ge 1$,
\begin{align*}
\int_{\R^3} f^n (t,v) (1+|v|^2) \dv  = \int_{\R^3} \fin^n (v) (1+|v|^2)\dv\,, 
\\
\int_{\R^3} f^n \ln f^n (t,v) \dv + \tfrac12 \int_0^t \iint_{\R^3 \times \R^3} |F^n (s,v,w)|^2 \ds \dv \dw \le \int_{\R^3} \fin^n  \ln \fin^n  (v) \dv\,, 
\\
\int_{\R^3} h^\kappa_+ (f^n) (t,v) \dv + \tfrac12 \int_0^t \iint_{\R^3 \times \R^3} \left|F^{n,\kappa}_+ (s,v,w)\right|^2 \ds \dv \dw \le \int_{\R^3} h^\kappa_+ (\fin^n)  (v) \dv\,.
\end{align*}
In particular, \eqref{e:dissip-all} holds true. 

It is also explained in \cite{ggiv} (see the discussion following equation (18) in \cite{ggiv}) that \( f \ln_+ f \le f \ln (f/G) - f + G\) where $G (v)= \frac1{(2 \pi)^{3/2}}e^{-|v|^2/2}$, and in particular,
\[ 
f \ln_+ f \le f \ln f + \tfrac12(|v|^2+3\ln(2\pi)-2)f + G\,.
\] 
This inequality together with the previous estimates yield the following ones: for all $t>0$,
\begin{align}
\label{e:second-m} 
&\int_{\R^3} f^n (t,v) (1+|v|^2) \dv  \le 2(M_0+E_0)\,,
\\
\label{e:dpsin} 
&\int_{\R^3} f^n \ln_+ f^n (t,v) \dv + \tfrac12 \int_0^t \iint_{\R^3 \times \R^3} |F^n (s,v,w)|^2 \ds \dv \dw \le \bar{H}_0:= H_0 + E_0+3\ln(2\pi)M_0+1\,.  
\end{align}
Applying the Dunford-Pettis theorem (Theorem 4.30 in \cite{brezis}) shows that we can extract a subsequence $f^{n_i}$ such that
\( 
f^{n_i} \rightharpoonup f \text{ in } L^1 ((0,T) \times B_R)
\) 
for some $H$-solution $f$ of \eqref{e:landau} with initial data $\fin$.

\paragraph{Integrability.} Proposition~\ref{p:lede} with $\delta = \bar\delta$ (given by Lemma~\ref{l:diff-mat}), with $F=f=f^n$ and (a sequence of test functions converging to) $\Psi =1$ implies that for all $t>0$,
\[
4c_0 \int_{\R^3} \left|\nabla_v \sqrt{f^n (t,v)}\right|^2 \frac{\ka (|v|+R_0)}{(1+|v|)^3} \dv \le \tfrac12 \iint_{\R^3 \times \R^3} |F^n(t,v,w)|^2 \dv \dw  + \frac{40 k (\bar\delta)}{\bar\delta^3} (2M_0)^2\,,
\]
where $c_0$ is also given by Lemma~\ref{l:diff-mat}. We conclude from \eqref{e:dpsin} that
\begin{equation}
\label{e:grad-estim}
4c_0 \int_0^T \int_{\R^3} \left|\nabla_v \sqrt{f^n (t,v)}\right|^2 \frac{\ka (|v|+R_0)}{(1+|v|)^3} \dv \dt \le \bar{\bar{H}}_0:= \bar{H}_0 +\frac{40 k (\bar\delta)}{\bar\delta^3} (2M_0)^2 T\,,
\end{equation}
from which we get \eqref{e:gradfn}.

The Sobolev embedding inequality corresponding to $H^1(\R^3)\subset L^6(\R^3)$, applied to $\chi\sqrt{f^n}$ where the function $\chi\in C^\infty(\R^3)$ is such that $\un_{B_R}\le\chi\le\un_{B_{R+1}}$ with $|\nabla\chi|\le 2$, implies that, for all $R>0$, 
\begin{equation}\label{fL1L3<}
\|f^n\|_{L^1 ((0,T), L^3 ( B_R))} \le C \sqrt{\frac{\bar{\bar{H}}_0 (2+R)^3}{4c_0k(R_0)}+8M_0}
\end{equation}
for some absolute constant $C>0$. Then, interpolating with $f^n \in L^\infty ((0,T),L^1(\R^3))$ yields 
\[ 
\|f^n\|_{L^{5/3} ((0,T) \times B_R)} \le C \left( T \int_{\R^3} \fin (v) \dv  \right)^{\frac25} \left(\frac{\bar{\bar{H}}_0 (2+R)^3}{4c_0k(R_0)}+8M_0\right)^{\frac3{10}}\,.
\]
We have thus established \eqref{e:fn-equi}. 

\paragraph{Compactness.}
As far as \eqref{e:ae-conv} is concerned, we study the gradient and the time derivative: 
\[
\nabla_v f^n = 2 \sqrt{f^n} \nabla_v \sqrt{f^n} 
\]
and in particular, we deduce from \eqref{e:grad-estim} that
\begin{align*}
\int_0^T \| \nabla_v f^n(t,\cdot ) \|_{L^1 (B_R)}^2 \dt & \le 4 \int_0^T \left\{  \int_{\R^3} f^n(t,v) \dv \right\} \left\{ \int_{B_R} |\nabla_v \sqrt{f^n}|^2 \dv \right\} \dt   
\\
&\le\left\{ \int_{\R^3} \fin (v) \dv \right\} \frac{\bar{\bar{H}}_0(1+R)^3}{c_0k(R_0)}\,.
\end{align*}
Now for $\phi \in C^{\infty}_c ( B_R)$, we have
\begin{multline*}
\partial_t\int_{B_R} f^n(t,v) \phi (v)  \dv  =  \frac1n  \int_{B_R} f^n (t,v) \Delta_v \phi (v)  \dv 
\\
-\tfrac12\int_{B_R} \sqrt{a_n (v-w)} \sqrt{f^n (t,v) f^n (t,w)} (\nabla_v \phi (v) - \nabla_w \phi (w)) \cdot F_n(t,v,w)  \dv \dw\,,
\end{multline*}
from which we get 
\begin{align*}
&\left| \int_{B_R} \partial_t f^n(t,v) \phi (v)  \dv \right| 
\\
&\le   \frac{1}n \left\{ \int_{\R^3} \fin (v) \dv \right\}  \| \Delta_v \phi\|_{L^\infty ( B_R)}
\\
&  +  \| D^2 \phi\|_{L^\infty (B_R)} \int_{B_R \times \R^3} \sqrt{|a_n (v-w)|} |v-w| \sqrt{f^n (t,v) f^n (t,w)} |F_n(t,v,w)| \dv \dw   
\\
& \le \frac{1}n \left\{ \int_{\R^3} \fin (v) \dv \right\} \| \Delta_v \phi\|_{L^\infty (B_R)}
\\
& + \| D^2 \phi\|_{L^\infty ( B_R)} \left\{ \int_{B_R\times \R^3} |a_n (v-w)| |v-w|^2 f^n (t,v) f^n (t,w)  \dv \dw  \right\}^{\frac12}\left\{ \iint_{\R^3 \times \R^3} |F^n (t,v,w)|^2 \dv \dw  \right\}^{\frac12}\,. 
\end{align*}

We estimate a part of the right-hand side of the previous inequality. In order to do so, we first remark that $\ka (|z|)|z| \le \ka (1) (1+|z|^2)$. Indeed, for $|z| \le 1$, we have $\ka (|z|)|z| \le \ka (1)$ since $\ka$ is nondecreasing while, for $|z|\ge 1$, 
we have $\ka (|z|)|z| \le \ka (1)|z|^2$, since $\ka(r)/r$ is nonincreasing. With such a piece of information at hand, we now get the following estimate,
\begin{align*}
\int_{B_R\times \R^3} |a_n(v-w)| & |v-w|^2 f^n (t,v) f^n (t,w)  \dv \dw  
\\
& \le \ka(1) \int_{B_R\times \R^3} (1+|v-w|^2) f^n (t,v) f^n (t,w)  \dv \dw  
\\
& \le \ka(1) \int_{B_R\times \R^3} f^n (t,w) (1+ 2R^2+2|w|^2) f^n (t,v)   \dv \dw  
\\
& \le \ka(1) 2(1+R^2) \int_{\R^3}  f^n (t,v) \dv \int_{\R^3} f^n (t,w) (1+|w|^2) \dw  
\\
& \le \ka(1) 2(1+R^2) 4 (M_0+E_0)^2\,,
\end{align*}
where we used in particular \eqref{e:second-m} to get the last inequality. Using \eqref{e:dpsin}, we then conclude that 
\[
\left\| \int_{B_R} \partial_t f^n(t,v) \phi (v) \dv\right\|_{L^2 ([0,T])} \le C \|D^2 \phi \|_{L^\infty (B_R)}
\]
with $C = 2M_0 + \sqrt{\ka(1) 2(1+R^2)} 2 (M_0+E_0) \sqrt{2\bar H_0}$. 

We then deduce  that the sequence $\{f^n\}$ is bounded in $L^2 ((0,T);W^{1,1}(B_R))$ while $\{ \partial_t f^n \}$ is bounded in $L^2 ((0,T), W^{-2,1} (B_R))$. Then Aubin's lemma (see \cite{jsimon}) allows us to conclude that $\{ f^n \}$ is relatively compact 
in $L^2((0,T),L^1 (B_R))$, and in particular, we can extract a subsequence such that \eqref{e:ae-conv} holds true.

\paragraph{Weak limits.} We can also assume that along the subsequence $n_i$, we have $F^{n_i} \rightharpoonup F$ in $L^2 ((0,T) \times \R^3 \times \R^3)$ for some $F$ to be determined. In order to identify $F$, we argue as in \cite{V98} after 
observing that
\begin{equation}\label{PtCvgF}
\sqrt{f^n(v) f^n (w) a_n (v-w)}  \left( \frac{\nabla_v f^n (v)}{f^n(t,v)} - \frac{\nabla_w f^n (w)}{f^n(t,w)} \right) = 2 \Pi (v-w) (\nabla_v -\nabla_w) \bigg(\sqrt{\alpha_n (v-w) f^n (v)f^n(w)} \bigg)\,,
\end{equation}
where $\alpha_n(z):=X(n|z|)\tfrac{k(|z|)}{|z|}$ (so that $a_n(z)=\alpha_n(z)\Pi(z)$). We have omitted the time dependence for the sake of clarity. Thanks to \eqref{e:ae-conv}, we know that, for all $\bar R>0$, the sequence 
$\sqrt{\alpha_{n_i} (v-w) f^{n_i} (t,v)f^{n_i} (t,w)}$ converges a.e. to $\sqrt{\alpha(v-w) f(t,v) f(t,w)}$ in $(0,T) \times B_{\bar R} \times B_{\bar R}$, where $\alpha(z):=\tfrac{k(|z|)}{|z|}$. Uniform integrability is obtained by considering its $L^{2+2\eps}$
norm,
\begin{align*}
\iiint_{(0,T) \times B_{\bar R} \times B_{\bar R}} & |\alpha_n (v-w)|^{1+\eps} (f^n)^{1+\eps} (t,v) (f^n)^{1+\eps} (t,w) \dt \dv \dw 
\\
&\le (\ka (2R))^{1+\eps} \iint_{[0,T] \times B_{\bar R}} (f^n)^{1+\eps} (t,\cdot) \ast_v \frac{\un_{B_{2R}}}{|\cdot|^{1+\eps}}(v) (f^n)^{1+\eps}(t,v)\dt \dv.
\end{align*}
We use Lemma~\ref{l:source-nl} with $\delta=2R$ and $\rho=1+\eps$, with $\tfrac1p=\tfrac{2-\eps}{3(1+\eps)}$ and $\tfrac1q=\tfrac{1+4\eps}{2(1+\eps)}$ so that $\tfrac1p+\tfrac2{3q}=1$, and with $\beta=1+\eps$ and $\tfrac1\theta=\tfrac{1-2\eps}{2(1+\eps)}$
so that $\tfrac1\beta+\tfrac2{3\theta}=\tfrac{4-2\eps}{3(1+\eps)}>1$ for small enough $\eps>0$, coupled to Lemma \ref{l:interpolation} and \eqref{e:gradfn} to obtain an upper bound for the term above. This implies that the convergence \eqref{PtCvgF} holds 
in $L^2 ((0,T) \times B_{\bar R} \times B_{\bar R})$. By uniqueness of the limit in the sense of distributions,
\begin{equation}\label{LimF}
\begin{aligned}
F(t,v,w) = & 2 \Pi (v-w) (\nabla_v -\nabla_w) \bigg(\sqrt{\alpha(v-w) f(t,v)f(t,w)} \bigg)
\\
= & \sqrt{f(t,v) f (t,w) a (v-w)}  \left( \frac{\nabla_v f (t,v)}{f(t,v)} - \frac{\nabla_w f (t,w)}{f(t,w)} \right)\,.
\end{aligned}
\end{equation}

Observing that
$$
\begin{aligned}
|F^{\kappa, \inn, n}_+|\le|F^n|\un_{f^n(t,v)>\kappa}\un_{f^n(t,w)>\kappa}&+2\sqrt{\alpha(v-w)}\sqrt{f^n(t,w)\un_{f^n(t,w)\le\kappa}}\left|\nabla_v\sqrt{f^n(t,v)}\right|\un_{f^n(t,v)>\kappa}
\\
&+2\sqrt{\alpha(v-w)}\sqrt{f^n(t,v)\un_{f^n(t,v)\le\kappa}}\left|\nabla_w\sqrt{f^n(t,w)}\right|\un_{f^n(t,w)>\kappa}\,,
\end{aligned}
$$
we first deduce that $F^{\kappa, \inn, n}_+$ is bounded in $L^2([0,T]\times\R^3\times\R^3)$. (To bound the second and third terms on the right-hand side of the inequality above, use \eqref{e:gradfn} and the bound 
$$
\alpha(z)\le k(1)\left(\tfrac{\un_{|z|<1}}{|z|}+1\right)\in (L^1+L^\infty)(\R^3)\,,\text{ while }f^n\un_{f^n\le\kappa}\in (L^\infty_tL^1_v\cap L^\infty_{t,v})((0,+\infty)\times\R^3)
$$ 
with Young's convolution inequality.) On the other hand
$$
\begin{aligned}
F^{\kappa, \inn, n_i}_+=&2\Pi(v-w)(\nabla_v-\nabla_w)\sqrt{\alpha^\inn_{n_i}(v-w)}\left(\sqrt{f^{n_i}(t,v)}-\sqrt{\kappa}\right)_+\left(\sqrt{f^{n_i}(t,w)}-\sqrt{\kappa}\right)_+
\\
&+2\sqrt{\alpha^\inn_{n_i}(v-w)}\left(\sqrt{f^{n_i}(t,w)}\wedge\sqrt{\kappa}\right)\Pi(v-w)\nabla_v\left(\sqrt{f^{n_i}(t,v)}-\sqrt{\kappa}\right)_+
\\
&-2\sqrt{\alpha^\inn_{n_i}(v-w)}\left(\sqrt{f^{n_i}(t,v)}\wedge\sqrt{\kappa}\right)\Pi(v-w)\nabla_w\left(\sqrt{f^{n_i}(t,w)}-\sqrt{\kappa}\right)_+\,,
\end{aligned}
$$
and one shows that the first term on the right-hand side of this equality converges to 
$$
2\Pi(v-w)(\nabla_v-\nabla_w)\sqrt{\alpha^\inn(v-w)}\left(\sqrt{f(t,v)}-\sqrt{\kappa}\right)_+\left(\sqrt{f(t,w)}-\sqrt{\kappa}\right)_+
$$ 
in $\mathcal D'((0,+\infty)\times\R^3\times\R^3)$ by the same argument as was used on $F^{n_i}$. Because of the $L^2$ bound on $F^{\kappa, \inn, n}_+$ proved above, it is enough to test the sequence $F^{\kappa, \inn, n}_+$ on a dense
family of test functions, of the form $\phi(t,v,w)=\psi(t,v)\chi(w)$. Then one easily checks that
$$
\int_{\R^3}\sqrt{\alpha^\inn_{n_i}(v-w)}\left(\sqrt{f^{n_i}(t,w)}\wedge\sqrt{\kappa}\right)\chi(w)\dd w\to\int_{\R^3}\sqrt{\alpha^\inn_{\delta}(v-w)}\left(\sqrt{f(t,w)}\wedge\sqrt{\kappa}\right)\chi(w)\dd w
$$
a.e. on $(0,+\infty)\times\R^3$, and, using again the inequality $\alpha(z)\le k(1)\left(\tfrac{\un_{|z|<1}}{|z|}+1\right)$, that
$$
\sup_{(t,v)\in(0,+\infty)\times\R^3}\left|\int_{\R^3}\sqrt{\alpha^\inn_{n_i}(v-w)}\left(\sqrt{f^{n_i}(t,w)}\wedge\sqrt{\kappa}\right)\chi(w)\dd w\right|<\infty\,.
$$
By dominated convergence and weak-strong convergence (Proposition 3.5 (iv) in \cite{brezis}), we conclude from \eqref{e:gradfn-conv} that
$$
\begin{aligned}
\int_0^T\iint_{\R^3\times\R^3}\phi(t,v,w)\sqrt{\alpha^\inn_{n_i}(v-w)}\left(\sqrt{f^{n_i}(t,w)}\wedge\sqrt{\kappa}\right)\Pi(v-w)\nabla_v\left(\sqrt{f^{n_i}(t,v)}-\sqrt{\kappa}\right)_+\dd v\dd w\dd t
\\
\to\int_0^T\iint_{\R^3\times\R^3}\phi(t,v,w)\sqrt{\alpha^\inn_{\delta}(v-w)}\left(\sqrt{f(t,w)}\wedge\sqrt{\kappa}\right)\Pi(v-w)\nabla_v\left(\sqrt{f(t,v)}-\sqrt{\kappa}\right)_+\dd v\dd w\dd t&\,.
\end{aligned}
$$
With this, we easily identify the weak limit of $F^{\kappa, \inn, n_i}_+$ and prove that \eqref{e:dissipk-conv} holds.
\end{proof}

\subsection{Construction of suitable weak solutions}\label{sec:lei}

\begin{prop}[Local entropy estimates]\label{p:ee}
Let $\Psi\in C^\infty_c([0,T] \times \R^3)$, and set $\varphi := \Psi^2$. There exists  a negligible set $\mathcal{N} \subset (0,T)$ such that any $H$-solution $f$ obtained in Proposition~\ref{p:approximate} satisfies the truncated
entropy inequality \eqref{e:lee} for all $0<t_1<t_2<T$ with $t_1,t_2 \notin \mathcal{N}$ and all $\kappa\in(1,+\infty)\cap\mathbb Q$.
\end{prop}

Specifically, for each $\delta \in (0,\bar\delta)$ with $\bar\delta$ given by Lemma~\ref{l:diff-mat}, with $\ain(z)$ and $\aout (z)$ defined in \eqref{e:aout}, the local entropy dissipation term is
\begin{equation}
\left. 
\label{e:diff-drift}
\begin{aligned}
- &\iint_{[t_1,t_2] \times \R^3} ( A \nabla_v f - (\nabla_v \cdot A) f) \cdot \left( \frac{\nabla_v f}{f} \varphi + \ln \left(\frac{f}{\kappa}\right) \nabla_v \varphi \right) \un_{\{ f \ge \kappa\} }  \dt \dv 
\\
=&  - \tfrac12 \int_{t_1}^{t_2} \iint_{\R^3 \times \R^3} |F^{\kappa,\inn}_+|^2 \varphi (t,v) \dt \dv \dw -  \int_{t_1}^{t_2}\int_{\R^3} A^\out (t,v) \colon \frac{(\nabla_v \fkp)^{\otimes 2}}{f} (t,v) \varphi (t,v) \dt \dv 
\\
& + \int_{t_1}^{t_2}\int_{\R^3} \fkp (t,v) \varphi (t,v)\left(8\pi\ka(0) \kappa+2\int_{\R^3}\frac{k'(|v-w|)}{|v-w|^2}(f(t,w)\wedge\kappa)\dw\right) \dt \dv
\\
& + \int_{t_1}^{t_2} \left(\mathcal{E}_\out (f)+ \mathcal{E}_\inn (f) \right)(t) \dt\,,
\end{aligned} 
\right\}
\end{equation}
with $\fkp$ defined in \eqref{f-ka+} and $F^{\kappa,\inn}_+$ as in \eqref{e:fkappain}, while
\[ 
A^\out (t,v) = \int_{\R^3} \aout(v-w) f(t,w) \dw \,.
\]
Moreover
\[
\mathcal{E}_\out (f)(t) :=   \int_{\R^3}  h^\kappa_+ (f)(t,v)A^\out(t,v) \colon\nabla^2_v \varphi(t,v)\dv + 2 \int_{\R^3}   h^\kappa_+ (f)(t,v)(\nabla_v \cdot A^\out) (t,v)\cdot \nabla_v \varphi (t,v)  \dv\,,
\]
while   
\begin{align*}
\mathcal{E}_\inn (f)(t) :=& - \iint_{\R^3 \times \R^3}\!\! F(t,v,w) \cdot \sqrt{1\!-\!X\left(\tfrac{|v-w|}\delta\right)}  \sqrt{\ain(v\!-\!w)}\sqrt{ f (t,v) f(t,w) } \ln_+ (f/\kappa) (t,v) \nabla_v \varphi (t,v) \dv \dw 
\\
&  + \tfrac12 \iint_{\R^3 \times \R^3} \left(F^{\kappa,\inn}_+ (t,v,w) \cdot \sqrt{\ain (v-w)} \frac{\nabla_v \fkp}{\sqrt{f}}(t,w)\right)\sqrt{f(t,v)}  (\varphi (t,w) - \varphi (t,v)) \dv \dw  
\\
& - \iint_{\R^3 \times \R^3} (\nabla_z \cdot \ain) (v-w) \cdot \nabla_v \varphi (t,v) (f \wedge \kappa)(t,w) \fkp (t,v) \dv \dw 
\\
& - \iint_{\R^3 \times \R^3} \nabla_z^2 : a^\out(v-w) \fkp (t,w)  \fkp (t,v) \varphi (t,v)  \dv \dw\,,
\end{align*}
with $F$ given by \eqref{e:Fn}. The proof of this rather technical result can be found in Appendix \ref{a:trunc-entropineq}.

\begin{proof}[Proof of Proposition~\ref{p:existence}]
Applying Proposition \ref{p:approximate} provides a $H$-solution $f$; this $H$-solution satisfies condition 1 in the definition of suitable weak solutions (Definition \ref{d:suitable}) because of \eqref{e:dissip-all}. Condition 2 in Definition \ref{d:suitable},
i.e. the local truncated entropy inequality \eqref{e:lee}, follows from Proposition~\ref{p:ee}. Hence all $H$-solutions obtained as limit points as $n\to\infty$ of the sequence $f^n$ of solutions of \eqref{e:landau-reg} with initial data \eqref{e:finn} are
suitable weak solutions of the Landau equation \eqref{e:landau} with initial data \eqref{e:ic}.
\end{proof}

\begin{remark}\label{rmk:VillaniSuitable}
The construction of suitable weak solutions leading to Proposition \ref{p:existence} shows that any Villani solution of \eqref{e:landau}, with the terminology of the introduction (see section \ref{ss:MainRes}), is a suitable weak solution.

\end{remark}

\subsection{Suitable local entropy estimates}\label{ss:SuitLocEntr}

In the next lemma, we use  the local entropy dissipation estimate (Proposition~\ref{p:lede}) to get the suitable form of the general local truncated entropy estimate from the definition of suitable weak solutions, see Definition~\ref{d:suitable},
that we can effectively use in the proof of the De Giorgi lemma. 

We also recall that the constants $\bar\delta$ and $R_0$ given by Lemma~\ref{l:diff-mat} depend only on the hydrodynamic bounds $m_0,M_0,E_0$ and $H_0$. 

\begin{lemma}[Suitable entropy estimates]\label{l:suitable}
Let $\fin\colon \R^3 \to [0,+\infty)$ measurable satisfy \eqref{ini-bnd-bis} and $v_0 \in \R^3$. There exists a constant $C_0\geq 1$ depending only on $|v_0|, m_0,M_0,E_0,H_0,\ka(\bar\delta),\ka(R_0)$ such that any suitable weak solution $f$ 
of the Landau equation in $(0,T) \times \R^3$  satisfies, for all $t_0 \in (0,T)$, all $r \in (0,1)$ and $\delta \in (0,\bar\delta]$ such that $t_0 >(r+\delta)^2 $ and $\kappa \in[1,+\infty)\cap\mathbb{Q}$,
\begin{equation}\label{e:suitable}
\left.\begin{aligned}
&\esssup_{t\in (t_0-r^2,t_0]}   \int_{B_r(v_0)}   h_+^\kappa (f (t,v))  \dv  +  \int_{Q_r(t_0,v_0)} \frac{|\nabla_v f^\kappa_+(t,v)|^2}{f(t,v)} \dt \dv 
\\
&\le C_0 (\kappa + \delta^{-2}) \int_{Q_{r+\delta} (t_0,v_0)} f \bigg( \ln_+ (f/\kappa)+\ln_+^2 (f/\kappa) \bigg)  (t,v) \dt \dv 
\\
& +  C_0 \int_{Q_{r+\delta} (t_0,v_0)} \left(  f(t,\cdot) \ast_v \frac{\ka}{\delta^2|\cdot|}(v) \right) f \bigg( \ln_+ (f/\kappa)+\ln_+^2 (f/\kappa) \bigg)  (t,v)   \dt \dv 
\\
& +  C_0 \int_{Q_{r+\delta} (t_0,v_0)} \left\{  (f \wedge \kappa) (t,\cdot) \ast_v \left( \frac{\ka'}{|\cdot|^2} +  \frac{\ka}{\delta|\cdot|^2} \right) (v) \right\} f \ln_+ (f/\kappa) (t,v)  \dt \dv\,.
\end{aligned}
\right\}
\end{equation}
\end{lemma}

\begin{remark}
The three convolution products appearing in the right-hand side of \eqref{e:suitable} are well defined. Indeed, $\ka/|\cdot|$ and $\ka'/|\cdot|^2$ and $\ka/|\cdot|^2$ are bounded at infinity. Moreover, $\ka'/|\cdot|^2$ and $\ka/|\cdot|^2$ are 
integrable near the origin and $f \wedge \kappa \le \kappa$. To finish with, $k/|\cdot|$ is bounded from above by $1/|\cdot|$ around the origin (up to some positive constant). 
\end{remark}

\begin{proof}
Let  $\Psi \colon \R^3 \to \R$ be $C^\infty$ and satisfy $\un_{B_r(v_0)}\le\Psi\le\un_{B_{r+\delta}(v_0)}$, together with
\begin{equation}
\label{e:diff-psi}
|\nabla_v \Psi| \le c_\ast \delta^{-1} \quad \text{ and } \max_{|e|=1} |\nabla^2_v \Psi \colon e^{\otimes 2}| \le c_\ast \delta^{-2}
\end{equation}
for some absolute constant $c_\ast \ge 1$. We then consider $\varphi = \Psi^2$. 

Since $f$ is a suitable weak solution, we know that there exists a negligible set $\mathcal{N} \subset (0,T)$ such that for all $t_1,t_2 \in (0,T) \setminus \mathcal{N}$ with $t_1 < t_2$, 
\begin{equation}\label{e:dg1lee-init}
\int_{\R^3}   h_+^{\kappa} (f (t_2,v))  \varphi (v) \dv \le  \int_{\R^3}   h_+^{\kappa} (f (t_1,v))  \varphi (v) \dv + \mathcal{T}_1 + \mathcal{T}_2\,,
\end{equation}
with
\begin{align*}
\mathcal{T}_1 &=   - \iint_{[t_1,t_2] \times \R^3} ( A \nabla_v f - (\nabla_v \cdot A) f) \cdot (h_+^{\kappa})''(f) (\nabla_v f)  \Psi^2 \dt \dv\,,
\\
\mathcal{T}_2 &= - \iint_{[t_1,t_2] \times \R^3} ( A \nabla_v f - (\nabla_v \cdot A) f) \cdot (h_+^{\kappa})'(f) \nabla_v \Psi^2  \dt \dv\,.
\end{align*}

\paragraph{The term $\mathcal{T}_1$.} Since all the terms below involve the same time variable $t$, we occasionally drop it in long formulas for the sake of simplicity. We systematically use the notation \eqref{f-ka+}. Thus
\begin{align*}
\mathcal{T}_1 =   &  -  \int_{t_1}^{t_2} \iint_{\R^3 \times \R^3} a (v-w) f(v)f(w)\left( \frac{\nabla_v \fkp}{f}(v) -\frac{\nabla_w f}{f}(w) \right) \cdot \frac{\nabla_v \fkp}{f} (v)  \Psi^2 (v) \dt \dv  \dw 
\\
\intertext{we remark that $f = f \wedge \kappa + \fkp$ and we split $\nabla_w f (w)$ accordingly,} 
= & -  \int_{t_1}^{t_2} \iint_{\R^3 \times \R^3} a (v-w) f(v)f(w)\left( \frac{\nabla_v \fkp}{f}(v) -\frac{\nabla_w \fkp}{f}(w) \right) \cdot \frac{\nabla_v \fkp}{f} (v)  \Psi^2 (v) \dt \dv  \dw
\\
& +  \int_{t_1}^{t_2} \iint_{\R^3 \times \R^3} a (v-w) \nabla_w (f \wedge \kappa)(w) \cdot \nabla_v \fkp (v)  \Psi^2 (v) \dt \dv  \dw 
\\
\intertext{we split $\Psi^2(v)$ as $\Psi(v) (\Psi (w) + [\Psi (v)-\Psi(w)])$,}
= & -  \int_{t_1}^{t_2} \iint_{\R^3 \times \R^3} a (v-w) f(v)f(w)\left( \frac{\nabla_v \fkp}{f}(v) \Psi(v)-\frac{\nabla_w \fkp}{f}(w)\Psi(w)\right)   \cdot \frac{\nabla_v \fkp}{f} (v)  \Psi (v) \dt \dv  \dw 
\\
&  +  \int_{t_1}^{t_2} \iint_{\R^3 \times \R^3} \bigg( a (v-w) (\Psi(v) - \Psi (w)) \Psi (v) \bigg)   \nabla_w \fkp (w) \cdot \nabla_v \fkp (v)   \dt \dv  \dw 
\\
& +  \int_{t_1}^{t_2} \iint_{\R^3 \times \R^3} \bigg( a (v-w) \Psi^2 (v) \bigg) \nabla_w (f \wedge \kappa)(w) \cdot \nabla_v \fkp (v)   \dt \dv  \dw\,.
\end{align*}
After symmetrizing the first term with respect to $v$ and $w$, we recognize the localized entropy dissipation. We also symmetrize the second term and integrate by parts  the second and third terms with respect to $v$ and $w$ and we get 
\begin{align*}
\mathcal{T}_1 =   - \tfrac12 \int_{[t_1,t_2] \times \R^3 \times \R^3} |F_\Psi^\kappa (t,v,w)|^2  \dt \dv \dw  & +  \int_{t_1}^{t_2} \iint_{\R^3 \times \R^3} \Gamma_1 (v,w) \fkp(t,w) \fkp (t,v) \dt \dv \dw 
\\
& +   \int_{t_1}^{t_2} \iint_{\R^3 \times \R^3} \Gamma_2 (v,w) (f \wedge \kappa)(t,w) f^\kappa_+ (t,v)  \dt \dv \dw  \,,
\end{align*}
with
\begin{equation}
\label{e:fpsik}
F_\Psi^\kappa (t,v,w) :=  \sqrt{a (v-w)}\sqrt{ f(t,v)f(t,w)} \left( \frac{\nabla_v \fkp}{f}(v) \Psi(v)-\frac{\nabla_w \fkp}{f}(w)\Psi(w)\right)
\end{equation}
and
\begin{align*}
\Gamma_1 (v,w) := & \tfrac12 (-\nabla_z^2 \colon a) (v-w) (\Psi (v) -\Psi (w))^2  
\\
& + (-\nabla_z \cdot a) (v-w)  (\nabla_v \Psi(v) + \nabla_v \Psi(w))(\Psi(v)-\Psi (w)) 
\\
& - a(v-w) \nabla_w \Psi (w) \cdot\nabla_v \Psi (v) \,,
\\
\Gamma_2 (v,w) & := (-\nabla_z^2 : a) (v-w) \Psi^2 (v) + (-\nabla_z \cdot a) (v-w) \cdot \nabla_v (\Psi^2) (v)\,.
\end{align*}

We manipulate the second and third terms separately. As far as the second term is concerned, we use formulas \eqref{e:div-a}-\eqref{e:hess-a} for derivatives of $a$ and the fact that $\ka'(r)\le \ka (r)/r$ to get
\begin{align*}
|\Gamma_1 (v,w)| \le & \tfrac12 \|\nabla_v \Psi\|_\infty^2 |\nabla^2_z \colon a(v-w)| |v-w|^2(\un_{B_{r+\delta}(v_0)} (v) + \un_{B_{r+\delta}(v_0)} (w))
\\
& + 2\|\nabla_v \Psi\|_\infty^2 |\nabla_z \cdot a(v-w)| |v-w|   (\un_{B_{r+\delta}(v_0)} (v) + \un_{B_{r+\delta}(v_0)} (w)) 
\\
& + \|\nabla_v \Psi\|_\infty^2 |a(v-w)| \un_{B_{r+\delta}(v_0)} (v) \un_{B_{r+\delta}(v_0)} (v) 
\\
\le & \|\nabla_v \Psi\|_\infty^2 \ka' (|v-w|) (\un_{B_{r+\delta}(v_0)} (v) + \un_{B_{r+\delta}(v_0)} (w)) 
\\
&  +  \|\nabla_v \Psi\|_\infty^2 \frac{4 \ka (|v-w|)}{|v-w|} (\un_{B_{r+\delta}(v_0)} (v) + \un_{B_{r+\delta}(v_0)} (w)) 
\\
& + \|\nabla_v \Psi\|_\infty^2 \frac{\ka (|v-w|)}{ |v-w|} \un_{B_{r+\delta}(v_0)} (v) \un_{B_{r+\delta}(v_0)} (w) 
\\
\le & \phantom{+} \frac{6 c_\ast^2}{\delta^2} \frac{\ka (|v-w|)}{|v-w|} (\un_{B_{r+\delta}(v_0)} (v) + \un_{B_{r+\delta}(v_0)} (w))\,.
\end{align*}

With such an estimate at hand, we can write
\begin{align*}
\int_{t_1}^{t_2} \iint_{\R^3 \times \R^3} \Gamma_1 (v,w)  \fkp (t,w)  \fkp (t,v) \dt \dv \dw \le & \frac{12 c_\ast^2}{\delta^2} \int_{t_1}^{t_2} \int_{B_{r+\delta}(v_0)} \left\{ \fkp (t,\cdot) \ast_v \frac{\ka}{|\cdot|} (v) \right\} \fkp (t,v) \dt \dv  
\\
\le & \frac{12 c_\ast^2}{\delta^2} \int_{t_1}^{t_2} \int_{B_{r+\delta}(v_0)} \left\{  f (t,\cdot) \ast_v \frac{\ka}{|\cdot|} (v) \right\} \fkp (t,v) \dt \dv\,.
\end{align*}

Recalling that $(-\nabla_z^2 \colon a) = 8 \pi \ka (0) \delta_0 +2 \ka' (|z|)|z|^{-2}$, see \eqref{e:hess-a}, we have
\begin{align*}
\left|\int_{\R^3} \Gamma_2 (v,w) (f \wedge \kappa) (t,w) \dw - 8 \pi \ka (0) (f \wedge \kappa) (t,v) \Psi^2 (v) - \left\{ (f \wedge \kappa)(t,\cdot) \ast_v \frac{2k'}{|\cdot|^2}(v) \right\} \Psi^2 (v)\right|
\\
\le \left\{ (f \wedge \kappa)(t,\cdot) \ast_v \frac{2k}{|\cdot|^2}(v) \right\} |\nabla_v (\Psi^2)(v)|&\,.
\end{align*}

Combining the different estimates for $\mathcal{T}_1$, we get
\begin{equation}\label{e:t1-suitable}
\left.
\begin{aligned}
\mathcal{T}_1 \le   & - \frac12 \int_{t_1}^{t_2} \iint_{\R^3 \times \R^3} |F_\Psi^\kappa (t,v,w)|^2 \dt \dv \dw + 8 \pi \ka (0) \kappa \int_{Q_{r+\delta} (t_0,v_0)} \fkp (t,v) \dv 
\\
& + \frac{12 c_\ast}{\delta^2} \int_{Q_{r+\delta} (t_0,v_0)} \left\{  f (t,\cdot) \ast_v \frac{k}{|\cdot|} (v) \right\} \fkp (t,v) (t,v) \dt \dv
\\
& + 2 \int_{Q_{r+\delta} (t_0,v_0)} \left\{  (f \wedge \kappa) (t,\cdot) \ast_v \left( \frac{k'}{|\cdot|^2} + 2c_\ast \frac{ k}{\delta|\cdot|^2} \right) (v) \right\} \fkp (t,v)  \dt \dv\,.
\end{aligned}
\right\}
\end{equation}

\paragraph{The term $\mathcal{T}_2$.} We proceed as above for the $\mathcal{T}_1$ term by splitting it as follows: $\mathcal{T}_2 = \mathcal{T}_2^++ \mathcal{T}_2^\wedge$ with
\begin{align*}
\mathcal{T}_2^+ & =  -  \int_{t_1}^{t_2} \iint_{\R^3 \times \R^3} a (v-w) f(v)f(w)\left( \frac{\nabla_v \fkp}{f}(v) -\frac{\nabla_w \fkp}{f}(w) \right) \cdot \ln_+(f/\kappa) (v)  \nabla_v (\Psi^2) (v) \dt \dv  \dw\,,
\\
\mathcal{T}_2^\wedge &=   \int_{t_1}^{t_2} \int_{B_{r+\delta}(v_0)}  \left\{ \int_{\R^3} (\nabla_z \cdot a) (v-w)   (f \wedge \kappa)(t,w) \dw \right\}\cdot  f \ln_+ (f/\kappa) (t,v)  \nabla_v (\Psi^2) (v) \dt \dv\,.
\end{align*}
We obtained $\mathcal{T}_2^\wedge$ through an integration by parts with respect to $w$. 

\paragraph{The term $\mathcal{T}_2^+$.} We make  $F_\Psi^\kappa$ appear, see \eqref{e:fpsik} above, so that
\begin{align*}
\mathcal{T}_2^+  = & - 2 \int_{t_1}^{t_2} \iint_{\R^3 \times \R^3} a (v-w) f(t,v)f(t,w) \left( \frac{\nabla_v \fkp}{f}(t,v) \Psi (v) -\frac{\nabla_w \fkp}{f}(t,w) \Psi(v) \right) 
\\
& \hspace{7cm}          \cdot \ln_+(f/\kappa) (t,v)  \nabla_v \Psi (v) \dt \dv \dw 
\\
= & - 2 \int_{t_1}^{t_2} \iint_{\R^3 \times \R^3} F_{\Psi}^\kappa (t,v,w) \cdot \sqrt{a (v-w) f(t,v)f(t,w)} \ln_+(f/\kappa) (t,v) \nabla_v \Psi (v) \dt \dv  \dw 
\\
& + 2 \int_{t_1}^{t_2} \iint_{\R^3 \times \R^3} a (v-w) f(v) \nabla_w \fkp(t,w) (\Psi(v) -\Psi(w)) \cdot \ln_+(f/\kappa) (t,v) \nabla_v \Psi (v) \dt \dv \dw\,.
\end{align*}
We now use the Cauchy-Schwarz inequality in the first term and integrate by parts with respect to $w$ in the second term in order to get
\begin{align*}
 \mathcal{T}_2^+ \le & \phantom{+} \frac14 \int_{t_1}^{t_2} \iint_{\R^3 \times \R^3} |F_\Psi^\kappa (t,v,w)|^2 \dt \dv \dw 
 \\
& + 4  \int_{t_1}^{t_2} \int_{B_{r+\delta}(v_0)} \left\{ \int_{\R^3} a (v-w) f(t,w) \dw \right\}:f \ln_+^2(f/\kappa) (t,v) \nabla_v \Psi (v)^{\otimes 2} \dt \dv 
\\
& + 2 \int_{t_1}^{t_2} \int_{\R^3} \left\{ \int_{\R^3} \Gamma_3 (v,w) \fkp(t,w) \dw \right\}\cdot f\ln_+ (f/\kappa) (t,v)  \nabla_v \Psi (v) \dt \dv\,,
\end{align*}
with
\[ 
\Gamma_3 (v,w) := \bigg( (\nabla_z \cdot a) (v-w) (\Psi(v) -\Psi(w)) + a (v-w) \nabla_w \Psi (w) \bigg). 
\]
Recalling \eqref{e:div-a} and  $a (z) = \frac{\ka(|z|)}{|z|} \Pi (z)$ and $|\nabla_v \Psi | \le c_\ast \delta^{-1}$, we have
\[
2|\Gamma_3 (v,w)|  \le 2 \frac{c_\ast}{\delta}  \left( \frac{2\ka(|v-w|)}{|v-w|}  +\frac{\ka (|v-w|)}{|v-w|} \right) \le 6 \frac{c_\ast}{\delta}  \frac{\ka(|v-w|)}{|v-w|}\,.
\]
This implies that 
\begin{multline*}
2 \int_{t_1}^{t_2} \int_{\R^3} \left\{ \int_{\R^3} \Gamma_3 (v,w) \fkp(t,w) \dw \right\} \cdot f\ln_+ (f/\kappa) (t,v)  \nabla_v \Psi (v) \dt \dv 
\\
\le  6 \frac{c_\ast^2}{\delta^2} \int_{t_1}^{t_2} \int_{B_{r+\delta}(v_0)} \left( f(t,\cdot) \ast_v \frac{k}{|\cdot|}(v) \right) f\ln_+ (f/\kappa) (t,v) \dt \dv\,.
\end{multline*}
Moreover, we have 
\begin{multline*}
4\int_{t_1}^{t_2} \int_{B_{r+\delta}(v_0)} \left\{ \int_{\R^3} a (v-w) f(t,w) \dw \right\}:f \ln_+^2(f/\kappa) (t,v) \nabla_v \Psi (v)^{\otimes 2} \dt \dv 
\\
\le 4 \frac{c_\ast^2}{\delta^2} \int_{t_1}^{t_2} \int_{B_{r+\delta}(v_0)}\left( f(t,\cdot) \ast_v \frac{k}{|\cdot|}(v) \right) f\ln_+^2 (f/\kappa) (t,v) \dt \dv\,.
\end{multline*}

We thus derive from the previous upper bound for $\mathcal{T}_2^+$ the following one,
\begin{equation}\label{e:t2in}
\left.
\begin{aligned}
\mathcal{T}_2^+\le & \phantom{+} \frac14 \int_{t_1}^{t_2} \iint_{\R^3 \times \R^3} |F_\Psi^\kappa (t,v,w)|^2 \dt \dv \dw 
\\
& + 6 \frac{ c_\ast^2}{\delta^2} \int_{Q_{r+\delta} (t_0,v_0)} \left( f(t, \cdot) \ast_v \frac{k}{|\cdot|}(v) \right)  f \left( \ln_+ (f/\kappa) + \lnp^2 (f/\kappa) \right) (t,v)  \dt \dv\,.
\end{aligned}
\right\}
\end{equation}

\paragraph{The term $\mathcal{T}_2^\wedge$.} 
We use $|\nabla_z \cdot a| = \frac{2 k}{|\cdot|^2}$ to write,
\begin{equation}
\label{e:t2wedge}
\mathcal{T}_2^\wedge \le 2 \frac{c_\ast}{\delta} \int_{t_1}^{t_2} \int_{B_{r+\delta}(v_0)}  \left\{ (f \wedge \kappa)(t,\cdot) \ast_v  \frac{2k}{|\cdot|^2}(v) \right\}f \ln_+ (f/\kappa) (t,v)   \dt \dv\,.
\end{equation}

\paragraph{Intermediate local entropy estimates.}
Combining \eqref{e:dg1lee-init} with \eqref{e:t1-suitable}, \eqref{e:t2in} and \eqref{e:t2wedge} yields
\begin{equation}\label{e:lee-before-mean}
\left.
\begin{aligned}
& \int_{B_r(v_0)}   h_+^\kappa (f (t_2,v))  \dv  + \frac14 \int_{t_1}^{t_2} \int_{\R^3 \times \R^3} |F_\Psi^\kappa (t,v,w)|^2\dt \dv \dw 
\\
\le &   \phantom{+}   \int_{B_{r+\delta}(v_0)}   h_+^\kappa (f (t_1,v))   \dv 
\\
& + c_{\ast \ast} \kappa  \int_{Q_{r+\delta} (t_0,v_0)} f \bigg( \ln_+ (f/\kappa)+\ln_+^2 (f/\kappa) \bigg)  (t,v) \dt \dv 
\\
& + c_{\ast \ast} \int_{Q_{r+\delta} (t_0,v_0)} \left(  f(t,\cdot) \ast_v \frac{k}{\delta^2|\cdot|}(v)\right) f \bigg( \ln_+ (f/\kappa)+\ln_+^2 (f/\kappa) \bigg)  (t,v)  \dt \dv 
\\
& +  c_{\ast \ast} \int_{Q_{r+\delta} (t_0,v_0)} \left\{  (f \wedge \kappa) (t,\cdot) \ast_v \left( \frac{k'}{|\cdot|^2} +  \frac{ k}{\delta|\cdot|^2} \right) (v) \right\} f \ln_+ (f/\kappa) (t,v)  \dt \dv
\end{aligned} 
\right\}
\end{equation}
for some absolute constant $c_{\ast\ast}$. In several places, we have used the elementary inequality $\fkp \le f \ln_+ (f/\kappa)$ (since $h_+ (r) = r \lnp r - (r-1)_+ \ge 0$) to unify the various terms involved. As a result, this inequality is not sharp, 
but nevertheless sufficient for our purpose.

\paragraph{Local entropy dissipation and macroscopic estimate.}
We apply Proposition~\ref{p:lede} with $F=\fkp$ and $\delta = \bar\delta$ given by Lemma~\ref{l:diff-mat}. Recalling that $a (z)\ge a^\out_{\bar\delta}(z)$, we get, for $F_\Psi^\kappa$ defined in \eqref{e:fpsik},
\begin{align}
\nonumber \frac12 \int_{\R^3 \times \R^3} & |F_\Psi^\kappa (t,v,w)|^2 \dv \dw 
\\
\nonumber  & \ge \frac{c_0 \ka(R_0)}{(1+|v_0|+r)^3} \int_{B_r(v_0)} \frac{|\nabla_v f^\kappa_+ (t,v)|^2}{f(t,v)} \dv -  \frac{40 \ka(\bar\delta)}{\bar\delta^{3}} \left\{\int_{B_{r+\delta} (v_0)} \fkp(t,v) (1+c_\ast\bar\delta\delta^{-1})\dv \right\}^2 
\\
\intertext{use now $r,\bar\delta \in (0,1)$ and get}
\label{e:dissip-suitable}
&\ge \frac{c_0 \ka (R_0)}{(2+|v_0|)^3} \int_{B_r(v_0)} \frac{|\nabla_v f^\kappa_+ (t,v)|^2}{f(t,v)} \dv -  \frac{40 \ka(\bar\delta)(1 + c_\ast \bar\delta)^2 M_0}{\bar\delta^3\delta^2} \int_{B_{r+\delta} (v_0)} \fkp(t,v) \dv\,.
\end{align}
  
We now combine \eqref{e:lee-before-mean} with  \eqref{e:dissip-suitable}  and get
\begin{equation}
\left.
\begin{aligned}
\int_{B_r(v_0)} &  h_+^\kappa (f (t_2,v))  \dv  + \int_{t_1}^{t_2}  \int_{B_r(v_0)} \frac{|\nabla_v f^\kappa_+ (t,v)|^2}{f(t,v)} \dv  \dt \dv  \le \bar C_0   \int_{B_{r+\delta}(v_0)}   h_+^\kappa (f (t_1,v))   \dv  
\\
& + \bar C_0 (\kappa + \delta^{-2}) \int_{Q_{r+\delta} (t_0,v_0)} f \bigg( \ln_+ (f/\kappa)+\ln_+^2 (f/\kappa) \bigg)  (t,v) \dt \dv 
\\
& + \bar C_0 \int_{Q_{r+\delta} (t_0,v_0)} \left(  f(t,\cdot) \ast_v \frac{k}{\delta^2|\cdot|}  \right) f \bigg( \ln_+ (f/\kappa)+\ln_+^2 (f/\kappa) \bigg)  (t,v)   \dt \dv 
\\
& +  \bar C_0 \int_{Q_{r+\delta} (t_0,v_0)} \left\{  (f \wedge \kappa) (t,\cdot) \ast_v \left( \frac{k'}{|\cdot|^2} +  \frac{ k}{\delta|\cdot|^2} \right) (v) \right\} f \ln_+ (f/\kappa) (t,v)  \dt \dv
\end{aligned} 
\right\}
\label{e:lee-before-mean-bis}
\end{equation}
with $\bar C_0\ge1$ depending only on $|v_0|$, $m_0,M_0,E_0$, $H_0$ and $\ka (\delta_0)$, $\ka (R_0)$.

\paragraph{Mean in time.} We now compute the mean with respect to time of the previous estimate. More precisely,  choosing $t_0-(r+ \delta)^2<t_1 <t_0 -r^2 <t_2 < t_0$ in \eqref{e:lee-before-mean-bis}, and averaging in $t_1$, we get 
\begin{align*}
\int_{B_r(v_0)}   h_+^\kappa (f (t_2,v))  \dv  &+ \frac1{\delta(2r+\delta)}\int_{t_0 -(r+ \delta)^2}^{t_0-r^2} \left\{ \int_{t_1}^{t_2} \mathcal{D} (\fkp, f^\kappa_+) \dt \right\} \dd t_1 
\\
&\le  \frac{\bar C_0}{\delta(2r+\delta)} \int_{t_0-(r+ \delta)^2}^{t_0} \int_{B_{r+\delta}(v_0)}   h_+^\kappa (f (t_1,v)) \dd t_1  \dv +  \mathcal{S}\,,
\end{align*}
for a.e. $t_2 \in (t_0-r^2,t_0]$, with 
\[  
\mathcal{D} (f^\kappa_+,f^\kappa_+)(t) =\int_{B_r(v_0)} f^{-1}(t,v) |\nabla_v f^\kappa_+(t,v)|^2 \dv \,,
\]
and
\begin{align*}
\mathcal{S} &=  \bar C_0 (\kappa + \delta^{-2}) \int_{Q_{r+\delta} (t_0,v_0)} f \bigg( \ln_+ (f/\kappa)+\ln_+^2 (f/\kappa) \bigg)  (t,v) \dt \dv 
\\
& + \bar C_0 \int_{Q_{r+\delta} (t_0,v_0)} \left(  f(t,\cdot) \ast_v \frac{k}{\delta^2|\cdot|}  \right) f \bigg( \ln_+ (f/\kappa)+\ln_+^2 (f/\kappa) \bigg)  (t,v)   \dt \dv 
\\
& +  \bar C_0 \int_{Q_{r+\delta} (t_0,v_0)} \left\{  (f \wedge \kappa) (t,\cdot) \ast_v \left( \frac{k'}{|\cdot|^2} +  \frac{ k}{\delta|\cdot|^2} \right) (v) \right\} f \ln_+ (f/\kappa) (t,v)  \dt \dv\,.
\end{align*}

Dropping the entropy dissipation term, we get 
\[
\int_{B_r(v_0)}   h_+^\kappa (f (t_2,v))  \dv \le   \frac{\bar C_0}{\delta(2r+\delta)} \int_{t_0-(r+ \delta)^2}^{t_0} \int_{B_{r+\delta}(v_0)}   h_+^\kappa (f (t_1,v)) \dd t_1  \dv + \mathcal{S}\,,
\]
while dropping the entropy term in the left-hand side and picking $t_2 = t_0$ yields,
\begin{align*}
\int_{t_0-r^2}^{t_0}  \mathcal{D} (\fkp, \fkp) \dt = &  \frac1{\delta(2r+\delta)} \int_{t_0 -(r+ \delta)^2}^{t_0-r^2} \left\{ \int_{t_0-r^2}^{t_0} \mathcal{D} (\fkp, \fkp) \dt \right\} \dd t_1 
\\
\le & \frac1{ \delta (2r+\delta)}\int_{t_0 -(r+ \delta)^2}^{t_0-r^2} \left\{ \int_{t_1}^{t_0} \mathcal{D} (\fkp, \fkp) \dt \right\} \dd t_1 
\\
\le & \frac{\bar C_0}{\delta(2r+\delta)} \int_{t_0-(r+ \delta)^2}^{t_0} \int_{B_{r+\delta}(v_0)}   h_+^\kappa (f (t_1,v)) \dd t_1  \dv  +  \mathcal{S}
\\
\le & \frac{\bar C_0}{2\delta^2} \int_{t_0-(r+ \delta)^2}^{t_0} \int_{B_{r+\delta}(v_0)}   h_+^\kappa (f (t_1,v)) \dd t_1  \dv  + \mathcal{S}\,.
\end{align*}
Combining the last two inequalities yields,
\[
\esssup_{t\in (t_0-r^2,t_0]}   \int_{B_r(v_0)}   h_+^\kappa (f (t,v))  \dv  + \int_{t_0-r^2}^{t_0} \mathcal{D} (\fkp, \fkp) \dt  \le  \frac{2\bar C_0}{ \delta^2} \int_{Q_{r+\delta}(t_0,v_0)}   h_+^\kappa (f (t,v)) \dd t  \dv + 2 \mathcal{S}\,.
\]
Since $\bar C_0$ (see  the definition of $\mathcal{S}$) depends only on $|v_0|$, $m_0,M_0,E_0$, $H_0$, $\ka(\delta_0)$, $\ka (R_0)$, we obtained the desired estimate for some constant $C_0$ depending only on these quantities. 
\end{proof}

\bigskip
From now on, we stick to the case $\ka(r)=r^{\gamma+3}$. Recall from the introduction the definition of the scaled solution $f_\eps(t,v)=\eps^2 f(t_0+\eps^2 t,v_0+\eps v)$ for the Landau-Coulomb equation. For the Landau equation 
with $\gamma \in [-3,-2)$, we define the scaling transformation as
\begin{equation}\label{f-eps-gam}
f_\eps(t,v):=\eps^{\gamma+5} f(t_0+\eps^2 t,v_0+\eps v)\,.
\end{equation}
Since we used macroscopic bounds in the previous proof, we need to understand how the suitable entropy estimates we obtained in the previous lemma are modified by the scaling transform \eqref{f-eps-gam}, where we assume that 
$(t_0,v_0) \in (0,T) \times \R^3$ with $0<\eps<1\wedge\sqrt{t_0}$, and that $(t,v)\in Q_1(0,0)$.

\begin{lemma}[Scaled suitable entropy estimates]\label{l:scaled-ee}
Let $0< \eps < \tfrac12\wedge\sqrt{t_0}$. The scaled suitable weak solution $f_\eps$ defined in $Q_1(0,0)$ by \eqref{f-eps-gam} satisfies, for all $\kappa_\eps\in [1,2]\cap\mathbb{Q}$,  all $r_\eps \in (0,2]$ and all $\delta_\eps \in (0,\bar\delta]\subset(0,1)$,
\begin{equation}\label{e:scaled}
\left.
\begin{aligned}
& \esssup_{t\in (-r_\eps^2,0]}   \int_{B_{r_\eps}}   h_+^{\kappa_\eps} (f_\eps (t,v))  \dv  +  \int_{Q_{r_\eps}} \frac{|\nabla_v f^{\kappa_\eps}_{\eps +}(t,v)|^2}{f_\eps(t,v)} \dt \dv 
\\
&\le  C_0 \left(\ka(0)\kappa_\eps+ \frac1{\delta_\eps^2} \right) \int_{Q_{r_\eps+\delta_\eps} } f_\eps \ln_+ (f_\eps/\kappa_\eps)(t,v) \dt \dv 
\\
& + C_0 \int_{Q_{r_\eps + \delta_\eps}} \frac{1}{\delta_\eps^2} Z[f_\eps](t,v) \; f_\eps \bigg(\ln_+ (f_\eps/\kappa_\eps)+\ln_+^2 (f_\eps/\kappa_\eps)\bigg)(t,v) \dt \dv 
\\
& +\frac{C_0}{\delta_\eps^2} \int_{Q_{r_\eps+\delta_\eps}} \left(  f_\eps(t,\cdot) \ast_v \frac{\un_{B_1}}{|\cdot|}  \right) f_\eps \bigg(\ln_+ (f_\eps/\kappa_\eps)+\ln_+^2 (f_\eps/\kappa_\eps)\bigg)(t,v)   \dt \dv 
\end{aligned}
\right\}
\end{equation}
where we recall that $f^{\kappa_\eps}_{\eps +} := (f_\eps -\kappa_\eps)_+$, while $k(r)=r^{\gamma+3}$ and
\begin{equation}\label{Def-Zfeps}
Z[f_\eps](\bar t,\cdot) := f_\eps(\bar t,\cdot)\ast_v \frac{\ka (|\cdot|)\un_{B_1^c}}{|\cdot|}\,,
\end{equation}
and $C_0>1$ depends only on $|v_0|, m_0,M_0,E_0,H_0$  and $\ka (\delta_0), \ka (R_0)$. 
\end{lemma}

\begin{proof}
Without loss of generality and for the sake of clarity, we assume that the unscaled solution $f$ is defined on $[-1,\infty)\times\R^3$, and consider the case $t_0=0$ and $v=0$. We set $\kappa = \eps^{-\gamma-5} \kappa_\eps$, while $r = \eps r_\eps \in (0,1)$ 
and $\delta = \eps \delta_\eps \in (0,\bar\delta]$ (for $0 < \eps < 1/2$)  and we apply the entropy estimates of the previous lemma.  We first compute,
\[
\eps^{\gamma+2} \esssup_{s \in (t_0-r^2,t_0]}    \int_{B_r(v_0)}   h_+^{\kappa} (f (s,w))  \dw =  \esssup_{t \in (-r_\eps^2,0]}    \int_{B_{r_\eps}}   h_+^{\kappa_\eps} (f_\eps (t,v))  \dv\,.
\]
Each term in \eqref{e:suitable} will be multiplied by $\eps^{\gamma+2}$. Recall that $\ka(\eps)=\eps^{\gamma+3}$, so $\ka(0)\neq 0$ only if $\gamma=-3$.

The first term in the right-hand side \eqref{e:suitable} becomes
\begin{align*}
& \eps^{\gamma+2}(\ka(0) \kappa +\delta^{-2})\int_{Q_{r+\delta}} f \ln_+(f/\kappa) (t,v) \dt \dv
\\
=& \phantom{+} \eps^{\gamma+2}(\ka(0)\eps^{\gamma+5} \kappa_\eps +\eps^{-2}\delta_\eps^{-2})\int_{Q_{r+\delta}} \eps^{-\gamma-5} f_\eps(\eps^{-2}t,\eps^{-1}v) \ln_+(f_\eps/\kappa_\eps) (\eps^{-2}t,\eps^{-1}v) \dt \dv
\\
=& \phantom{+} \eps^{\gamma+2}(\ka(0)\eps^{-2} \kappa_\eps +\eps^{-2}\delta_\eps^{-2})\int_{Q_{r_\eps+\delta_\eps}} \eps^{-\gamma-5} f_\eps \ln_+(f_\eps/\kappa_\eps) (\bar t,\bar v) \eps^5 \dd \bar t \dd \bar v
\\
=& \phantom{+} (\ka(0) \kappa_\eps +\delta_\eps^{-2})\int_{Q_{r_\eps+\delta_\eps}}  f_\eps \ln_+(f_\eps/\kappa_\eps) (\bar t,\bar v) \dd \bar t \dd \bar v\,.
\end{align*}
Moreover, the second term in the right-hand side \eqref{e:suitable} becomes
\begin{align*}
&  \eps^{\gamma+2} \frac{1}{\delta^2} \int_{Q_{r+\delta}} \left(  f(t,\cdot) \ast_v \frac{\ka}{|\cdot|} (v)\right) f \bigg(\ln_+ (f/\kappa)+\ln_+^2 (f/\kappa)\bigg)(t,v)   \dt \dv 
\\
=  &  \phantom{+} \frac{\eps^{\gamma+2}}{\eps^2 \delta_\eps^2}  \int_{Q_{r+\delta}} \left( \int_{\R^3} \eps^{-\gamma-5} f_\eps\left(\frac{t}{\eps^2},\frac{w}{\eps}\right)  \frac{\ka (|v-w|)}{|v-w|} \dw \right) 
	\eps^{-\gamma-5} f_\eps \bigg(\ln_+ (f_\eps/\kappa_\eps)+\ln_+^2 (f_\eps/\kappa_\eps)\bigg)\left(\frac{t}{\eps^2},\frac{v}{\eps}\right)   \dt \dv 
\\
=  & \phantom{+} \eps^{-\gamma-10} \frac1{\delta_\eps^2}  \int_{Q_{r_\eps+\delta_\eps}} \left( \int_{\R^3}  f_\eps(\bar t,\bar w) \frac{\ka (\eps |\bar v- \bar w|)}{\eps|\bar v- \bar w|} \eps^3 \dd \bar w   \right)  
	f_\eps \bigg(\ln_+ (f_\eps/\kappa_\eps)+\ln_+^2 (f_\eps/\kappa_\eps)\bigg)(\bar t,\bar v) \eps^5  \dd \bar t \dd \bar v
\\
=  &  \phantom{+} \eps^{-\gamma-3} \frac1{\delta_\eps^2} \int_{Q_{r_\eps+\delta_\eps}} \left( \int_{\R^3}  f_\eps(\bar t,\cdot) \ast_{\bar v}\frac{\ka (\eps |\cdot|)}{|\cdot|}    \right) 
	f_\eps \bigg(\ln_+ (f_\eps/\kappa_\eps)+\ln_+^2 (f_\eps/\kappa_\eps)\bigg)(\bar t,\bar v)  \dd \bar t \dd \bar v
\\
=  &  \phantom{+}  \frac1{\delta_\eps^2} \int_{Q_{r_\eps+\delta_\eps}} \left( \int_{\R^3}  f_\eps(\bar t,\cdot) \ast_{\bar v}\frac{\ka ( |\cdot|)}{|\cdot|}    \right) 
	f_\eps \bigg(\ln_+ (f_\eps/\kappa_\eps)+\ln_+^2 (f_\eps/\kappa_\eps)\bigg)(\bar t,\bar v)  \dd \bar t \dd \bar v\,.
\end{align*}

We split the convolution product as follows:
\begin{align*}
f_\eps(\bar t,\cdot) \ast_{\bar v} \frac{\ka ( |\cdot|)}{|\cdot|}&= f_\eps(\bar t,\cdot) \ast_{\bar v} \frac{\ka ( |\cdot|) \un_{B_1}}{|\cdot|}+f_\eps(\bar t,\cdot) \ast_{\bar v} \frac{\ka (|\cdot|) \un_{B_1^c}}{|\cdot|} 
\\
\intertext{in the first term, we use that $\ka ( r) \le \ka (1)=1$ for $r \in (0,1]$,}
&\le  f_\eps(\bar t,\cdot) \ast_{\bar v} \frac{\un_{B_1}}{|\cdot|}+f_\eps(\bar t,\cdot) \ast_{\bar v} \frac{\ka ( |\cdot|) \un_{B_1^c}}{|\cdot|}\,.
\end{align*}
In the right-hand side of the previous inequality, we recognize the function $Z[f_\eps]$ introduced in the statement of the lemma. In particular, the two terms in the right-hand side above are respectively the third and second ones appearing in the
right-hand side of \eqref{e:scaled}.

Let us have a closer look at the last term in \eqref{e:suitable}:
\begin{align*}
\eps^{\gamma+2} & \int_{Q_{r+\delta}}\!\left\{  (f \wedge \kappa) (t,\cdot) \ast_v \left( \frac{\ka'}{|\cdot|^2} +  \frac{\ka}{\delta|\cdot|^2} \right) (v) \right\} f \ln_+ (f/\kappa) (t,v)  \dt \dv 
\\
= & \eps^{\gamma+2}  \int_{Q_{r+\delta}} \left\{ \int_{\R^3}\frac{f_\eps \wedge \kappa_\eps}{\eps^{\gamma+5}}\left(\frac{t}{\eps^2},\frac{w}{\eps}\right)\left( \frac{\ka'(|v\! -\! w|)}{|v\!-\!w|^2} +  \frac{\ka(|v\!-\!w|)}{\delta|v\!-\!w|^2} \right) \dw \right\} 
	\eps^{-\gamma-5} f_\eps \ln_+ \left(\frac{f_\eps}{\kappa_\eps}\right)\left(\frac{t}{\eps^2},\frac{v}{\eps}\right)\dt \dv 
\\
= & \eps^{-\gamma -3} \int_{Q_{r_\eps+\delta_\eps}} \left\{ \int_{\R^3}  (f_\eps \wedge \kappa_\eps) (\bar t,\bar w)  \left( \frac{\eps \ka'(\eps|\bar v -\bar w|)}{|\bar v-\bar w|^2} +  \frac{\ka(\eps|\bar v-\bar w|)}{\delta_\eps  |\bar v-\bar w|^2} \right)  \dd \bar w \right\} 
  	f_\eps \ln_+ (f_\eps/\kappa_\eps) (\bar t,\bar v)  \dd \bar t \dd \bar v
\\
= &  \int_{Q_{r_\eps+\delta_\eps}} \left\{ \int_{\R^3}  (f_\eps \wedge \kappa_\eps) (\bar t,\bar w)  \left( \frac{ \ka'(|\bar v -\bar w|)}{|\bar v-\bar w|^2} +  \frac{\ka(|\bar v-\bar w|)}{\delta_\eps  |\bar v-\bar w|^2} \right)  \dd \bar w \right\} 
  	f_\eps \ln_+ (f_\eps/\kappa_\eps) (\bar t,\bar v)  \dd \bar t \dd \bar v\,.
\end{align*}

We now split the integral with respect to $w$ in two parts. In order to estimate the integral outside $B_1$, we first remark that $\ka'(|z|)/ |z| \le \ka (|z|)/|z|^2 \le \ka ( |z|)/|z|$ for $|z| \ge 1$, and that $\delta_\eps < 1$. Hence
\begin{align*}
\int_{\R^3}(f_\eps \wedge \kappa_\eps) (\bar t,\bar w) &  \left( \frac{\ka'(|\bar v -\bar w|)}{|\bar v-\bar w|^2} +  \frac{\ka(|\bar v-\bar w|)}{\delta_\eps |\bar v-\bar w|^2} \right)  \dd \bar w 
\\
\le & \kappa_\eps  \int_{B_1}   \left( \frac{\ka'(|z|)}{|z|^2} +  \frac{\ka(1)}{\delta_\eps  |z|^2} \right) \dd z + 2 \int_{\R^3 \setminus B_1(\bar v)} f_\eps (\bar t , \bar w)  \frac{\ka( |\bar v-\bar w|)}{\delta_\eps |\bar v-\bar w|}  \dd \bar w\,.
\end{align*}
We now compute, for each $\delta,\eps\in(0,1)$
\[ 
\int_{B_1}   \left( \frac{\ka'(|z|)}{|z|^2} +  \frac{\ka(1)}{\delta_\eps  |z|^2} \right) \dd z = 4 \pi  \int_0^1 ( \ka'( r) + 1/\delta_\eps ) \dd r = 4\pi ( \ka (1) - \ka (0)+ 1/\delta _\eps) \le 8 \pi \delta_\eps^{-1} \,.  
\]
Using now that $\kappa_\eps \le 2$, we finally get
\[
\int_{\R^3}  (f_\eps \wedge \kappa_\eps) (\bar t,\bar w) \left( \frac{ \ka'(|\bar v -\bar w|)}{|\bar v-\bar w|^2} +  \frac{\ka(|\bar v-\bar w|)}{\delta_\eps |\bar v-\bar w|^2} \right)  \dd \bar w 
	\le \frac{16 \pi }{\delta_\eps} + \frac2\delta_\eps f_\eps(\bar t,\cdot)\ast_v \frac{\ka ( |\cdot|)\un_{B_1^c}}{|\cdot|}(\bar v).  
\]
We recognize in the right-hand side the function $Z[f_\eps]$ from the statement. The proof of Lemma \ref{l:scaled-ee} is now complete. 
\end{proof}

\subsection{Local mass estimates}
\label{s:mass}

We recall that $\eps$-scaled suitable weak solutions are defined by \eqref{f-eps-gam} for some suitable weak solution $f$ of the Landau equation.

\begin{lemma}[Local mass estimates]\label{l:mass}
Let $f_\eps$ be an $\eps$-scaled suitable weak solution of the Landau equation in $Q_1(0,0)$ for some $\eps\in (0,1)$. Then, for all $\lambda \in (0,1/4)$, 
\begin{multline*}
\frac{1}{C_1 \lambda^3}    \| f_\eps \|_{L^\infty_t L^1_v (Q_{2\lambda})} \le  \left(1+ \frac1{\lambda^4}   \| F_\eps\|_{L^2 (Q_2 \times \R^3)}\right)   \| f_\eps \|_{L^\infty_t L^1_v (Q_2)} 
\\
+  \frac1{\lambda^8}\left(  \| F_\eps\|_{L^2 (Q_2 \times \R^3)}^2 +1 \right) \left\|\nabla_v \sqrt{f_\eps} \right\|_{L^2 (Q_2)}^2  +  \frac{1}{\lambda^8} \eps^{\gamma+2}  \| F_\eps\|_{L^2 (Q_2 \times \R^3)}^2\,,
\end{multline*}
where
\[ 
F_\eps(t,v,w) := \sqrt{a(v-w)} \sqrt{f_\eps(t,v)f_\eps(t,w)} (\nabla_v (\ln f_\eps)(t,v) - \nabla_w (\ln f_\eps)(t,w)), 
\]
with $\sqrt{a}$ defined in \eqref{sqrta} and where the constant $C_1$ depends only on $M_0$.
\end{lemma}

\begin{proof}
The proof is split in several steps. We first introduce the test function that will be used. We next write the corresponding local inequalities and estimate error terms appearing in their right-hand side. 

\paragraph{Test function.} Let $\Psi\in C^\infty(\R^3)$ satisfy $\un_{B_{1/2}}\le\Psi\le 1_{B_1}$ and let $\varphi := \Psi^2$. It is possible to choose $\Psi$ such that $|\nabla_v \Psi | \le 3$, so that $|\nabla_v \varphi | \le 6$. 
Let also $\M$ be defined in $[-1,0] \times \R^3$ as follows: 
\begin{equation}\label{Def-M}
\M (t,v) := (\lambda^2-t)^{-3/2} e^{-\frac{|v|^2}{4 (\lambda^2-t)}}\,.
\end{equation}
In particular, $\partial_t\M = - \Delta_v\M$. We shall use the lower bound
\begin{equation}\label{e:M-lower}
\M \ge \frac1{(5\lambda^2)^{3/2}}e^{-1} = \frac{\bar c_0}{\lambda^3}\qquad \text{ in } Q_ {2\lambda}(0,0)
\end{equation}
for some absolute constant $\bar c_0>0$, and the upper bounds
\begin{equation}\label{e:gradM}
\M  \le \lambda^{-3}\quad \text{ in } Q_1 (0,0)\,,\qquad\qquad|\nabla_v\M| \le \bar c_1 \lambda^{-4} \quad\text{ in } Q_1(0,0)\,,
\end{equation}
where $\bar c_1$ is another absolute constant. We notice that $\nabla_v \varphi (v) \neq 0$ only for $v \in B_{1}(0) \setminus B_{1/2}(0)$ and
\[ 
|\nabla_v\M| = \frac1{2|v|^4} \frac{|v|^5}{(\lambda^2-t)^{5/2}} e^{-\frac{|v|^2}{4 (\lambda^2 -t)}}\le \bar c_2 /6\quad\text{ outside } B_{1/2}\,,
\]
where  $\bar c_2$ is an absolute constant. Since $|\nabla_v \varphi | \le 6$, we have,
\begin{equation}\label{e:test1}
|\nabla_v \varphi \cdot \nabla_v\M | \le \bar c_2\qquad \text{ in } Q_1(0,0)\,.
\end{equation}
Using again that $|\nabla_v \varphi | \le 6$, we estimate $\nabla_v (\varphi M)$ in $Q_1(0,0)$ as follows:
\[
|\nabla_v (\varphi\M) | \le 6\M + |\nabla_v\M| \le  (6 + \bar c_1 \lambda^{-1}) \lambda^{-3}\,.
\]
Since $\lambda \in (0,1)$, we have
\begin{equation}\label{e:test2}
|\nabla_v (\varphi\M) | \le \bar c_3 \lambda^{-4} \,,
\end{equation}
with $\bar c_3 = 6+\bar c_1$.

\paragraph{Lower and upper bounds.} Since $f_\eps$ is an $\eps$-scaled suitable weak solution of the Landau equation in $Q_1(0,0)$, we know that there exists a negligible set $\mathcal{N}\subset (-1,0)$ such that, 
for all $t_1,t_2 \in (-1,0) \setminus \mathcal{N}$, 
\begin{equation}
\label{e:lme0}
\left. 
\begin{aligned}
\int_{\R^3}   f_\eps (t_2,v)  \varphi(v) \M (t_2,v) \dv \le & \phantom{+}  \int_{\R^3}   f_\eps (t_1,v)  \varphi(v) \M (t_1,v)  \dv - \iint_{[t_1,t_2] \times \R^3} f_\eps  \varphi \Delta_v \M  \dt \dv
\\
&- \iint_{[t_1,t_2] \times \R^3} ( A_\eps \nabla_v f_\eps - (\nabla_v \cdot A_\eps) f_\eps) \cdot  \nabla_v (\varphi \M)  \dt \dv\,,
\end{aligned}
\right\}
\end{equation}
where
\[
A_\eps (t,v):= \int_{\R^3} a(v-w) f_\eps(t,w) \dw\,.
\]
We use \eqref{e:M-lower} for $t_2 \in (-4\lambda^2,0]$  to get, 
\begin{equation}\label{LowMasst2}
\int_{\R^3}   f_\eps (t_2,v)  \varphi(v) \M (t_2,v) \dv  \ge \frac{\bar c_0}{\lambda^3} \int_{B_{2\lambda}}   f_\eps (t_2,v) \dv \,,
\end{equation}
while, for $t_1 \in (-1,-1+\lambda^2]$, we have $\lambda^2 -t_1 \ge 1$ so that $\M(t_1,v) \le 1$.
In particular,
\begin{equation}\label{UpMasst1}
\int_{\R^3}   f_\eps (t_1,v)  \varphi(v) \M (t_1,v)  \dv \le  \int_{B_1} f_\eps (t_1,v) \dv\,.
\end{equation}

\paragraph{Time derivative of the test function.}
We continue with the second term in the right-hand side of \eqref{e:lme0}. Integrating by parts with respect to $v$, we get
\[
- \iint_{[t_1,t_2] \times \R^3} f_\eps  \varphi \Delta_v\M  \dt \dv = \iint_{[t_1,t_2] \times \R^3} f_\eps  (\nabla_v \varphi \cdot \nabla_v\M ) \dt \dv + \iint_{[t_1,t_2] \times \R^3}  ( \nabla_v f_\eps \cdot \nabla_v\M ) \varphi \dt \dv\,.
\]
Moreover, \eqref{e:gradM} implies
\begin{align*}
\iint_{[t_1,t_2] \times \R^3} ( \nabla_v f_\eps \cdot \nabla_v\M ) \varphi \dt \dv & \le 2 \bar c_1\lambda^{-4}  \iint_{Q_1}  \sqrt{f_\eps} \left|\nabla_v \sqrt{f_\eps} \right| \dt \dv 
\\
& \le 2 \bar c_1 \lambda^{-4} \left\| \nabla_v \sqrt{f_\eps}\right\|_{L^2 (Q_1)} \| f_\eps \|_{L^\infty_t L^1_v (Q_1)}^{\frac12}\,. 
\end{align*}
Using \eqref{e:test1}, we conclude that
\begin{equation}\label{e:lme1}
- \iint_{[t_1,t_2] \times \R^3} f_\eps  \varphi \Delta_v\M  \dt \dv \le  \bar c_2  \| f_\eps \|_{L^\infty_t L^1_v (Q_1)} + 2 \bar c_1 \lambda^{-4} \| f_\eps \|_{L^\infty_t L^1_v (Q_1)}^{\frac12} \left\| \nabla_v \sqrt{f_\eps}\right\|_{L^2 (Q_1)} . 
\end{equation}

\paragraph{Diffusion and drift terms.}
We finally estimate the last term in \eqref{e:lme0} as follows:
\begin{align*}
- \iint_{[t_1,t_2] \times \R^3} & ( A_\eps \nabla_v f_\eps - (\nabla_v \cdot A_\eps) f_\eps) \cdot  \nabla_v (\varphi\M)  \dt \dv 
\\
& = - \int_{t_1}^{t_2} \int_{\R^3 \times \R^3} a(v-w) f_\eps(t,v) f_\eps(t,w) \left(\frac{\nabla_v f_\eps}{f_\eps} (t,v) - \frac{\nabla_w f_\eps}{f_\eps}(t,w)\right) \cdot \nabla_v (\varphi\M)  \dt \dv 
\\
& = - \int_{t_1}^{t_2} \int_{\R^3 \times \R^3} F_\eps(t,v,w) \cdot \sqrt{a(v-w)} \sqrt{f_\eps(t,v) f_\eps(t,w)}   \nabla_v (\varphi\M)  \dt \dv \dw\,,
\\
\intertext{and we use \eqref{e:test2} to get}
& \le \bar c_3 \lambda^{-4} \| F_\eps\|_{L^2 (Q_1 \times \R^3)} \left\{ \int_{Q_1} \left( f_\eps(t,\cdot) \ast_v  \frac{\ka(|\cdot|)}{ |\cdot|}(v) \right) f_\eps(t,v) \dt \dv \right\}^{\frac12}
\\
& \le \bar c_3 \lambda^{-4} \| F_\eps\|_{L^2 (Q_1 \times \R^3)} \left\|  f_\eps \ast_v  \frac{\ka(|\cdot|)}{|\cdot|} \right\|_{L^1_t L^\infty_v (Q_1)}^{\frac12} \| f_\eps \|_{L^\infty_t L^1_v (Q_1)}^{\frac12}
\\
& \le \bar c_3 \lambda^{-4} \| F_\eps\|_{L^2 (Q_1 \times \R^3)}\left(  \left\|  f_\eps \ast_v  \frac{\un_{B_1}}{|\cdot|} \right\|_{L^1_t L^\infty_v (Q_1)} +M_0 \eps^{\gamma+2}\right)^{\frac12} \| f_\eps \|_{L^\infty_t L^1_v (Q_1)}^{\frac12}\,.
\end{align*}
The last inequality is obtained as follows:
\begin{align*}
f_\eps(t,\cdot) \ast_v  \frac{\ka(|\cdot|)}{|\cdot|} (v) & = \int_{B_1} \frac{\ka (|z|)}{|z|} f_\eps (t,v+z) \dz + \int_{\R^3 \setminus B_1} \frac{\ka (|z|)}{|z|} f_\eps (t,v+z) \dz\,,
\\
\intertext{and we use that $\ka(|z|)=|z|^{\gamma+3}\le 1$ in $B_1$ since $\gamma\ge-3$, and that $\ka (|z|) / |z| = |z|^{\gamma+2} \le 1$ outside $B_1$ since $\gamma<-2$ to get}
& \le \left(  f_\eps(t,\cdot) \ast_v  \frac{\un_{B_1}}{|\cdot|} (v) \right) + M_0 \eps^{\gamma+2} \,. 
\end{align*}
The change of variables $\bar z = \eps z$ yields the desired bound for the second term above. To control the short-range contribution, we use the second inequality in Lemma \ref{l:source-nl} with $\rho=1$ and $\delta=1$, with Lemma \ref{l:interpolation}:
\[
\begin{aligned}
- \iint_{[t_1,t_2] \times \R^3}  ( A_\eps \nabla_v f_\eps - (\nabla_v \cdot A_\eps) f_\eps) \cdot  \nabla_v (\varphi M)  \dt \dv&
\\
\le \bar c_3 \lambda^{-4} \| F_\eps\|_{L^2 (Q_1 \times \R^3)} \left(\left(\tfrac{8\pi}3\right)^{2/3}\|f_\eps\|_{L^1_tL^3_v(Q_2)}+ M_0 \eps^{\gamma+2}\right)^{\frac12} \| f_\eps \|_{L^\infty_t L^1_v (Q_1)}^{\frac12}&
\\
\le \bar c_3 \lambda^{-4} \| F_\eps\|_{L^2 (Q_1 \times \R^3)} \left( \left(\tfrac{8\pi}3\right)^{2/3}\Cinterpol\left(\| f_\eps \|_{L^\infty_t L^1_v (Q_2)} + \left\| \nabla_v \sqrt{f_\eps} \right\|_{L^2(Q_2)}^2\right)
+ M_0 \eps^{\gamma+2}\right)^{\frac12} \| f_\eps \|_{L^\infty_t L^1_v (Q_1)}^{\frac12}&\,.
 \end{aligned}
 \]

\paragraph{Conclusion.} Combining this estimate with \eqref{e:lme0}, \eqref{LowMasst2}, \eqref{UpMasst1} and \eqref{e:lme1} and considering the essential supremum with respect to $(t_1,t_2) \in (-1,-1+\lambda^2]\times(-4\lambda^2,0]$ leads to the following inequality, for some absolute 
constant $\bar{C}_1$,
\begin{align*}
\frac{1}{\bar{C_1} \lambda^3}    \| f_\eps \|_{L^\infty_t L^1_v (Q_{2\lambda})} \le &    \| f_\eps \|_{L^\infty_t L^1_v (Q_2)} +  \lambda^{-4} \| f_\eps \|_{L^\infty_t L^1_v (Q_2)}^{\frac12}\left\| \nabla_v \sqrt{f_\eps}\right\|_{L^2 (Q_2)} 
\\
& +\!  \lambda^{-4} \| F_\eps\|_{L^2 (Q_2 \times \R^3)}\left(  \| f_\eps \|_{L^\infty_t L^1_v (Q_2)}^{\frac12} \!+\!  \| \nabla_v \sqrt{f_\eps}\|_{L^2(Q_2)}  \!+\! \sqrt{M_0 \eps^{\gamma+2}}   \right) \| f_\eps \|_{L^\infty_t L^1_v (Q_2)}^{\frac12}\,,
\\
\intertext{and we use the elementary inequality $ab \le a^2 + b^2/4$  to get}
\le& 4 \| f_\eps \|_{L^\infty_t L^1_v (Q_2)} +  \frac1{4\lambda^{8}}   \left\| \nabla_v \sqrt{f_\eps}\right\|_{L^2 (Q_2)}^2 + \lambda^{-4} \| F_\eps\|_{L^2 (Q_2 \times \R^3)}  \| f_\eps \|_{L^\infty_t L^1_v (Q_2)}  
\\
& + \frac{1}{4 \lambda^8} \| F_\eps\|_{L^2 (Q_2 \times \R^3)}^2  \left\|\nabla_v \sqrt{f_\eps} \right\|_{L^2 (Q_2)}^2 + \frac1{4\lambda^8}  M_0 \eps^{\gamma+2}\| F_\eps\|_{L^2 (Q_2 \times \R^3)}^2\,.
\end{align*}
Rearranging terms yields the desired estimate with a constant $C_1$ depending only on $M_0$. 
\end{proof}


\section{De Giorgi lemma and partial regularity}
\label{s:degiorgi}


\subsection{De Giorgi lemma}

\begin{lemma}[De Giorgi]\label{l:dg1}
Let $f_\eps$ be an $\eps$-scaled suitable weak solution of the Landau equation around $(t_0,v_0)$  for $\eps \in (0,1)$, as in \eqref{f-eps-gam}, such that $\fin$ satisfies \eqref{ini-bnd-bis} and
\[ 
\esssup_{Q_1} \left( f_\eps \ast_v \frac{\ka ( |\cdot|)\un_{B_1(0)^c}}{|\cdot|} \right) \le Z_\eps 
\]
for some $Z_\eps \ge 1$. There exists $\etadg \in  (0,1)$ depending only on $m_0,M_0,E_0, H_0, |v_0|$ and $\ka(\delta_0),\ka(R_0)$ such that
\[
\begin{aligned}
\esssup_{t \in (-4,0]} \int_{B_2} (f_\eps-1)_+ (t,v) \dv + \iint_{Q_2} \left|\nabla_v \sqrt{f_\eps}\right|^2 \un_{\{ f_\eps \ge 1\}}\dt \dv\le \etadg \;  Z_\eps^{-\frac32}
\\
\implies f_\eps \le 2 \quad\text{ a.e. in }Q_{1/2}&\,.
\end{aligned}
\]
\end{lemma}

\begin{proof}
The proof relies on an iterative procedure similar to De Giorgi's solution of Hilbert's 19\textsuperscript{th} problem \cite{DeG}, and sketched below for the reader's convenience.

\paragraph{Sketch of the iterative procedure.}
For shrinking cylinders of radius $r_j$ and  increasing truncation levels $\kappa_j$, we study the local quantity
\[ 
U_j := \esssup_{t \in (t_j,0]} \int_{B_{r_j}} (f_\eps-\kappa_j)_+ (t,v) \dv + \iint_{Q^j} \left| \nabla_v \sqrt{f_\eps} \right|^2 \un_{\{ f_\eps \ge \kappa_j \}}  \dt \dv\,,
\]
where $t_j := -r_j^2$ and $Q^j :=Q_{r_j}$. We aim at proving that $U_j$ goes to $0$ as $j \to \infty$ by establishing the following inductive estimate,
\[ 
U_{j+1} \le C_2^j \left( Z_\eps U_j^{5/3} +  U_j^{4/3} \right)\,,\qquad j\ge 0\,,
\]
for some constant $C_2>0$. Since
\begin{equation}\label{e:U0}
U_0 = \esssup_{t \in (-1,0]} \int_{B_1} (f_\eps-1)_+(t,v) \dv + \iint_{Q_1} \left|\nabla_v \sqrt{f_\eps}\right|^2 \un_{\{ f_\eps \ge 1\}}\dt \dv\le \etadg  Z_\eps^{-\frac32},
\end{equation}
we can choose $\etadg$ small enough so that $U_j \to 0$ as $j\to\infty$. With Fatou's lemma, this implies in particular that $f_\eps \le 2$ a.e. in $Q_{1/2}$. This will be explained in detail at the end of the proof. 

\bigskip
We choose $r_j$ and $\kappa_j$ such that $r_j$ decreases from $1$ to $1/2$ and that $\kappa_j$ increases from $1$ to $2$,
\[ 
r_j := \frac12(1 + 2^{-j}), \qquad \kappa_j :=   2- 2^{-j}, \qquad \text{ for } j \ge 0\,. 
\]
In particular $Q^j = (t_j,0] \times B_{r_j}$.

\paragraph{Entropy inequalities.}
We now apply Lemma~\ref{l:scaled-ee} with $r_\eps = r_{j+1}$ and $\delta_\eps = \delta_j = r_j-r_{j+1}=2^{-j-2}$ and $\kappa_\eps = \kappa_{j+\frac12} = \frac12 (\kappa_j+\kappa_{j+1})$. Remark that in this case 
$r_\eps + \delta_\eps/2 = (r_j + r_{j+1})/2=:r_{j+\frac12}$ and $r_\eps + \delta_\eps = r_j$, while $-(r_\eps+\delta_\eps)^2 = t_j$. Thus, we get from \eqref{e:scaled} with $\kappa_\eps = \kappa_{j+\frac12} $ and from the 
assumption on $f_\eps \ast_v (\un_{B_1^c} |\cdot|^{-2})$ that,
\begin{align}
\nonumber 
& \esssup_{t \in (t_{j+1},0]} \int_{B_{r_{j+1}}} h^{\kappa_{j+\frac12}}_+ (f_\eps(t,v)) \dv +  \int_{Q^{j+1}} \frac{|\nabla_v f_{\eps,+}^{\kappa_{j+\frac12}}(t,v)|^2}{f_\eps(t,v)} \dt \dv 
\\
\label{e:rhs1}  \le  &\,  3 C_0 Z_\eps 2^{2j +4} \iint_{Q^j} f_\eps\bigg( \ln_+ (f_\eps/{\kappa_{j+\frac12}}) + \ln_+ (f_\eps/{\kappa_{j+\frac12}})^2 \bigg)(t,v) \dt \dv 
\\
\label{e:rhs4}  & + C_0 2^{2j+4}  \iint_{Q^j} \left(  f_\eps(t,\cdot) \ast_v \frac{\un_{B_1}}{|\cdot|}(v)  \right) f_\eps \bigg(  \ln_+ (f_\eps/{\kappa_{j+\frac12}}) + \ln_+^2 (f_\eps/{\kappa_{j+\frac12}}) \bigg)(t,v)  \dt \dv
\end{align}
for some constant $C_0$ depending on $|v_0|, m_0,M_0,E_0, H_0$ and $\ka(\eps_0),\ka(R_0)$, and with $\bar\delta,R_0$ given by Lemma~\ref{l:diff-mat}. We have used that $Z_\eps \ge 1$.

\paragraph{Lower bound of the left-hand side.}
We first remark that, since $h_+^{\kappa}$ is convex and nonnegative for each $\kappa >0$, we have
\[ 
\left(h_+^{\kappa_{j+\frac12}}\right)'(\kappa_{j+1}) (f_\eps-\kappa_{j+1}) \le h_+^{\kappa_{j+\frac12}} (f_\eps) - h_+^{\kappa_{j+\frac12}} (\kappa_{j+1}) \le h_+^{\kappa_{j+\frac12}} (f_\eps)\,.
\]
Since $\left(h_+^{\kappa_{j+\frac12}}\right)'(\kappa_{j+1}) = \ln (\kappa_{j+1}/\kappa_{j+\frac12}) \ge \ln (1+ \frac{2^{-j}}{8}) \ge  \cln 2^{-j}$ with
\[
\cln = \frac18 \min_{r \in [0,1/2]} r^{-1} \ln (1+r) = \ln (3/2)/4 <1 \qquad \text{(absolute)}\,,
\] 
we have
\[
\cln 2^{-j} f_{\eps,+}^{\kappa_{j+1}}\le h_+^{\kappa_{j+\frac12}} (f_\eps)\,,
\]
and in particular,
\begin{equation}\label{e:lower-bound-Uk}
\cln 2^{-j} U_{j+1} \le \esssup_{t \in (t_{j+1},0]} \int_{B_{r_{j+1}}} h^{\kappa_{j+\frac12}}_+ (f_\eps(t,v)) \dv +  \iint_{Q_{r_{j+1}}} \frac{|\nabla_v f_{\eps,+}^{\kappa_{j+\frac12}}(t,v)|^2}{f_\eps(t,v)} \dt \dv\,.
\end{equation}
 
We now estimate the two terms in the right-hand side of the local entropy inequality by powers of $U_j$.

\paragraph{Nonlinearization of the term \eqref{e:rhs1}.} We first remark that
\[
(1+r_+) \bigg( \ln (1+r_+) + \ln (1+r_+)^2 \bigg) \le \clnb (r_+ + r_+^{5/3})
\]
for some absolute constant $\clnb$. Since $\kappa_{j+\frac12} \ge 1$, applying this inequality to $r = f_\eps/\kappa_{j+\frac12}-1$ yields
\begin{align*}
f_\eps(t,v) \bigg(\ln_+ (f_\eps/{\kappa_{j+\frac12}}) (t,v) + \ln_+ (f_\eps/{\kappa_{j+\frac12}})^2 (t,v) \bigg) \le& \clnb \kappa_{j+\frac12} \left( \left( \frac{f_\eps}{\kappa_{j+\frac12}}-1 \right)_+ +\left( \frac{f_\eps}{\kappa_{j+\frac12}}-1 \right)_+^{5/3} \right)
\\
\le& \clnb \left( f_{\eps,+}^{\kappa_{j+\frac12}} + \left(f_{\eps,+}^{\kappa_{j+\frac12}}\right)^{5/3} \right)\,.
\end{align*}
Set $\kappa_{j+\frac14}:=\tfrac12(\kappa_{j+\frac12}+\kappa_j)$. Observe that $f_{\eps,+}^{\kappa_{j+\frac12}}\le f_{\eps,+}^{\kappa_{j+\frac14}}$ since $\kappa_{j+\frac12} \ge \kappa_{j+\frac14}$, and
\[ 
f_{\eps,+}^{\kappa_{j+\frac12}} \le (f_\eps-\kappa_{j+\frac14}) \un_{f_\eps-\kappa_{j+\frac14} \ge \kappa_{j+\frac12} - \kappa_{j+\frac14}} \le \frac{(f_\eps-\kappa_{j+\frac14})^{5/3}_+}{(\kappa_{j+\frac12}-\kappa_{j+\frac14})^{2/3}}
= 2^{\frac23 (j+3)} (f_\eps-\kappa_{j+\frac14})^{5/3}_+\,.
\]
Combining the two previous inequalities yields
\[
f_\eps(t,v) \bigg(\ln_+ (f_\eps/{\kappa_{j+\frac12}}) (t,v)  + \ln_+ (f_\eps/{\kappa_{j+\frac12}})^2 (t,v) \bigg)\le 8\clnb 2^{\frac23 j} \left(f_{\eps,+}^{\kappa_{j+\frac14}}\right)^{\frac53}\,.
\]
This implies in turn that the term in \eqref{e:rhs1} is bounded from above as follows:
\[
\iint_{Q^j} f_\eps (t,v)  \bigg((\ln (f_\eps/{\kappa_{j+\frac12}}))_+ (t,v)  + (\ln (f_\eps/{\kappa_{j+\frac12}}))_+^2 (t,v) \bigg) \dt \dv \le  8\clnb 2^{\frac23 j} \left\| f_{\eps,+}^{\kappa_{j+\frac14}} \right\|_{L^{5/3}(Q^j)}^{\frac53}\,. 
\]
Using Lemma~\ref{l:nonlin} with $p=q=5/3$ while $I= (-r_j^2,0]$ and $\mathcal{B} = B_{r_j}$ with $r_j\le 1$, we find that
\[
\left\| f_{\eps,+}^{\kappa_{j+\frac14}} \right\|_{L^{5/3}(Q^j)}\le\Cnonlin 2^{19}2^{4j}U_j\,,
\]
and therefore
\begin{equation}\label{e:rhs1-result}
\iint_{Q^j} f_\eps (t,v)  \bigg((\ln (f_\eps/{\kappa_{j+\frac12}}))_+ (t,v)  + (\ln (f_\eps/{\kappa_{j+\frac12}}))_+^2 (t,v) \bigg) \dt \dv \le \left(\clnb 2^{\frac{104}3} \Cnonlin^{\frac53} \right) 2^{\frac{22}3 j} U_j^{\frac53}\,.
\end{equation}

\paragraph{Nonlinearization of the term \eqref{e:rhs4}.}
Similarly
\[  
(1+r_+) \bigg( \ln (1+r_+) + \ln (1+r_+)^2 \bigg) \le \clnt (r_+ + r_+^{4/3})
\]
for some absolute constant $\clnt \ge 1$. Since $\kappa_{j+\frac12} \ge 1$, applying this inequality to $r = f_\eps/\kappa_{j+\frac12}-1$ yields
\begin{align*}
f_\eps(t,v)  \bigg((\ln (f_\eps/{\kappa_{j+\frac12}}))_+ (t,v) + (\ln (f_\eps/{\kappa_{j+\frac12}}))_+^2 (t,v) \bigg) &\le \clnt \kappa_{j+\frac12} \left( \left( \frac{f_\eps}{\kappa_{j+\frac12}}-1 \right)_+ +\left( \frac{f_\eps}{\kappa_{j+\frac12}}-1 \right)_+^{4/3} \right)
\\
& \le \clnt \left( f_{\eps,+}^{\kappa_{j+\frac12}} + \left(f_{\eps,+}^{\kappa_{j+\frac12}} \right)^{4/3} \right)\,.
\end{align*}
As above,
\[ 
f_{\eps,+}^{\kappa_{j+\frac12}} \le (f_\eps-\kappa_{j+\frac14}) \un_{\{ f_\eps-\kappa_{j+\frac14} \ge \kappa_{j+\frac12} - \kappa_{j+\frac14} \}} \le \frac{(f_\eps-\kappa_{j+\frac14})^{4/3}_+}{(\kappa_{j+\frac12}-\kappa_{j+\frac14})^{1/3}}
= 2^{\frac13 (j+3)} (f_\eps-\kappa_{j+\frac14})^{4/3}_+\,,
\]
so that
\begin{equation}
\label{e:flogf}
f_\eps(t,v) \bigg((\ln (f_\eps/{\kappa_{j+\frac12}}))_+ (t,v)  + (\ln (f_\eps/{\kappa_{j+\frac12}}))_+^2 (t,v) \bigg) \le 4\clnt 2^{\frac13 j} \left(f_{\eps,+}^{\kappa_{j+\frac14}}\right)^{\frac43}\,.
\end{equation}
Next we use the fact that $Q^j \subset Q_1$ to write,
\begin{align*}
\iint_{Q^j} \left(  f_\eps(t,\cdot) \ast_v \frac{\un_{B_1}}{|\cdot|}(v)  \right) & f_\eps \bigg( \ln_+ (f_\eps/{\kappa_{j+\frac12}})+\ln_+^2 (f_\eps/{\kappa_{j+\frac12}})\bigg)(t,v)  \dt \dv 
\\
&\le \left\|  f_\eps \ast_v \frac{\un_{B_1}}{|\cdot|} \right\|_{L^5 (Q_1)} \left\|f_\eps ( \ln_+ (f_\eps/{\kappa_{j+\frac12}}) + \ln_+^2 (f_\eps/{\kappa_{j+\frac12}})) \right\|_{L^{5/4}(Q^j)}\,.
\end{align*}
On the other hand
\[
\begin{aligned}
\left\|  f_\eps(t,\cdot) \ast_v \frac{\un_{B_1}}{|\cdot|} \right\|_{L^5 (Q_1)}\le\left\|(f_\eps(t,\cdot)-\tfrac32)_+ \ast_v \frac{\un_{B_1}}{|\cdot|} \right\|_{L^5 (Q_1)}+\tfrac32|Q_1|^{1/5}\left\|\frac1{|\cdot|} \right\|_{L^1 (B_1)}
\\
\le C_{HLS,1}\left\|(f_\eps(t,\cdot)-\tfrac32)_+\right\|_{L^5_tL^{15/13}_v(Q_2)}+\tfrac32|Q_1|^{1/5}\cdot 2\pi&\,,
\end{aligned}
\]
by the first inequality in Lemma \ref{l:source-nl} with $\delta=\rho=1$, with $p=\tfrac{15}{13}$ and $q=5$, and with $\beta=\theta=\tfrac54$. Moreover, Lemma \ref{l:nonlin} implies that
\[
\begin{aligned}
\left\|(f_\eps(t,\cdot)-\tfrac32)_+\right\|_{L^5_tL^{15/13}_v(Q_2)}&\le\Cnonlin\cdot 3^5\cdot 2^3\left( \| (f_\eps(t,\cdot)-1)_+\|_{L^\infty_t L^1_v (I \times \mathcal{B})} + \left\|\nabla_v(\sqrt{f_\eps(t,\cdot)}-1)_+ \right\|_{L^2 (I \times  \mathcal{B})}^2 \right)
\\
&\le 3^5\cdot 2^3\Cnonlin
\end{aligned}
\]
by the assumption in the lemma, since $\eta_{DG}<1$ while $Z_\eps\ge1$. 

Summarizing, we have proved that
\[
\begin{aligned}
\iint_{Q^j} \left(  f_\eps(t,\cdot) \ast_v \frac{\un_{B_1}}{|\cdot|}(v)  \right) & f_\eps \bigg( \ln_+ (f_\eps/{\kappa_{j+\frac12}})+\ln_+^2 (f_\eps/{\kappa_{j+\frac12}})\bigg)(t,v)  \dt \dv 
\\
&\le(2^3\cdot 3^5\cdot C_{HLS,1}\Cnonlin+2^{2/5}\cdot 3^{4/5}\cdot\pi^{6/5})\cdot 4\clnt 2^{\frac13 j}\left\|\left(f_{\eps,+}^{\kappa_{j+\frac14}}\right)^{\frac43} \right\|_{L^{5/4}(Q^j)}
\\
&=(2^3\cdot 3^5\cdot C_{HLS,1}\Cnonlin+2^{2/5}\cdot 3^{4/5}\cdot\pi^{6/5})\cdot 4\clnt 2^{\frac13 j}\left\|f_{\eps,+}^{\kappa_{j+\frac14}}\right\|^{4/3}_{L^{5/3}(Q^j)}\,.
\end{aligned}
\]
Using Lemma~\ref{l:nonlin} as in the case of the term \eqref{e:rhs1} shows that
\begin{equation}\label{e:rhs2-result}
\begin{aligned}
\iint_{Q^j} \left(  f_\eps(t,\cdot) \ast_v \frac{\un_{B_1}}{|\cdot|}(v)  \right) f_\eps \bigg( \ln_+ (f_\eps/{\kappa_{j+\frac12}})+\ln_+^2 (f_\eps/{\kappa_{j+\frac12}})\bigg)(t,v)  \dt \dv 
\\
\le(2^3\cdot 3^5\cdot C_{HLS,1}\Cnonlin+2^{2/5}\cdot 3^{4/5}\cdot\pi^{6/5}) \cdot\clnt\cdot2^{\frac{82}3}\Cnonlin ^{4/3}2^{\frac{17}3 j}U_j^{4/3}&\,.
\end{aligned}
\end{equation}

\paragraph{Convergence.}
We then combine entropy inequalities with \eqref{e:lower-bound-Uk}, \eqref{e:rhs1-result} and \eqref{e:rhs2-result} and get
\[
\begin{aligned}
\cln 2^{-j} U_{j+1} \le &3 C_0 Z_\eps 2^{2j +4}\left(\clnb 2^{\frac{104}3} \Cnonlin^{\frac53} \right) 2^{\frac{22}3 j} U_j^{\frac53}
\\
&+ C_0 2^{2j+4}(2^3\cdot 3^5\cdot C_{HLS,1}\Cnonlin+2^{2/5}\cdot 3^{4/5}\cdot\pi^{6/5}) \cdot\clnt\cdot2^{\frac{82}3}\Cnonlin ^{4/3}2^{\frac{17}3 j}U_j^{4/3}\,.
\end{aligned}
\]

Since $\cln, \clnb, \clnt, C_{HLS,1}$ are absolute, while $\Cnonlin$ depends only on an upper bound on the size of the ball in the velocity space (here one can take this upper bound equal to $2$ since $Q^j\subset Q_2$), and since $C_0$ depends only on 
$|v_0|,m_0,M_0,E_0,H_0$, we conclude that there exists $C_2\ge 1$ depending only on $|v_0|,m_0,M_0,E_0,H_0$ such that
\[
\text{ for each } j \ge 0, \qquad  U_{j+1} \le C_2^{j+1} \left( U_j^{4/3} + Z_\eps U_j^{5/3} \right)\,.
\]    
We know from the assumption \eqref{e:U0} that $U_0$ satisfies
\[ 
U_0 \le \etadg  Z_\eps^{-3/2} 
\]
where $\etadg$ remains to be chosen small enough so that $U_j\to 0$ as $j\to \infty$. 

We remark that the sequence $V_j:= Z_\eps^{3/2} U_j$ satisfies
\[ 
\text{ for each }j \ge 0, \qquad  V_{j+1} \le (C_2)^{j+1} \left( V_j^{4/3} + V_j^{5/3} \right)\,,
\]
and $V_0 \le \etadg$. Indeed, we used that $Z_\eps \ge 1$, so that $Z_\eps^{3/2} U_j^{4/3} \le Z_\eps^2 U_j^{4/3} = V_j^{4/3}$. 

We then choose $\etadg <1/2$  and claim that $V_j\le 1$ for all $j\ge 0$. We argue by contradiction: pick the first $N \ge 0$ such that $V_{N+1} > 1$ and $V_N \le 1$.  We then have $V_N^{5/3}\le V_N^{4/3}$, so that
\[ 
1< V_{N+1} \le 2C_2^{N+1} V_N^{\frac43}\le (2C_2)^{N+1} V_N^{\frac43} \le (2 C_2)^{p(N+1)} V_0^{(4/3)^{N+1}}
\]
with $p (N+1) = \sum_{i=0}^{N} (N+1-i)(4/3)^{i}$. A classical computation (recalled in Appendix, see Lemma~\ref{l:iteration}) shows that
\[
p(N+1)  \le 9\cdot(4/3)^{N+2} = 12\cdot(4/3)^{N+1}.
\]
This implies that 
\[ 
1<V_{N+1} \le (2 C_2)^{12\cdot(4/3)^{N+1}} (V_0)^{(4/3)^{N+1}} \le \left( (2 C_2)^{12} \etadg \right)^{(4/3)^{N+1}}.
\]
Choosing now
\[
\etadg = \min (1/2, (2  C_2)^{-12}/2)
\]
yields a contradiction, since it implies that $V_{N+1} \le (\frac12)^{4/3} <1$. Therefore $V_j \le 1$ for all $j\ge 0$, so that $V_j^{5/3}\le V_j^{4/3}$ and hence $V_{j+1}\le 2C_2^{j+1}V_j^{4/3}$, which implies in turn that
\[ 
V_{j+1} \le (2 C_2)^{p(j+1)} V_0^{(4/3)^{j+1}} \le  \left( (2 C_2)^{12} V_0 \right)^{(4/3)^{j+1}} \le (1/2)^{(4/3)^{j+1}} . 
\]
In particular $V_j \to 0$ (and thus $U_j \to 0$) as $j \to +\infty$. 
\end{proof}

\subsection{Smallness of the local mass}

The goal of this section is to prove that, if the dissipation of a suitable weak solution is sufficiently small after zooming in at scale $\eps$, then  the local mass will be small
enough to apply the De Giorgi lemma (Lemma \ref{l:dg1}) at this scale. Such a result is a consequence of the local mass estimates (Lemma \ref{l:mass}) derived in the previous section. 

\begin{lemma}[Smallness of the local mass]\label{l:eventual-mass}
Let $\ka (r) = r^{\gamma+3}$ with $\gamma \in [-3,-2)$, set $\gammas:=-(\gamma+2)\in(0,1]$, and let $a(z)$ be given by \eqref{Def-a}. Let $f$ be a suitable weak solution of the Landau equation \eqref{e:landau} in $[0,T)\times\R^3$. There exist
a  constant $\lambda \in (0,1/4) \cap \mathbb{Q}$ depending only on $\gamma$ and the mass upper bound $M_0$, and a constant $\eta_D \in (0,1)$ depending only on $\lambda$ and $\etadg$ from Lemma~\ref{l:dg1} such that the following 
implication holds: assuming  
\[
\limsup_{j \to \infty} \left\{ \eps_j^{-\frac{7 \gammas}2} \int_{Q_{2 \eps_j}(t_0,v_0)} \left( \left|\nabla_v \sqrt{f} (t,v)\right|^2 + \int_{\R^3} |F(t,v,w)|^2 \dw \right) \dt \dv \right\} \le \frac\etad{2}
\]
for some $(t_0,v_0) \in (0,T) \times \R^3$ with $\eps_j = \lambda^j$ and $4\eps_j^2 < t_0$, then $f_{\eps_j} (t,v) := \eps_j^{\gamma+5} f(t_0 + \eps_j^2 t, v_0 + \eps_j v)$ satisfies
\[
\limsup_{j \to \infty} \bigg\{Z_{\eps_j}^{3/2} \,  \|f_{\eps_j} \|_{L^\infty_t L^1_v (Q_2)} \bigg\}  \le \frac\etadg3
\]
with $Z_{\eps_j}:=M_0 \eps_j^{-\gammas}$. (One can choose $\etad := (2/9)\lambda^8 M_0^{-3/2} \etadg$.)
\end{lemma}

We recall that  $F (t,v,w) = \sqrt{a (v-w)} \sqrt{ f (t,v) f(t,w)} \left( \nabla_v (\ln f)(t,v) - \nabla_w (\ln f)(t,w)\right)$ and also that $\sqrt{a} (z) = |z|^{-1/2}\ka (|z|)^{\frac12}\Pi (z)$. 

\begin{proof}
Without loss of generality, we reduce our attention to the case $(t_0,v_0)= (0,0)$ by considering the shifted distribution function $\tilde{f}(t,v) = f(t_0+v,v_0+v)$. We know from the assumption that there exists $K \in \mathbb{N}$ such that, for all $k \ge K$,
\[ 
\int_{Q_{2 \eps_j}}\left( \left|\nabla_v \sqrt{f} (t,v)\right|^2 + \int_{\R^3} |F(t,v,w)|^2 \dw \right) \dt \dv \le \eta_D \, \eps_j^{\frac{7 \gammas}2}\,.
\]
This means that the scaled function $f^j(t,v) := \eps_j^{\gamma+5} f(\eps_j^2 t, \eps_j v)$ satisfies 
\[ 
\int_{Q_2} \left( \left|\nabla_v \sqrt{f^j} (t,v)\right|^2 + \int_{\R^3} |F^j (t,v,w)|^2 \dw \right) \dt \dv \le \eta_D \, \eps_j^{\frac{5\gammas}2}\,,
\]
where the function $F^j (t,v,w)$ is defined as $\sqrt{a (v-w) f^j (t,v) f^j(t,w)} \left( \nabla_v (\ln f^j)(t,v) - \nabla_w (\ln f^j)(t,w)\right)$.

Observe that 
\[     
\frac1{C_1 \lambda^{3-\gammas}} \| f^{j+1} \|_{L^\infty_t L^1_v (Q_2)} = \frac1{C_1 \lambda^3} \| f^j \|_{L^\infty_t L^1_v (Q_{2\lambda})}\,.
\]
Lemma~\ref{l:mass} implies that there exists a constant $C_1>0$ (depending only on $M_0$) such that, for all $\lambda \in (0,1/4)$,
\begin{align*}
\frac1{C_1 \lambda^{3-\gammas}} \| f^{j+1} \|_{L^\infty_t L^1_v (Q_2)}  & \le \left(1+\lambda^{-4} \sqrt{\etad}\eps_j^\frac{5\gammas}4\right) \| f^j \|_{L^\infty_t L^1_v (Q_2)} + \lambda^{-8} (\etad\eps_j^\frac{5\gammas}2 +1) \etad \eps_j^\frac{5\gammas}2 
+ \lambda^{-8} \eps_j^{\gamma+2}\etad\eps_j^\frac{5\gammas}2 
\\
& \le(1+\lambda^{-4} \sqrt{\etad} ) \| f^j \|_{L^\infty_t L^1_v (Q_2)} + 3\lambda^{-8} \etad \eps_j^\frac{3\gammas}2\,,
\end{align*}
where we used that $\etad \le 1$, and that $\eps_j \le 1$ to get the second line. 

We pick $\lambda \in (0,1/4)\cap\mathbb{Q}$ such that $C_1 \lambda^{3-\frac{5\gammas}2}\le\frac14$. Remark that $\lambda$ depends only on $M_0$ and $\gammas$ since $C_1$ depends only on $M_0$. Then we pick $\etad\le\lambda^8$. 
The previous inequality yields
\[
\frac1{\lambda^{\frac{3\gammas}2}}  \| f^{j+1} \|_{L^\infty_t L^1_v (Q_2)}  \le  \frac12 \| f^j \|_{L^\infty_t L^1_v (Q_2)} +  \frac{3 \etad}{4 \lambda^8}  \eps_j^{\frac{3 \gammas}2}\,.
\]
We used that $\etad \le 1$ and $\eps_j \le 1$. Recalling that $\eps_j = \lambda^j$, we can rewrite the previous inequality as,
\[
\eps_{j+1}^{-\frac{3\gammas}2}  \| f^{j+1} \|_{L^\infty_t L^1_v (Q_2)} \le \frac12  \eps_j^{-\frac{3\gammas}2} \| f^j \|_{L^\infty_t L^1_v (Q_2)} +  \frac{3 \etad}{4 \lambda^8} \,,
\]
which implies in turn that
\[  
\limsup_{j \to \infty} \bigg\{ \eps_j^{-\frac{3\gammas}2} \,  \|f_{\eps_j} \|_{L^\infty_t L^1_v (Q_2)} \bigg\}  \le  \frac{3\etad}{2 \lambda^8}\,.
\]
We then pick $\etad \le \frac{2\lambda^8}{3} (\etadg/3) M_0^{-3/2}$ and get the sought estimate.  
\end{proof}

We deduce from the previous lemma that the De Giorgi lemma can be applied at scale $\eps_j = \lambda^j$ for $j$ large enough if the entropy dissipation is small enough. 

\begin{lemma}[Improved De Giorgi lemma]\label{l:improved-dg}
Let $\ka (r) :=r^{\gamma+3}$ with $\gamma \in [-3,-2)$, set $\gammas := -(2+\gamma)\in (0,1]$, and let $a(z)$ be given by \eqref{Def-a}. Let $f$ be a suitable weak solution of the Landau equation \eqref{e:landau} in $[0,T)\times\R^3$ with initial datum 
$\fin$ satisfying \eqref{ini-bnd-bis}. Given $(t_0,v_0)\in (0,T)\times\R^3$, there exists $\etadgp>0$ depending only on $m_0,M_0,E_0, H_0, |v_0|$ such that 
\[
\begin{aligned}
\limsup_{j \to \infty} \left\{  \eps_j^{- \frac{7\gammas}2} \int_{Q_{2 \eps_j}(t_0,v_0)} \left( \left|\nabla_v \sqrt{f} (t,v)\right|^2 + \int_{\R^3} |F(t,v,w)|^2 \dw \right) \dt \dv \right\} \le 2 \etadgp
\\
\implies f\le 2 \eps_{j_0}^{-\gamma-5}\text{ a.e. in }Q_{\eps_{j_0}/2} (t_0,v_0)
\end{aligned}
\]
for some positive integer $j_0$. 
\end{lemma}

\begin{proof}
Let $\etadgp := \frac14 \min (\etad , \etadg M_0^{-3/2})$ with $\etadg$ given by Lemma~\ref{l:dg1} and $\etad$ given by Lemma~\ref{l:eventual-mass}. The $\eps_j$-scaled suitable weak solution $f_j$ defined by 
$f_j (t,v) := \eps_j^{\gamma+5} f(t_0 + \eps_j^2 t, v_0 + \eps_j v)$ satisfies, for $(t,v) \in Q_1(0,0)$, 
\begin{align*}
f_j(t,\cdot) \ast \frac{\ka ( |\cdot|) \un_{B_1^c}}{|\cdot|} (v) &= \int_{|w| \ge 1} f_j(t,v-w)\frac{\ka (|w|)}{|w|} \dw 
\\
&\le  k(1) \int_{|w| \ge 1} \eps_j^{\gamma+5} f (t_0+\eps_j^2 t, v_0+\eps_j v - \eps_j w) \dw 
\\
& = \eps_j^{\gamma+2} \int_{|\bar w| \ge \eps_j}  f (t_0+\eps_j^2 t,v_0+ \eps_j v - \bar w)  \dd \bar w \le 2M_0 \eps_j^{-\gammas}
\end{align*}
according to \eqref{e:hydro}. Therefore
\[ 
\esssup_{(t,v)\in Q_1(0,0)} \left(  f_j(t,\cdot) \ast_v \frac{\ka (|\cdot|) \un_{B_1^c}}{|\cdot|}(v)\right) \le Z_j\,,\qquad \text{ with }Z_j := M_0 \eps_j^{-\gammas}\,.
\]

We know from Lemma~\ref{l:eventual-mass} that $f_j$ also satisfies,
\[
\limsup_{j \to \infty} \left\{  Z_j^{\frac32} \|f_j\|_{L^\infty_t L^1_v (Q_2)} \right\}  \le \frac{\etadg}{3}
\]
with $\eps_j := \lambda^j \in (0,1)$. In particular, for $j$ large enough, we have
\[ 
Z_j^{\frac32}  \|f_j\|_{L^\infty_t L^1_v (Q_2)}    \le \left( \tfrac{1}{3} +\tfrac1{12} \right) \etadg .
\]

Moreover, the assumption of the lemma implies that, for $j$ large enough, 
\[
\eps_j^{-\frac{3\gammas}2}  \int_{Q_2}   \left|\nabla_v \sqrt{f_j} (t,v)\right|^2 \dt \dv \le \eps_j^{-\frac{5\gammas}2}  \int_{Q_2}   \left|\nabla_v \sqrt{f_j} (t,v)\right|^2 \dt \dv \le \left(2+\tfrac13 \right) \etadgp\le \left(\tfrac12 + \tfrac1{12}\right) \etadg M_0^{-3/2}
\]
since $\eps_j\in(0,1)$. In particular,
\[ 
\|f_j\|_{L^\infty_t L^1_v (Q_2)} + \int_{Q_2}   \left|\nabla_v \sqrt{f_j} (t,v)\right|^2 \dt \dv \le \left( \tfrac{1}{3}+\tfrac12 +\tfrac2{12} \right) \etadg Z_j^{-3/2} = \etadg Z_j^{-3/2}.
\]
We deduce from Lemma~\ref{l:dg1} that $f_j \le 2$ in $Q_{1/2}(0,0)$, that is to say $f \le 2 \eps_j^{-\gamma-5}$ in $Q_{\eps_j/2} (t_0,v_0)$. 
\end{proof}

\subsection{Partial regularity}
\label{sec:partial}

Our first main result in this paper is the following partial regularity theorem.

\smallskip
\noindent
\textbf{Theorem A'.} 
{\it 
Let $f$ be a global suitable weak solution of the Landau equation \eqref{e:landau} for $\gamma\in[-3,-2)$, with initial data \eqref{e:ic} satisfying \eqref{ini-bnd}. Then $\Pfinal(\Sing[f])=0$ with $\ms:=\tfrac72|2+\gamma|$.
}

\smallskip
In view of Remark \ref{rmk:VillaniSuitable}, Theorem A' implies Theorem A of section \ref{ss:MainRes}.

\smallskip
The proof of Theorem A' uses a Vitali-type covering lemma for parabolic cylinders, for which we refer to \cite[Lemma~6.1]{ckn}. To state it, we introduce the following notation: for any parabolic cylinder $Q_r(t,v)$, we set
\[
5 \cdot Q_r(t,v):=Q_{5r}(t,v)\,.
\]

\begin{lemma}[Vitali covering lemma]\label{l:vitali}
Let $X$ be a bounded set of $\R\times \R^3$. Assume that $X \subset \bigcup_{\alpha} Q^\alpha$ for some collection $\{Q^\alpha\}_\alpha$ of parabolic cylinders with sizes $r_\alpha\in(0,1)$. Then there exists a countable subcollection
$\{ Q^{\alpha_j} \}_{j=1}^{\infty}$ of $\{Q^\alpha\}_\alpha$ whose elements are pairwise disjoint and such that  
$$
X \subset \bigcup_{j=1}^{+\infty} 5 \cdot Q^{\alpha_j}\,.
$$
\end{lemma}

\begin{proof}[Proof of Theorem A']
Let $R>0$ and $\delta \in (0,\bar\delta]\subset (0,1)$. Let $\mathcal{V} \subset(0,T) \times B_R$ be a neighborhood of the singular set $\Sing\left[f\big|_{(0,T) \times B_R}\right]$. We know from Lemma~\ref{l:improved-dg} that, for each point  
$(\tau,\zeta) \in\Sing\left[f\big|_{(0,T) \times B_R}\right]$, there exists $r\equiv r (\tau,\zeta) \in (0,\delta)$, such that
\begin{equation}\label{e:qrj}
\int_{Q_r(\tau,\zeta)}\left( \left|\nabla_v \sqrt{f} (t,v)\right|^2 + \int_{\R^3} |F(t,v,w)|^2 \dw \right) \dt \dv > \etadgp  (r/2)^{\ms}
\end{equation}
where $\ms:=  \frac{7 \gammas}2 =  \frac72 |2+\gamma|$. 

We can also reduce $r>0$ such that  $5r<\delta$ and $Q_r (\tau,\zeta) \subset \mathcal{V}$ for all $(\tau,\zeta)$. Hence $\Sing\left[f\big|_{(0,T) \times B_R}\right]$ is covered by   cylinders $Q_r(\tau,\zeta) \subset \mathcal{V}$ 
such that $5 r < \delta$ and \eqref{e:qrj} holds true. 

The Vitali-type covering lemma (Lemma~\ref{l:vitali}) implies that there exists two sequences $(t_j,v_j) \in (0,T) \times B_R$ and $r_j \in (0,\delta)$  such that 
\[
j \neq k\implies Q_{r_j}(t_j,v_j) \cap Q_{r_k} (t_k,x_k) = \varnothing\,,
\] 
such that  
\begin{equation}\label{CovSing}
Q_{r_j} (t_j,v_j) \subset \mathcal{V}\quad\text{ and }\quad\Sing\left[f\big|_{(0,T) \times B_R}\right]\subset \bigcup_{j=1}^{\infty} 5\cdot Q_{r_j} (t_j,v_j)\,,
\end{equation}
and such that \eqref{e:qrj} holds true for the parabolic cylinder $Q_{r_j}(t_j,v_j)$.  Because the cylinders $Q^{\alpha_j} := Q_{r_j} (t_j,v_j)$ are pairwise disjoint, we can write
\begin{equation}\label{Hmeas}
\begin{aligned}
\sum_{j=1}^\infty r_j^{\ms} & \le \frac{2^{\ms}}{\etadgp}   \sum_{j=1}^\infty   \int_{Q^{\alpha_j}}\left( \left|\nabla_v \sqrt{f} (t,v)\right|^2 + \int_{\R^3} |F(t,v,w)|^2 \dw \right) \dt \dv 
\\
& \le \frac{2^{\ms}}{\etadgp}      \int_{\mathcal{V}}   \left( \left|\nabla_v \sqrt{f} (t,v)\right|^2 + \int_{\R^3} |F(t,v,w)|^2 \dw \right) \dt \dv\,.
\end{aligned}
\end{equation}
Since $\delta$ is arbitrarily small and $\ms\le\tfrac72<5$, \eqref{CovSing} and \eqref{Hmeas} imply that $\Sing\left[f\big|_{(0,T) \times B_R}\right]$ is Lebesgue-negligible in $(0,+\infty) \times \R^3$. The inclusion \eqref{CovSing} and the definition 
of $\Pfinal$ with the inequality \eqref{Hmeas} imply that 
\[
\Pfinal \left(\Sing\left[f\big|_{(0,T) \times B_R}\right]\right) \le \sum_{j=1}^\infty (5r_j)^{\ms} \le \frac{10^{\ms}}{\etadgp}\int_{\mathcal{V}}   \left( \left|\nabla_v \sqrt{f} (t,v)\right|^2 + \int_{\R^3} |F(t,v,w)|^2 \dw \right) \dt \dv\,.
\]
Since $\mathcal{V}$ is an arbitrary neighborhood of $\Sing\left[f\big|_{(0,T) \times B_R}\right]$, and  since $\Sing\left[f\big|_{(0,T) \times B_R}\right]$  is Lebesgue-negligible, the last right-hand side of the inequality above can be made arbitrarily 
small, so that $\Pfinal\left(\Sing\left[f\big|_{(0,T) \times B_R}\right]\right)=0$.  Since $R,T>0$ are arbitrary, we conclude that $\Pfinal (\Sing\left[f\right])=0$.
\end{proof}


\section{The axisymmetric case}\label{s:axi-sym}


Our second main result in this paper is the following local regularity result for axisymmetric and radial suitable weak solutions of the Landau equation. The present section is devoted to the study of such solutions. In this section again, we restrict our 
attention to the case $k(r)=r^{\gamma+3}$ with $\gamma\in[-3,-2)$. 

\smallskip
\noindent
\textbf{Theorem B'.} 
{\it 
Let $f$ be a global suitable weak solution of the Landau equation \eqref{e:landau} for $\gamma\in[-3,-2)$, with initial data \eqref{e:ic} satisfying \eqref{ini-bnd}. Assume that $f$ is axisymmetric, viz. there exists a measurable function $\mathcal F$ on 
$(0,+\infty)\times(0,+\infty)\times\R$ and $(\bar v,\omega)\in\R^3\times\mathbb{S}^2$ such that
$$
f(t,v)=\mathcal F(t,|(v-\bar v)\times\omega|,(v-\bar v)\cdot\omega)\quad\text{ for a.e. }(t,v)\in(0,+\infty)\times\R^3.
$$
Then $\Sing[f]\subset(0,+\infty)\times(\bar v+\R\omega)$.
}

\smallskip
The case of radial, suitable weak solutions is a straightforward consequence of Theorem B'.

\smallskip
\noindent
\textbf{Corollary C'.} 
{\it 
Let $f$ be a global suitable weak solution of the Landau equation \eqref{e:landau} for $\gamma\in[-3,-2)$, with initial data \eqref{e:ic} satisfying \eqref{ini-bnd}. Assume that $f$ is radial, viz. there exists $\bar v\in\R^3$ and a measurable function 
$\mathfrak F$ on $(0,+\infty)\times(0,+\infty)$ such that
$$
f(t,v)=\mathfrak F (t,|v-\bar v|)\quad\text{ for a.e. }(t,v)\in(0,+\infty)\times\R^3.
$$
Then $\Sing[f]\subset(0,+\infty)\times\{\bar v\}$.
}

\begin{proof}[Proof of Corollary C']
Indeed, for all $\omega\in\mathbb S^2$, one has 
\[
\mathfrak F(t,|v-\bar v|)=\mathcal F(t,|(v-\bar v)\times\omega|,(v-\bar v)\cdot\omega) \qquad \text{ with }\mathcal F(t,r,z):=\mathfrak F\left(t,\sqrt{r^2+z^2}\right)\,.
\]
In other words, any radial suitable weak solution $f$ around $\bar v$ is axisymmetric around the axis $\bar v+\mathbb{R}\omega$ for all $\omega\in\mathbb S^2$. By Theorem B',
\[
\mathbf S[f]\subset\bigcap_{\omega\in\mathbb S^2}(0,+\infty)\times(\bar v+\R\omega)=(0,+\infty)\times\{\bar v\}\,.
\]
\end{proof}

\smallskip
As explained above, Theorem B' and Corollary C' imply Theorem B and Corollary C of section \ref{ss:MainRes} in view of Remark \ref{rmk:VillaniSuitable}.

\smallskip
In order to prove that the solution $f$ considered in Theorem B' is locally bounded, except along the symmetry axis, we first prove that long-range interactions are bounded away from the symmetry axis. We then prove that axisymmetric solutions 
enjoy better integrability away from the axis.

\subsection{Long-range interactions}

\begin{lemma}[Long-range interactions]\label{l:lri}
Assume that $f:(0,T) \times \R^3 \to[0,+\infty)$ measurable is axisymmetric: there exists $\mathcal F:(0,T) \times \R^2 \to[0,+\infty)$ such that $f(t,v) = \mathcal F (t,|v'|,v_3)$ for a.e. $v = (v',v_3) \in \R^2 \times \R$. Given $(t_0,v_0) \in (0,T)\times\R^3$ and
$\eps \in\left(0,\rho_0/2\wedge\sqrt{t_0}\right)$ with $\rho_0 = |v_0'|$, the scaled function 
$$
f_\eps (t,v) = \eps^{\gamma+5} f (t_0+\eps^2 t, v_0+\eps v)
$$
satisfies 
\[ 
\esssup_{Q_1} \left( \eps^{-\gamma-3} f_\eps \ast_v \frac1{|\cdot|} \right) \le \frac{C_\ast}{\rho_0}\sup_{t\ge 0}\left( \int_{\R^3} f(t,w) \dw +  \int_{\R^3} f \lnp f(t,w) \dw \right) + C_\ast \rho_0^2
\]
for some absolute constant $C_\ast>0$. 
\end{lemma}

\begin{proof}
Let $(t,  v) = (t_0,v_0) + (\eps^2 \bar t, \eps \bar v)$ for $(\bar t,\bar v) \in Q_1(0,0)$ and $v_0= (v_0',v_3^0)$ and $\rho_0 = |v_0'|$. Then,
\begin{align*}
\eps^{-\gamma-3} f_\eps(\bar t,\cdot) \ast_v \frac1{|\cdot|}(\bar v) &= \eps^{-\gamma-3}  \int_{\R^3} f_\eps (\bar t,\bar w) \frac{\dd\bar w}{|\bar v-\bar w|}  = \int_{\R^3} f (t, w) \frac{\dd w}{|v -  w|} 
\\
& \le  \int_{B_{\sigma_0}(v)} f (t, w) \frac{\dw}{| v -  w|}  + \frac1{\sigma_0} \int_{\R^3} f(t,w) \dw
\end{align*}
for some $\sigma_0$ to be chosen later. The second equality above follows from integrating by substitution with $w=v_0+\eps\bar w$. We are left with the task of estimating the first term on the last right-hand side above. 

We use the symmetry assumption on $f$ to do so. Write $v = (\rho_v , 0, v_3)$ and $ w = (\rho \cos \theta, \rho \sin \theta, w_3)$ and
\begin{align*}
| v -  w|^2 &= |v' - w'|^2  + | v_3-w_3|^2 = (|v'|-|w'|)^2 + | v_3-w_3|^2 + 2(|v'||w'|-v' \cdot w') 
\\
& = (\rho_v -\rho)^2 +  | v_3-w_3|^2  + 2 \rho_v \rho (1-\cos \theta)\,.
\end{align*}

Since $\bar v \in B_1$, we  have $\rho_v = |v_0'+ \eps \bar v'| \ge \rho_0 -\eps \ge \frac12 \rho_0$ for $\eps \le \frac12 \rho_0$. As far as $\rho$ is concerned, we have $\rho = |w'| \ge |v'|-|v'-w'| \ge \rho_v - \sigma_0 \ge \frac14 \rho_0$ for
$\sigma_0 \le \frac14 \rho_0$.  Similarly, \( \rho_v \le \frac32 \rho_0\) and \(\rho \le 2 \rho_0\) under the same condition. We have thus proved that
\begin{equation}\label{e:rho-lower}
\sigma_0 \le \frac14 \rho_0\implies\rho_v \in [\rho_0/2,3\rho_0/2] \quad \text{ and } \quad \rho \in [\rho_0/4,2\rho_0]\,.
\end{equation}

We next compute
\begin{multline}\label{e:pressure-estimate}
\int_{B_{\sigma_0}( v)} f ( t, w) \frac\dw{| v -  w|} = \int_{\R^3} f ( t, w) \frac{\un_{|v-w|^2\le\sigma_0^2}}{| v -  w|} \dw
\\
= 2 \int_{(0,+\infty)\times \R} \mathcal F (t, \rho,w_3) \left\{ \int_0^\pi\frac{\un_{\{|(\rho_v,v_3) -(\rho,w_3)|^2 + 2 \rho_v \rho (1-\cos \theta) \le \sigma^2_0\}}\dd \theta}{\sqrt{|(\rho_v,v_3) -(\rho,w_3)|^2 + 2 \rho_v \rho (1-\cos \theta) }} \right\} \rho \dd \rho \dw_3\,.
\end{multline}
Let us compute the inner integral. Let $A := |(\rho_v,v_3) -(\rho,w_3)|$ and $B :=\sqrt{\rho_v \rho}$, we have
\begin{align*}
\int_0^\pi \frac{\un_{\{|(\rho_v,v_3) -(\rho,w_3)|^2 + 2 \rho_v \rho (1-\cos \theta) \le \sigma^2_0\}}\dd \theta}{\sqrt{|(\rho_v,v_3) -(\rho,w_3)|^2 + 2 \rho_v \rho (1-\cos \theta) }} &=  \int_0^\pi \frac{\un_{\{A^2 + 2 B^2 (1-\cos \theta) \le \sigma_0^2\}}}{\sqrt{A^2 + 2 B^2 (1-\cos \theta) }}  \dd \theta 
\\
& = \frac1A \int_0^{\theta_\ast} \frac{ \dd \theta }{\sqrt{1 +  (2B/A)^2 \sin^2 (\theta/2) }} 
\end{align*}
with $\theta_\ast\in[0,\pi]$ such that $1 + (2B/A)^2 \sin^2(\theta_\ast/2) = (\sigma_0/A)^2$. Since $\sin\frac\theta{2}\ge\frac\theta\pi$ for all $\theta\in[0,\pi]$, one has
\begin{equation}\label{key<}
\begin{aligned}
\frac1A\int_0^{\theta_*}\frac{d\theta}{\sqrt{1+(\frac{2B}A)^2\sin^2\frac\theta{2}}}&\le\frac1A\int_0^{\theta_*}\frac{d\theta}{\sqrt{1+(\frac{2B\theta}{A\pi})^2}}\le\frac1A\int_0^{\tilde\theta_*}\frac{d\theta}{\sqrt{1+(\frac{2B\theta}{A\pi})^2}}
\\
&=\tfrac\pi{2B}\text{Arsinh}\left(\tfrac{2B\tilde\theta_*}{A\pi}\right)\le\tfrac\pi{2B}\ln\left(2\sqrt{\left(\tfrac{2B\tilde\theta_*}{A\pi}\right)^2+1}\right)=\tfrac\pi{2B}\left(\ln2+\ln\tfrac{\sigma_0}{A}\right)\,,
\end{aligned}
\end{equation}
where $\tilde\theta_*$ is defined by
\[
1+\left(\tfrac{2B\tilde\theta_*}{A\pi}\right)^2=\left(\tfrac{\sigma_0}{A}\right)^2=1+\left(\tfrac{2B}A\right)^2\sin^2\tfrac{\theta_*}2\ge 1+\left(\tfrac{2B\theta_*}{A\pi}\right)^2\,,\qquad\text{ so that }\theta_*\le\tilde\theta_*\,,
\]
which justifies the second inequality in \eqref{key<}. The third inequality in \eqref{key<} follows from the identity $\text{Arsinh}X=\ln\left(X+\sqrt{X^2+1}\right)$. Since we know from \eqref{e:rho-lower} that $B\ge\rho_0/2\sqrt{2}$, we conclude that
\begin{equation}\label{Ker<}
\int_0^\pi \frac{\un_{\{|(\rho_v,v_3) -(\rho,w_3)|^2 + 2 \rho_v \rho (1-\cos \theta) \le \sigma^2_0\}}\dd \theta}{\sqrt{|(\rho_v,v_3) -(\rho,w_3)|^2 + 2 \rho_v \rho (1-\cos \theta) }}
\le\tfrac{\pi\sqrt2}{\rho_0}\left(\ln2+\ln_+\tfrac{\sigma_0}{|(\rho_v,v_3)-(\rho,w_3)|}\right)\,.
\end{equation}
Inserting this bound in \eqref{e:pressure-estimate}, one finds that
\[
\begin{aligned}
\int_{B_{\sigma_0}(v)}f(t,w)\frac{dw}{|v-w|}&\le\tfrac{2\sqrt{2}\pi}{\rho_0}\int_{(0,+\infty)\times\mathbb R}\mathcal F(t,\rho,w_3)\left(\ln2+\ln_+\tfrac{\sigma_0}{|(\rho_v,v_3)-(\rho,w_3)|}\right)\rho d\rho dw_3
\\
& \le\tfrac{\sqrt{2}\ln 2}{\rho_0}\int_{\mathbb R^3}f(t,w)dw+\tfrac{2\sqrt{2}\pi}{\rho_0}\int_{(0,+\infty)\times\mathbb R}\mathcal F(t,\rho,w_3)\ln_+\left(\tfrac{\sigma_0}{|(\rho_v,v_3)-(\rho,w_3)|}\right)\rho d\rho dw_3\,.
\end{aligned}
\]

We now use Fenchel's inequality from convex analysis: for all $p,q >0$, 
\[ 
pq \le h_0 (p) + (h_0)^* (q) \,,
\]
where $h_0$ denotes the convex function $h_0(p):=p \ln p - (p-1)$ and $(h_0)^* (q) = e^q-1$ its Legendre-Fenchel conjugate.
In particular,
\[ 
pq \le p \lnp p + e^q\,.
\]
This implies that
\[  
\mathcal F (t, \rho,w_3) \ln\frac{\sigma_0}{|(\rho_v,v_3) -(\rho,w_3)|} \le (\mathcal F \lnp \mathcal F)(t, \rho,w_3)   +  \frac{\sigma_0}{|(\rho_v,v_3) -(\rho,w_3)|}\,. 
\]
With such an inequality at hand, we arrive at the following estimate,
\[
\int_{\R^3}  \mathcal F(t, \rho,w_3) \lnp\left(\frac{\sigma_0}{|(\rho_v,v_3) -(\rho,w_3)|}\right)\rho \dd \rho \dw_3  \le \tfrac1{2\pi} \int_{\R^3} f \lnp f(t,w) \dw +  4 \pi \rho_0 \sigma_0^2 .
\]
We used $\rho \le 2 \rho_0$ to estimate the second integral. Gathering together the previous estimates yields
\[
\int _{\R^3}f ( t, w) \frac{\dw}{| v -  w|}  \le \frac1{\sigma_0} \int_{\R^3} f(t,w) \dw +\frac{2\sqrt{2}\pi}{\rho_0} \left( \tfrac1{2\pi}\int_{\R^3} f(t,w) \dw + \frac1{2\pi} \int_{\R^3} f \lnp f(t,w) \dw +4 \pi \rho_0\sigma_0^2 \right)\,,
\]
since $\ln 2<1$. Choosing the constant $\sigma_0 = \rho_0/2\sqrt{2}\pi$ concludes the proof. (One can take $C_\ast=8\sqrt{2}\pi^2$.)
\end{proof}

\subsection{Improved integrability}

\begin{lemma}[Improved integrability]\label{l:better-int}
Let $f \in L^\infty_t L^1_v ((0,T) \times B_R)$ and $\nabla_v f \in L^2_t L^1_v ((0,T) \times B_R)$ for each $R>0$. Assume that $f$ is axisymmetric, that is to say of the form $f(t,(v',v_3)) = \mathcal F(t,|v'|,v_3)$ for some measurable $\mathcal F$ defined a.e.
on $(0,T)\times[0,+\infty)\times\R$. Then for all $t_0 \in (0,T)$, all $V_0 = (V_1^0,V_2^0)$, and all $r\in\left(0,\sqrt{t_0}\right)$ such that $V_1^0 > r + \rho_0$, one has 
\[
\mathcal F \in L^2_{t,v}(\mathcal{Q}_r (t_0,V_0))\,,
\]
where $\mathcal{Q}_r (t_0,V_0)\subset(0,+\infty)\times\R^2$ is the two-dimensional cylinder of size $r$ centered at $(t_0,V_0)$.
\end{lemma}

\begin{proof}
Observe first that $|\nabla_{v'} f (t,(v',v_3))| = |\partial_{V_1} \mathcal F(t,|v'|,v_3)|$. In particular,
\[ 
|\nabla_v f (t,(v',v_3))|^2 = |\partial_{V_1} \mathcal F (t,|v'|,v_3)|^2 + |\partial_{v_3} \mathcal F (t,|v'|,v_3)|^2 = |\nabla_V \mathcal F (t,|v'|,v_3)|^2\,. 
\]
If $V_1^0-r\ge \rho_0>0$, then we can write,
\begin{align*}
2 \pi  \int_{|V-V_0| \le r}  |\nabla_V \mathcal F (t,V)| \dd V & = \int_{-\pi}^\pi \int_{||v'|-V_1^0|^2 + |v_3-V_2^0|^2 \le r^2}  |\nabla_v f (t,(v',v_3))| \dd |v'|  \dv_3 \dd \theta  
\\
\intertext{using the inequality $|v'|\ge V_1^0-r>\rho_0$,}
&\le   \frac1{\rho_0} \int_{-\pi}^\pi \int_{||v'|-V_1^0| \le r,  |v_3-V_2^0| \le r}  |\nabla_v f (t,(v',v_3))| |v'| \dd |v'|  \dv_3 \dd \theta
\\
&\le   \frac1{\rho_0} \int_{-\pi}^\pi \int_{|v'| \le V_1^0+r,   |v_3-V_2^0| \le r}  |\nabla_v f (t,(v',v_3))| |v'| \dd |v'|  \dv_3 \dd \theta 
\\
\intertext{using the formula $\dv=|v'| \dd |v'|  \dv_3 \dd \theta$ for the volume element in cylindrical coordinates,}
& \le \frac1{\rho_0} \int_{|v-(0,V_2^0)|^2 \le (V_1^0 +r)^2 + r^2 } |\nabla_v f (t,v)| \dv\,. 
\end{align*}
Hence, if $r_1$ denotes $\sqrt{(V_1^0+r)^2 + r^2}$, we have
\[ 
\| \nabla_V \mathcal F \|^2_{L^2_t L^1_v (\mathcal{Q}_r (t_0,V_0))} \le (2\pi\rho_0)^{-2} \int_{t_0-r^2}^{t_0} \| \nabla_v f(t)\|_{L^1(B_{r_1} (0,V_2^0))}^2 \dt\,. 
\]
Similarly
\[
\| \mathcal F \|^2_{L^2_t L^1_v (\mathcal{Q}_r (t_0,V_0))} \le (2\pi\rho_0)^{-2} \int_{t_0-r^2}^{t_0} \| f(t,\cdot)\|_{L^1(B_{r_1} (0,V_2^0))}^2 \dt\,. 
\]
In particular, $\mathcal F \in L^2 (\mathcal{Q}_r (t_0,V_0))$ by the Sobolev embedding $W^{1,1}(B_R)\subset L^2(B_R)$ in space dimension $2$. 
\end{proof}

\subsection{Regularity away from the symmetry axis}

This subsection is devoted to the proof of Theorem B'.

\begin{proof}[Proof of Theorem B']
Let $f$ be a suitable weak solution of \eqref{e:landau}. We can assume without loss of generality that the direction of the symmetry axis is the third vector  $e_3 = (0,0,1)$ of the canonical basis of $\R^3$.  Let $t_0 \in (0,T)$ and $v_0=(v_0',v^0_3)$ 
be such that $v'_0 \neq 0$. We claim that, for $\eps \in\left(0,\rho_0/6\wedge\tfrac13\sqrt{t_0}\right)$,
\[ 
\esssup_{Q_1} \left( f_\eps \ast_v \frac{\ka ( |\cdot|) \un_{B_1^c}}{|\cdot|} \right) \le Z_0 
\]
for some constant $Z_0\ge 1$ depending only on $\rho_0= |v_0'|$ and the macroscopic bounds $M_0$, $E_0$, $H_0$. Indeed, for $|z| \ge 1$ and $\eps \in (0,1)$, we have 
\[ 
\eps^{\gamma+3} \frac{\ka ( |z|)}{|z|} = \eps^{\gamma+3} \frac1{|z|^{\gammas}} \le (1-\gammas)\eps+\gammas\frac1{|z|}\,,
\]
where we recall that $\gammas=|\gamma+2|=-(\gamma+2)\ge 0$. This inequality is trivial if $\gamma=-3$; if $\gamma\in(-3,-2)$, apply Young's inequality $ab\le\tfrac1\alpha a^\alpha+\tfrac1\beta b^\beta$ with $a=\eps^{\gamma+3}$ and $\alpha=\frac1{1-\gamma_*}=\frac1{\gamma+3}$, while $b=\frac1{|z|}$ and 
$\beta=\frac1\gammas$, so that $\frac1\alpha+\frac1\beta=1$. Hence
\[
\begin{aligned}
\esssup_{Q_1} \left( f_\eps\ast_v \frac{\ka ( |\cdot|)\un_{B_1^c}}{|\cdot|} \right)\le(1-\gammas)\sup_{-1<\bar t\le 0}\eps^{-\gamma-2}\int_{\R^3}f_\eps(\bar t,v)\dv + \gammas\esssup_{Q_1} \left( \eps^{-\gamma-3}f_\eps\ast_v\frac1{|\cdot|} \right)
\\
\le(1-\gammas)\sup_{0\le t_0-t<\eps^2}\int_{\R^3}f(t,w)\dw+\frac{C_*\gammas}{\rho_0}\sup_{t\ge 0}\int_{\R^3}f(1+\ln_+f)(t,w)\dw+C_*\gammas\rho_0^2
\\
\le2(1-\gammas)M_0+\frac{C_*\gammas}{\rho_0}(2M_0+\bar H_0)+C_*\gammas\rho_0^2&\,.
\end{aligned}
\]
We have used Lemma~\ref{l:lri} to prove the second inequality above, and \eqref{e:hydro} and \eqref{e:dpsin} to find upper bounds for $\int_{\R^3} f(t,w) \dw$ and $\int_{\R^3} f\lnp f(t,w) \dw$. Hence we can choose $Z_0:=1+2(1-\gammas)M_0+\frac{C_*\gammas}{\rho_0}(2M_0+\bar H_0)+C_*\gammas\rho_0^2$.

Thanks to De Giorgi's lemma (Lemma~\ref{l:dg1}), it is sufficient to prove that 
\begin{equation}\label{e:cond-dg1}
\esssup_{t \in (-4,0]} \int_{B_2} (f_\eps-1)_+ (t,v) \dv + \int_{Q_2} \left|\nabla_v \sqrt{f_\eps}\right|^2 \un_{\{ f_\eps \ge 1\}} \dt \dv \le \etadg \; Z_0^{-3/2}
\end{equation}
to conclude that $f$ is bounded in some neighborhood of $(t_0,v_0)$. 

In order to get such an estimate, we apply the suitable entropy estimates from Lemma~\ref{l:scaled-ee} with $r_\eps = 2$ and $\delta_\eps = 1/2\wedge\bar\delta$ and $\kappa_\eps =1/2$ in order to get 
\[
\begin{aligned}
\esssup_{t\in (-4,0]}   \int_{B_2(0)} h_+^{1/2} (f_\eps (t,v)) \dv + 4\int_{Q_2} \left|\nabla_v \sqrt{f_\eps}(t,v)\right|^2 \un_{f_\eps \ge 1/2} \dt \dv 
\\
\le  D_0 \int_{Q_{3} } f_\eps \bigg(\ln_+ (2f_\eps)+\ln_+^2 (2f_\eps) \bigg)(t,v) \dt \dv 
\end{aligned}
\]
for some constant $D_0$ depending only on $|v_0|$, $m_0$, $M_0$, $E_0$, $H_0$ and $\rho_0$. The constant $D_0$ involves 
\[
\esssup_{Q_3}Z[f_\eps]\le 2(1-\gammas)M_0+\gammas\esssup_{Q_3} \left(\eps^{-\gamma-3} f_\eps\ast_v\frac1{|\cdot|}\right)\,,
\]
as explained above. Besides
\[
\eps^{-\gamma-3} f_\eps(9\hat t,\cdot)\ast\frac1{|\cdot|}(3\hat v)= (3\eps)^{-\gamma-3} f_{3\eps}(\hat t,\cdot)\ast\frac1{|\cdot|}(\hat v)\,,
\]
so that, according to Lemma \ref{l:lri},
\[
\esssup_{Q_3} \left(\eps^{-\gamma-3} f_\eps\ast_v\frac1{|\cdot|}\right)=\esssup_{Q_1} \left((3\eps)^{-\gamma-3}f_{3\eps}\ast_v\frac1{|\cdot|}\right)\le \frac{C_*}{\rho_0}\sup_{t\ge 0}\int_{\R^3}f(1+\ln_+f)(t,w)\dw+C_*\rho_0^2\,.
\]
We then argue as in the proof of De Giorgi's lemma: on the left-hand side, we observe that
\[ 
h_+^{1/2} (f_\eps) \ge h_+^{1/2} (f_\eps) -h_+^{1/2} (1) \ge (h_+^{1/2})' (1) (f_\eps - 1) = (\ln 2) (f_\eps-1)\,.
\]
In particular, 
\begin{align}
\nonumber  
(\ln 2)  \esssup_{t\in (-4,0]}  & \int_{B_2}   (f_\eps (t,v)-1)_+ \dv + 4 \int_{Q_2} \left|\nabla_v \sqrt{f_\eps}(t,v)\right|^2 \un_{f_\eps \ge 1} \dt \dv 
\\
\nonumber 
& \le  D_0 \int_{Q_{3} } f_\eps \bigg(\ln (2f_\eps)+\ln^2 (2f_\eps) \bigg)(t,v) \un_{f_\eps (t,v) \ge 1/2} \dt \dv\,,
\\
\intertext{using that $\lnp (r) + \lnp^2 (r) \le C_\ast (1+r)$ for some absolute constant $C_\ast>0$ to get}
\label{e:ee-before-dg} 
& \le  D_0 C_\ast \int_{Q_{3} } (f_\eps+2f_\eps^2)(t,v)  \un_{f_\eps (t,v) \ge 1/2} \dt \dv \le 4 D_0 C_\ast \int_{Q_{3} } f_\eps^2(t,v)  \dt \dv\,.
\end{align}
Since $f$ is axisymmetric with respect to $e_3$, there exists $\mathcal F$ defined a.e. on $[0,T]\times[0,+\infty)\times\R$ and measurable such that $f(t,(v',v_3)) = \mathcal F (t,|v'|,v_3)$. In particular, since $dtdv=\eps^5d\bar td\bar v$ 
\begin{align*}
\int_{Q_3} f_\eps^2(t,v)  \dt \dv &= \eps^{2\gamma+10} \int_{Q_3} f(t_0 + \eps^2t, (v'_0+\eps v',v^0_3+\eps v_3))^2\dt \dv 
\\
&=\eps^{2\gamma+5} \int_{Q_{3\eps}}\mathcal F(t_0+\bar t,|v'_0+\bar v'|,v_3^0+\bar v_3)^2\dd\bar t\dd\bar v
\\
& \le \frac1\eps\int_{Q_{3\eps}}\mathcal F(t_0+\bar t,|v'_0+\bar v'|,v_3^0+\bar v_3)^2\dd\bar t\dd\bar v\,.
\end{align*}
\begin{figure}
\begin{center}
\includegraphics[width=12cm]{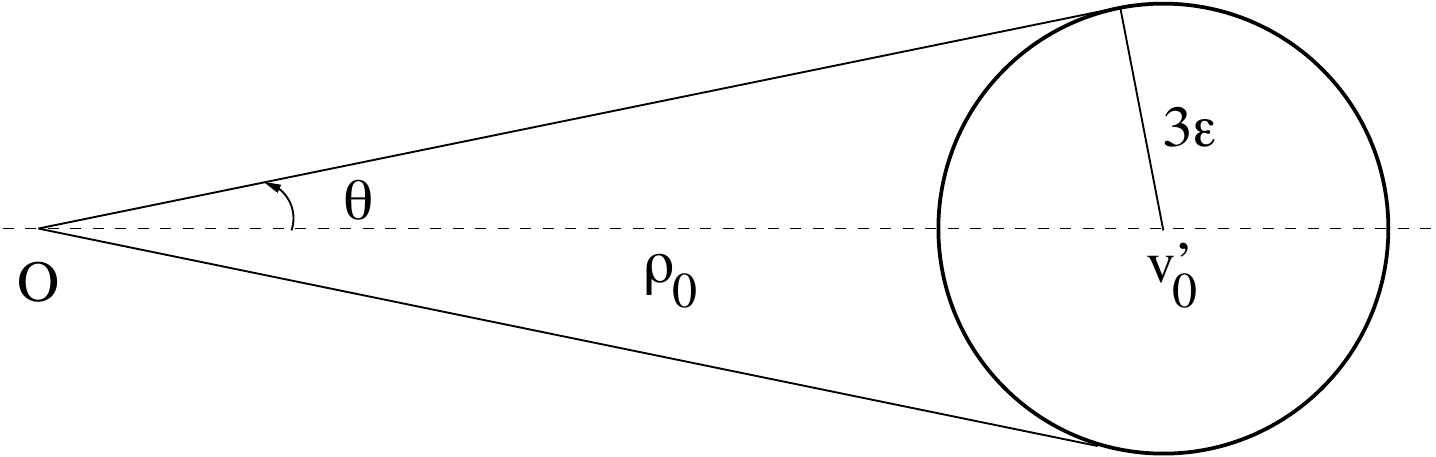}
\end{center}
\caption{The intersection of the orthogonal projection of the ball $B_{3\eps}(v_0)$ on the affine plane passing through the origin with direction $(\R e_3)^\perp$. The angle $\theta$ is $\arcsin\left(\frac{3\eps}{\rho_0}\right)$.}
\end{figure}

Now $(t_0,v_0)+Q_{3\eps}\subset\mathcal O_{3\eps}\subset(0,T)\times\R^3$, where the set $\mathcal O_{3\eps}$ is defined in the cylindrical coordinates $(\rho:=|v'_0+\bar v'|,\theta,w_3)$ for $(v'_0+\bar v',v^0_3+\bar v_3)\in\R^3$ by 
the inequalities
\[
t_0-9\eps^2<s\le t_0\,,\quad \rho_0-3\eps<\rho=|v'_0+\bar v'|<\rho_0+3\eps\,,\quad v_3^0-3\eps<w_3<v_3^0+3\eps\,,\quad|\theta|\le\arcsin\left(\frac{3\eps}{\rho_0}\right)\,;
\]
see Figure 1. Hence
\begin{equation}\label{AxiIneq}
\begin{aligned}
\int_{Q_3} f_\eps^2(t,v)  \dt \dv &\le\frac1\eps\int_{\mathcal O_{3\eps}}\mathcal F(t_0+\bar t,|v'_0+\bar v'|,v_3^0+\bar v_3)^2\dd\bar t\dd\bar v
\\
&=\frac1\eps\int_{-\arcsin(\frac{3\eps}{\rho_0})}^{+\arcsin(\frac{3\eps}{\rho_0})}d\theta\int_{t_0-9\eps^2}^{t_0}\int_{\rho_0-3\eps}^{\rho_0+3\eps}\int_{v_3^0-3\eps}^{v_3^0+3\eps}\mathcal F(s,\rho,w_3)^2\dd s\rho\dd\rho\dd w_3
\\
&=\frac2\eps\arcsin(\frac{3\eps}{\rho_0})\int_{t_0-9\eps^2}^{t_0}\int_{\rho_0-3\eps}^{\rho_0+3\eps}\int_{v_3^0-3\eps}^{v_3^0+3\eps}\mathcal F(s,\rho,w_3)^2\dd s\rho\dd\rho\dd w_3
\\
&\le\frac{3\pi}{\rho_0}\int_{t_0-9\eps^2}^{t_0}\int_{\rho_0-3\eps}^{\rho_0+3\eps}\int_{v_3^0-3\eps}^{v_3^0+3\eps}\mathcal F(s,\rho,w_3)^2\dd s\rho\dd\rho\dd w_3
\\
&\le\frac{3\pi}{\rho_0}(\rho_0+3)\int_{t_0-9\eps^2}^{t_0}\int_{\rho_0-3\eps}^{\rho_0+3\eps}\int_{v_3^0-3\eps}^{v_3^0+3\eps}\mathcal F(s,\rho,w_3)^2\dd s\dd\rho\dd w_3\,,
\end{aligned}
\end{equation}
using that $\arcsin y\le\tfrac\pi{2} y$ for all $y\in[0,1]$ to prove the penultimate inequality above.

\smallskip
We know from Lemma~\ref{l:better-int} that $\mathcal F \in L^2 (\mathcal{Q}_r (t_0,(\rho_0,v_3^0)))$ for all $r>0$ small enough. Therefore
\[ 
\int_{t_0-9\eps^2 }^{t_0} \int_{\rho_0-3\eps}^{\rho_0+3\eps} \int_{v_3^0-3\eps}^{v_3^0 + 3\eps} \mathcal F(s, \rho, \bar v_3)^2 \ds \dd \rho \dd \bar v_3 \to 0
\]
as $\eps \to 0^+$ by dominated convergence. In particular, there exists $\eps >0$ such that
\[  
\int_{Q_3 } f_\eps^2(t,v)  \dt \dv \le \frac{( \ln 2) \etadg Z_0^{-3/2}}{4 D_0 C_\ast} 
\]
where $C_\ast$ is the absolute constant from \eqref{e:ee-before-dg}. We conclude from \eqref{e:ee-before-dg} and \eqref{AxiIneq} that \eqref{e:cond-dg1} holds true. Now the De Giorgi lemma (Lemma \ref{l:dg1}) implies that $f$ must be bounded 
in some neighborhood of $(t_0,v_0)$. Hence $(t_0,v_0)\notin\mathbf S[f]$, as announced.
\end{proof}


\appendix

\section{Appendix: known facts and technical results}
\label{s:appendix}


\subsection{Bounded solutions are smooth}
\label{a:smooth}

Throughout this section, we assume that $k(r)=r^{\gamma+3}$ with $\gamma\in[-3,-2)$.

\begin{prop}[Smoothness] \label{p:smooth} 
Let $f$ be a suitable weak solution of the Landau equation \eqref{e:landau} in $[0,T) \times \R^3$ in the sense of Definition~\ref{d:suitable}.  Assume that $f$ is essentially bounded in $Q_{2r}(t_0,v_0)$. Then $f$ is a.e. equal to a $C^\infty$ function
in $Q_r (t_0,v_0)$.
\end{prop}

\begin{proof}
Desvillettes has proved in \cite{d15} that $H$-solutions are in fact weak solutions of the Landau equation. In particular, they satisfy
 \[ 
 \partial_t f = \nabla_v \cdot (A \nabla_v f) + \nabla_v \cdot g \,,
 \]
 where $g := f B$ and 
 \[
 A (t,v) := \int_{\R^3} f(t,v-z) |z|^{-1}k(|z|) \Pi (w) \dd z \quad \text{ and } \quad B (t,v) :=  2 \int_{\R^3} f(t,v-z) |z|^{-3}k(|z|) z \dd z.
 \]
Notice that for all $\beta_1 \in (1,2)$, the terms $A$ and $g$ are essentially bounded in $Q_{\beta_1 r}(t_0,v_0)$ as soon as $f$ is essentially bounded in $Q_{2r}(t_0,v_0)$. Moreover  \cite[Lemma~3.1]{zbMATH06668299} (see also 
Lemma~\ref{l:diff-mat}) implies that the  diffusion matrix is strictly elliptic,
\[  
A (t,v) \ge c (1+|v|)^{-3} I \,,
\]
where $c$ depends on $\int_{\R^3} \fin (v) \dv$, $\int_{\R^3} \fin (v) |v|^2 \dv$ and $\int_{\R^3} \fin \ln \fin (v) \dv$.  In particular, $A$ is uniformly elliptic in $Q_{\beta_1 r}(t_0,v_0)$ if $f$ is essentially bounded in $Q_{2r}(t_0,v_0)$:
\[
\text{there exists }\lambda, \Lambda >0\text{ such that, for all } (t,v) \in Q_{\beta_1 r}(t_0,v_0), \qquad \lambda I \le A(t,v) \le \Lambda I \,,
\]
for some ellipticity constants $\lambda, \Lambda$ depending only on $\gamma$, $c$, $|v_0|$, $r$, $\beta_1$ and $\|f\|_{L^\infty (Q_{2r}(t_0,v_0))}$ and hydrodynamic bounds. This implies that $f$ lies in a parabolic De Giorgi class (see \cite{L1} 
and \cite[p.~143-150]{MR1465184}), and thus is locally H\"older continuous in $Q_{\beta_2 r}(t_0,v_0)$ for all $\beta_2 \in (1,2)$. In particular,
\begin{equation}\label{e:local-holder}
[f]_{C^\alpha (Q_{\beta_2 r}(t_0,v_0))} \le C ( \|f\|_{L^\infty(Q_{\beta_1 r}(t_0,v_0))} + \| g\|_{L^\infty (Q_{\beta_1 r}(t_0,v_0))})
\end{equation}
for some $C>0$ depending only on $\lambda,\Lambda$ and $\beta_1>\beta_2,r$ and $\alpha \in (0,1)$ depending only on $\lambda,\Lambda$.
  
We next show that $A$ and $B$ are both $\alpha$-H\"older continuous in $Q_{\beta_3 r} (t_0,v_0)$ for all $\beta_3 \in (1,2)$. Indeed, pick $1 < \beta_3 <\beta_2 <2$ and let $\phi:\R^3 \to [0,1]$ be a $C^\infty$ function such that
\[
\un_{B_{(\beta_2-\beta_3)r/2}}\le\phi\le\un_{B_{(\beta_2-\beta_3)r}}\,.
\]
Write $A = A_\inn+A_\out$ with
\begin{align*}
A_\inn (t,v) & := \int_{\R^3} f(t,v-z) \phi (z) |z|^{-1}k(|z|) \Pi (z) \dd z\,,
\\
A_\out (t,v) &: = \int_{\R^3} f(t,v-z) (1-\phi) (z) |z|^{-1}k(|z|) \Pi (z) \dd z\,. 
\end{align*}
We can use the local H\"older estimate \eqref{e:local-holder} to prove that $A_\inn$ is $\alpha$-H\"older continuous in $Q_{\beta_3 r} (t_0,v_0)$ since $\text{supp}(\phi)\subset B_{(\beta_2-\beta_3)r}$ while $k(|\cdot|)/|\cdot|\in L^1_{loc}(\R^3)$. 
On the other hand, $(1-\phi)k(|\cdot|)/|\cdot|$ is bounded and smooth on $\R^3$, so that $A_\out$ is also $\alpha$-H\"older continuous in $Q_{\beta_3 r} (t_0,v_0)$ since $\int_{\R^3} f(t,v) \dv = \int_{\R^3} \fin (v) \dv<+\infty$. Summarizing, we see that
\[
[A]_{C^\alpha (Q_{\beta_3 r} (t_0,v_0))} \le C (\|f\|_{C^\alpha (Q_{\beta_2 r} (t_0,v_0))}+1)
\]
for some $C$ depending on the initial mass $\int_{\R^3} \fin (v) \dv$ and $\beta_2>\beta_3,\gamma, r$. The same reasoning yields
\[
[B]_{C^\alpha (Q_{\beta_3 r} (t_0,v_0))} \le C (\|f\|_{C^\alpha (Q_{\beta_2 r} (t_0,v_0))}+1)
\]
for some $C$ depending on the initial mass $\int_{\R^3} \fin (v) \dv$ and $\beta_2>\beta_3,\gamma,r$. This implies in particular
\[
[g]_{C^\alpha (Q_{\beta_3 r} (t_0,v_0))} \le C (\|f\|_{C^\alpha (Q_{\beta_2 r} (t_0,v_0))}+1)\,.
\]
Then, we can use the Schauder estimates for parabolic equations in divergence form (see \cite[p.~56]{MR1465184}), and conclude that the first derivatives of $f$ (and those of $A$ and $g$ as well) are $\alpha$-H\"older continuous in 
$Q_{\beta_4 r} (t_0,v_0)$ for all $\beta_4 \in (1,2)$, with a quantitative estimate. We then differentiate both sides of the Landau equation and conclude by a classical bootstrap argument.
\end{proof}

\subsection{Passing to the limit in local truncated entropy inequalities}\label{a:trunc-entropineq}

\begin{proof}[Proof of Proposition~\ref{p:ee}]
Let $f^n$ be the approximate solution from Proposition~\ref{p:approximate}. We first write the evolution of the localized truncated entropy functional for $f^n$ and we then explain how to pass to the limit to get \eqref{e:lee} and \eqref{e:diff-drift}. 
  
\paragraph{Getting started.}
Let $A_n (t,v)$ denotes the (approximate) diffusion matrix $\int_{\R^3} a_n(v-w) f^n(t,w) \dw$, and recast equation~\eqref{e:landau-reg} as
\[ 
\partial_t f^n  = \nabla_v \cdot ( A_n \nabla_v f^n - (\nabla_v \cdot A_n) f^n ) +\frac1{n} \Delta_v f^n\,.
\]
Then
\[
\left[ \int_{\R^3} h_+^\kappa (f^n)  \varphi  \dv \right]^{t_2}_{t_1}  =\int_{t_1}^{t_2} \left(\mathcal{T}_1^n +\mathcal{T}_2^n+\mathcal{T}_3^n+ \frac1n \mathcal{E}^n\right)\dt\,,
\]
\[ 
\text{ with } \quad 
\begin{cases}
\displaystyle\mathcal{T}_1^n & := - \displaystyle\int_{\R^3} (A_n \nabla_v f^n - (\nabla_v \cdot A_n)f^n) \cdot \un_{f^n\ge\kappa} \frac{\nabla_v f^n}{f^n}  \varphi \dv\,, 
\\
\mathcal{T}_2^n & = - \displaystyle\int_{\R^3} (A_n \nabla_v f^n - (\nabla_v \cdot A_n)f^n) \cdot \lnp \left(\frac{f^n}{\kappa}\right) \nabla_v \varphi \dv\,, 
\\
\mathcal{T}_3^n & = \displaystyle\int_{\R^3} h_+^\kappa (f^n) \partial_t \varphi \dv\,, 
\\
\mathcal{E}^n   & = -\displaystyle\int_{\R^3} \varphi \frac{|\nabla_v f^n|^2}{f^n} \un_{f^n\ge\kappa} \dv  -  \int_{\R^3} \lnp \left(\frac{f^n}\kappa\right) \nabla_v f^n \cdot \nabla_v \varphi  \dv\,,
\end{cases}
\]
for all $r>0$.

\paragraph{The error term.} We begin with the $\mathcal{E}^n$ term. Recalling that $\varphi = \Psi^2$, the inequality $2|a||b|\le|a|^2+|b|^2$ implies that
\begin{align}
\nonumber  
{}&\mathcal{E}^n =  -\int_{\R^3}  \Psi^2 (v) \frac{|\nabla_v f^n|^2}{f^n} \un_{E_\kappa^n} \dv - 2 \int _{\R^3}\sqrt{f^n} \lnp (f^n/\kappa)  \nabla_v \Psi  \cdot \Psi \frac{\nabla_v f^n (v)}{\sqrt{f^n}} \un_{E_\kappa^n} \dv\,,
\\
\label{e:e1} 
&\text{so that}\qquad|\mathcal{E}^n |\le 2\int_{\R^3}  \Psi^2 (v) \frac{|\nabla_v f^n|^2}{f^n} \un_{E_\kappa^n} \dv+\int_{\R^3}f^n \left(\lnp (f^n/\kappa)\right)^2  | \nabla_v \Psi|^2    \dv\,.
\end{align}
In view of \eqref{e:gradfn} and \eqref{e:fn-equi}, we conclude that $\tfrac1n\mathcal{E}^n\to 0$ in $L^1_{loc}((0,+\infty))$ as $n$ tends to $\infty$.

\paragraph{Entropy terms.} Let $\Psi$ and $\phi$ be supported in $[0,T] \times\bar{B}_R(0)$. Since $f^{n_i} \to f$ a.e. in $(0,T) \times \R^3$, see \eqref{e:ae-conv}, it is enough to check that $ h_+^\kappa (f^n)$ is uniformly integrable in 
$(0,T) \times B_{R}$ in order to conclude that,
\begin{eqnarray*}
\int_{t_1}^{t_2} \mathcal{T}_3^{n_i} \dt = \iint_{[t_1,t_2] \times \R^3} h_+^\kappa (f^{n_i}) \partial_t \varphi  \dt \dv &\to& \iint_{[t_1,t_2] \times \R^3} h_+^\kappa (f) \partial_t \varphi  \dt \dv = \int_{t_1}^{t_2} \mathcal{T}_3 (f) \dt\,,
\\
\text{for a.e. }t \in (0,T), \quad  \int_{\R^3} h_+^\kappa (f^{n_i}(t,\cdot))  \varphi  \dv &\to& \int_{\R^3} h_+^\kappa (f(t,\cdot))  \varphi  \dv \,,
\end{eqnarray*}
as $n_i \to \infty$. In order to get convergence a.e. in the time variable, it is necessary to consider a subsequence of $\{ f^{n_i} \}$. In the sequel, we may have to further extract subsequences from $\{ f^{n_i} \}$, but shall keep writing 
$f^{n_i}$ for the sake of simplicity.

Uniform integrability is an easy consequence of the fact that $h_+ (r) \le C (1 + (r-1)_+^{4/3})$ for some  absolute constant $C$. Indeed, it implies that for $\kappa \ge 1$, 
\[  
h_+^\kappa (f^n ) \le C (\kappa + (f^n-\kappa)_+^{4/3})\,,
\]
from which we deduce that for some modified absolute constant $\tilde C$, we have
\[ 
\bigg( h_+^\kappa (f^n ) \bigg)^{{5/4}} \le \tilde C (\kappa^{5/4} + (f^n-\kappa)_+^{5/3})\,. 
\]
Estimate \eqref{e:fn-equi} yields the uniform integrability of $h_+^\kappa (f^n)$ in the bounded set $(0,T) \times B_{R}$.    

\bigskip
Passing to the limit in $\mathcal{T}_1^n$ and $\mathcal{T}_2^n$ is more difficult.

\begin{lemma}[Passing to the limit in $\mathcal{T}_1^n$]\label{l:t1-lim}
Let $\aout$ and $\ain$ be defined as in \eqref{e:aout}. Then
\[
\limsup_{n_i \to \infty} \int_{t_1}^{t_2} \mathcal{T}_1^{n_i}(t) \dt \le \int_{t_1}^{t_2} \mathcal{T}_1(f)(t) \dt\,,
\]
 where
\begin{align*}
\mathcal{T}_1(f)(t) =& -\int_{\R^3} ( A(t,v) \nabla_v f(t,v) - (\nabla_v \cdot A(t,v)) f(t,v)) \frac{\nabla_v \fkp}{f}(t,v) \varphi (t,v)\dv
\\
:=& -\tfrac12 \iint_{\R^3 \times \R^3} \left|F^{\kappa,\inn}_+ (t,v,w)\right|^2 \varphi (t,v) \dv \dw - \int_{\R^3} A^\out (t,v) \colon \frac{(\nabla_v \fkp)^{\otimes 2}}{f} (t,v) \varphi (t,v)  \dv 
\\
& +8\pi\ka(0)\kappa\! \int_{\R^3} \fkp (t,v) \varphi (t,v) \dv\!+\!2\iint_{\R^3\times\R^3}\!\frac{k'(|v\!-\!w|)}{|v\!-\!w|^2}f^\kappa_+(t,v)(f(t,w)\wedge\kappa)\phi(t,v)\dv\dw \!+\! \mathcal{E}_1 (f)(t)
\end{align*}
with $A^\out (t,v): = \int_{\R^3} \aout(v-w) f(t,w) \dw$, and, recalling \eqref{e:fkappain} for the definition of $F^{\kappa,\inn}_+$, where
\begin{align*}
\mathcal{E}_1 (f)(t) := &\,\tfrac12 \iint_{\R^3 \times \R^3} F^{\kappa,\inn}_+ (t,v,w) \cdot \bigg( (\varphi (t,v) - \varphi (t,w)) \sqrt{\anin (v-w)} \bigg)\sqrt{f} (t,v) \frac{\nabla_w \fkp}{\sqrt{f}} (t,w)  \dv \dw 
\\
& - \int_{\R^3} (\nabla_v \cdot A^\out) (t,v)  \fkp (t,v) \cdot \nabla_v \varphi (t,v)  \dv   + \int_{\R^3} (-\nabla_v^2 : A^\out) (t,v) \fkp (t,v) \varphi (t,v)  \dv 
\\
& + \int_{\R^3 \times \R^3} (\nabla_z^2 : \aout )(v-w) (f \wedge \kappa) (t,w)  \fkp (t,v) \varphi (t,v)  \dv \dw 
\\
& - \iint_{\R^3 \times \R^3} (\nabla_z \cdot \ain) (v-w) \cdot \nabla_v \varphi (t,v) (f \wedge \kappa)(t,w) \fkp (t,v) \dv \dw\,.
\end{align*}
\end{lemma}

\medskip

\begin{lemma}[Passing to the limit in $\mathcal{T}_2^n$]\label{l:t2-lim}
One has 
\[
\lim_{n_i \to \infty} \int_{t_1}^{t_2} \mathcal{T}_2^{n_i}(t)\dt = \int_{t_1}^{t_2}\mathcal{T}_2(f)(t) \dt\,,
\] 
where
\begin{align*}
\mathcal{T}_2 (f)(t) = & - \int_{\R^3} ( A(t,v) \nabla_v f(t,v) - (\nabla_v \cdot A)(t,v) f(t,v)) \lnp (f/\kappa)(t,v) \cdot \nabla_v \varphi(t,v) \dv
\\
:= & - \iint_{\R^3 \times \R^3} F(t,v,w) \cdot G^\inn (t,v,w) \nabla_v \varphi (t,v) \dv \dw  +  \int_{\R^3}  h^\kappa_{+} (f) (t,v) A^\out (t,v)    \colon D^2_v \varphi (t,v)  \dv  
\\
&+ \int_{\R^3} (\nabla_v \cdot A^\out) (t,v)  \cdot \left ( f(t,v)  \lnp \left(\frac{f}{\kappa}\right) (t,v) + h^\kappa_+(f)(t,v)\right) \nabla_v \varphi (t,v)  \dv\,,
\end{align*}
with $F$ defined in \eqref{e:Fn} and 
\[
G^\inn (t,v,w) := (1-X(\delta^{-1}|v-w|))  \sqrt{a(v-w)} \sqrt{ f (t,v) f(t,w) } \lnp \left(\frac{f}{\kappa}\right) (t,v)\,.
\]
\end{lemma}

\medskip
Both these lemmas and the reasoning above lead to \eqref{e:lee} with \eqref{e:diff-drift}. Indeed combining the second term in the definition of $\mathcal E_1(f)$ with the last term in the definition of $\mathcal T_2(f)$ produces the second term
in the definition of $\mathcal E_\out(f)$ in \eqref{e:diff-drift}. Similarly, combining the third and fourth terms in the definition of $\mathcal E_1(f)$ produces the last term in the definition of $\mathcal E_\inn(f)$ in \eqref{e:diff-drift}.
\end{proof}

\medskip
We are therefore left with the task of proving Lemmas~\ref{l:t1-lim} and \ref{l:t2-lim}.

As a preliminary step, we study the sequence of diffusion matrices $A_n (t,v) = A^\inn(t,v)+A^\out(t,v)$, separating long-range from short-range interactions For all $\delta >0$, define
\begin{equation}\label{DefAinAout}
A^\inn_n (t,v) := \int_{\R^3} a^\inn_n (v-w) f^n (t,w) \dw\quad\text{ and }\quad A^\out_n (t,v) := \int_{\R^3} a^\out_n (v-w) f^n (t,w) \dw\,,
\end{equation}
recalling the definitions of $a^\inn_n$ and $a^\out_n$ in \eqref{e:anin}. Observe that
\[
n>2/\delta\text{ and }z\not=0\implies X(n|z|)X(\tfrac{|z|}\delta)=X(\tfrac{|z|}\delta)\implies a^\out_n(z)=\aout(z)\,.
\]

\begin{lemma}[Diffusion matrix -- long-range interactions]\label{l:anout}
Let $\delta > 2/n$. For a.e. $(t,v) \in (0,T) \times \R^3$, 
\[
|A^\out_n (t,v)| \le \frac{2\ka (\delta)M_0}{\delta}\,, \qquad |\partial_{v_l} A^\out_n (t,v)| \le \frac{8 \ka (\delta) M_0}{\delta^2}\,,  \qquad | \partial_{v_l,v_m} A^\out_n (t,v)| \le \frac{40 \ka (\delta) M_0}{ \delta^3}\,.
\]
Moreover, $\{A^\out_{n_i}\}$, $\{\partial_{v_l} A^\out_{n_i}\}$, $\{\partial_{v_l,v_m} A^\out_{n_i}\}$  converge respectively  to $A^\out$, $\partial_{v_l} A^\out$ and $\partial_{v_l,v_m} A^\out$ a.e. in $(0,T) \times B_R$ as $n_i \to \infty$. 
\end{lemma}

\begin{proof}
The estimate follows from \eqref{e:div-aout} and the fact that $X \in [0,1]$ supported in $[\tfrac12,+\infty)$. We also use that the function $\ka$ is nondecreasing and in particular $\ka (\delta/2) \le \ka (\delta)$. Convergence follows from a.e. convergence 
of $f^{n_i}$ to $f$, see \eqref{e:ae-conv}, while uniform integrability is obtained from \eqref{e:fn-equi}. 
\end{proof}

\begin{lemma}[Short-range interactions]\label{l:d2anin}
Let $\delta > 2/n$. The first and second derivatives of the matrix field $\anin$ defined in \eqref{e:anin} satisfy the following formulas:
\begin{align}
\label{e:d1anin}  
\nabla_z \cdot \anin (z) &= \ka (|z|)X(n|z|)(1-X(\delta^{-1}|z|)) \nabla_z \left(\frac2{ |z|} \right) = - 2\ka (|z|)X(n|z|)(1-X(\delta^{-1}|z|)) \frac1{|z|^2} \frac{z}{|z|}\,,
\\
\label{e:d2anin} 
\nabla_z^2 : \anin (z) & = \bigg(-\ka' (|z|) X(n|z|) (1-X(\delta^{-1}|z|)) - \ka (|z|) nX'(n|z|)  +  \ka (|z|) \delta^{-1}X'(\delta^{-1}|z|) \bigg) \frac2{|z|^2}\,.
\end{align}
\end{lemma}

\begin{proof}
If $\ka (|z|)$ is replaced with $\ka (|z|) X (n|z|) (1-X(|z|/\delta))$, we can apply  \eqref{e:div-a} and get \eqref{e:d1anin}. As far as second derivatives are concerned, we apply \eqref{e:hess-a} after updating $\ka$. Since $n>2/\delta$, for all $r>0$, 
if $X'(n r) \neq 0$ then $X(r/\delta)=0$, and if $X'(r/\delta) \neq 0$ then $X (n r)=1$. This implies that, for all $r>0$,
\[ 
(\ka (r) X (n r) (1-X(r/\delta)))' = X (n r) (1-X(r/\delta)) \ka' (r) + \ka (r) n X'(nr) - \ka (r) \delta^{-1} X'(\delta^{-1}r)\,.
\]
This yields the second announced formula.
\end{proof}

\medskip
We now turn to the proof of Lemma~\ref{l:t1-lim}.

\begin{proof}[Proof of Lemma~\ref{l:t1-lim}]
We first rewrite $\mathcal{T}_1^n$ as follows,
\[
\mathcal{T}_1^n \!=\! -\! \iint_{\R^3 \times \R^3}\!a_n (v\!-\!w) f^n (t,v) f^n (t,w) \left( \frac{\nabla_v f^n}{f^n} (t,v) \!-\! \frac{\nabla_w f^n}{f^n}(t,w)\right)\cdot \un_{f^n(t,v)\ge\kappa} \frac{\nabla_v f^n}{f^n} (t,v) \varphi (t,v)  \dv \dw\,.
\]
In view of this expression and recalling that $a_n = \anin + \aout$ for $\delta > 2/n$, it is natural to introduce the term $L^n_+ := \lnp \left(\frac{f^n}{\kappa} \right)$ with $\lnp (r) =  \max(\ln r, 0)$. With such a definition at hand, we can write,
\[ 
\mathcal{T}_1^n = \mathcal{T}_1^{+,n} +\mathcal{T}_1^{\wedge,n} + \mathcal{T}_1^{\out,n}\,,
\]
with
\begin{eqnarray*}
\mathcal{T}_1^{+,n}& := & - \iint_{\R^3 \times \R^3} \anin (v-w) f^n (t,v) f^n (t,w) \left( \nabla_v L^n_+ (t,v)- \nabla_w L^n_+ (t,w)\right) \cdot \nabla_v L^n_+ (t,v) \varphi (t,v)  \dv \dw \,,
\\
\mathcal{T}_1^{\wedge,n}& := &  \iint_{\R^3 \times \R^3} \anin (v-w)  \bigg( (1-\un_{f^n(t,w)\ge\kappa})\nabla_w f^n \bigg) (t,w) \cdot \bigg( \un_{f^n(t,v)\ge\kappa} \nabla_v f^n \bigg) (t,v) \varphi (t,v) \dv \dw\,,
\\
\mathcal{T}_1^{\out,n} & := & - \iint_{\R^3 \times \R^3} \aout (v-w) f^n (t,v) f^n (t,w) \left( \frac{\nabla_v f^n}{f^n} (t,v) - \frac{\nabla_w f^n}{f^n}(t,w)\right) \cdot \nabla_vL^n_+(t,v)\varphi (t,v)  \dv \dw\,.
\end{eqnarray*}

\paragraph{Short-range interactions.}
We use once again  that $\varphi = \Psi^2$. It allows us to write $\mathcal{T}_1^+$ as follows,
\begin{align}
\nonumber 
\mathcal{T}_1^{+,n} &
\\
\nonumber  
= & -\!\iint_{\R^3 \times \R^3}\!\anin (v\!-\!w) f^n (t,v) f^n (t,w) \left( \nabla_v L^n_+ (t,v) \Psi (t,v)\!-\! \nabla_w L^n_+ (t,w) \Psi (t,w) \right) \!\cdot\! \nabla_v L^n_+ (t,v) \Psi (t,v)  \dv \dw 
\\
\nonumber  
&  + \iint_{\R^3 \times \R^3} \anin (v-w) f^n (t,v) f^n (t,w)  \nabla_w L^n_+ (t,w) (\Psi (t,v) - \Psi (t,w)) \cdot \nabla_v L^n_+ (t,v) \Psi (t,v)  \dv \dw 
\\
\intertext{symmetrizing in $(v,w)$  and using that $f^n \nabla_v L^n_+ = \nabla f^{n,\kappa}_+$ where $f^{n,\kappa}_+ =  (f^n - \kappa)_+$,}
\label{e:firstterm}    
= & - \tfrac12 \iint_{\R^3 \times \R^3} \anin (v-w) f^n (t,v) f^n (t,w) : \left( \nabla_v L^n_+ (t,v) \Psi (t,v) - \nabla_w L^n_+ (t,w) \Psi (t,w) \right)^{\otimes 2} \dv \dw 
\\
\label{e:secondterm}  
&  + \tfrac12 \iint_{\R^3 \times \R^3} \anin (v-w)  (\Psi (t,v) - \Psi (t,w))^2 \nabla_w f^{n,\kappa}_+ (t,w) \cdot \nabla_v f^{n,\kappa} _+ (t,v)   \dv \dw\,.
\end{align}
We study the limit of each term separately. \bigskip

\paragraph{Short-range interactions (main term).}
As far as the first term is concerned, see \eqref{e:firstterm}, we write it as one half of the opposite of the $L^2$-norm of
\[
\sqrt{\anin (v-w)}\sqrt{ f^n (t,v) f^n (t,w)}  \left( \nabla_v L^n_+ (t,v) \Psi (t,v)- \nabla_w L^n_+ (t,w) \Psi (t,w) \right)\,.
\]
We shall prove that this sequence weakly converges in $L^2((0,T) \times B_R \times B_{R+\delta})$ to 
\[
\sqrt{\ain (v-w) }\sqrt{f (t,v) f (t,w)}  \left( \nabla_v L_+ (t,v) \Psi (t,v)- \nabla_w L_+ (t,w) \Psi (t,w) \right)
\]
with $L_+  = \lnp (f/\kappa)$. 

In order to prove this weak convergence, we simply observe that
\begin{multline*}
\sqrt{\anin (v-w) }\sqrt{f^n (t,v) f^n (t,w)}  \left( \nabla_v L^n_+ (t,v) \Psi (t,v)- \nabla_w L^n_+ (t,w) \Psi (t,w) \right) 
\\
= F^{\kappa,\inn,n}_+ (t,v,w) \Psi (t,v) + \sqrt{\anin (v-w) }\sqrt{f^n (t,v) f^n (t,w)} \nabla_w L^n_+ (t,w) (\Psi (t,v) - \Psi (t,w))
\end{multline*}
where
\[ 
F^{\kappa,\inn,n}_+ (t,v,w) :=  \sqrt{\anin (v-w) }\sqrt{f^n (t,v) f^n (t,w)}  \left( \nabla_v L^n_+ (t,v) - \nabla_w L^n_+ (t,w)  \right)\,. 
\]
We know from \eqref{e:dissipk-conv} that $F^{\kappa,\inn,n_i}_+ (t,v,w)$ converges towards $F^{\kappa,\inn}_+$ weakly in $L^2_{\text{loc}}((0,+\infty)\times\R^3\times\R^3)$. 

Moreover, recalling that $f^n \nabla_v L^n_+ = (f^n-\kappa)_+=\fnkp$ and $\anin(z) = X (n |z|)\ain (z)$, 
\begin{align*}
\sqrt{\anin (v-w)}\sqrt{ f^n (t,v) f^n (t,w)} & \nabla_w L^n_+ (t,w)  (\Psi (t,v) - \Psi (t,w)) 
\\
& =  \sqrt{X (n |v-w|)} \bigg( \sqrt{\ain (v-w)} (\Psi (t,v) - \Psi (t,w)) \bigg) \sqrt{f^n} (t,v) \frac{\nabla_w \fnkp}{\sqrt{f^n}} (t,w)\,.
\end{align*}

We claim  that
\begin{equation}\label{e:claim-L2weak}
\text{Claim :} \quad \sqrt{f^{n_i}} (t,v) \frac{\nabla_w \fnikp}{\sqrt{f^{n_i}}} (t,w)\text{ converges to }  \sqrt{f} (t,v) \frac{\nabla_w \fkp}{\sqrt{f}} (t,w) \text{ weakly in } L^2_{\text{loc}}((0,+\infty)\times\R^3\times\R^3) \,.
\end{equation}
 
In order to justify this claim, we first study the limit in the sense of distributions. Letting $g^n$ denote $g^n:=\frac{\nabla_w \fnkp}{\sqrt{f^n}}$, we pick $\phi \in C^\infty_c (Q_{T,R,\delta})$ with $Q_{T,R,\delta} = (0,T) \times B_R \times B_{R+\delta}$ 
and prove that 
\begin{equation}\label{e:distribution-limit}
\int_{Q_{T,R,\delta}} \sqrt{f^{n_i}} (t,v) g^{n_i} (t,w) \phi (t,v,w) \dt \dv \dw \rightarrow \int_{Q_{T,R,\delta}} \sqrt{f} (t,v) g (t,w) \phi (t,v,w) \dt \dv \dw
\end{equation}
as $n_i \to \infty$, where
\[
g(t,w)=2\nabla_w\sqrt{f^\kappa_+(t,w)}\,.
\]

(To see that 
\[
g^{n_i}\rightharpoonup g\quad\text{ in }L^2((0,T)\times B_{R+\delta})\,,
\]
use \eqref{e:ae-conv} and \eqref{e:fn-mass} to prove that $\sqrt{f^{n_i,\kappa}_+}\to\sqrt{f^\kappa_+}$ in $L^1_{loc}((0,+\infty)\times\R^3)$, together with the formula $g^n(t,w)=2\nabla_w\sqrt{f^n(t,w)}\un_{f^n(t,w)\ge\kappa}$ and \eqref{e:gradfn}.)

By Fubini's theorem, the left-hand side of \eqref{e:distribution-limit} becomes
\[  
\int_{(0,T) \times B_{R+\delta}} g^{n_i} (t,w)  \left\{ \int_{B_R} \sqrt{f^{n_i}} (t,v) \phi (t,v,w)  \dv \right\} \dt \dw. 
\]
We then prove that 
\[
\int_{B_R} \sqrt{f^{n_i}} (t,v) \phi (t,v,w)  \dv\to\int_{B_R} \sqrt{f} (t,v) \phi (t,v,w)  \dv\quad\text{ in }L^2 ((0,T) \times B_{R+\delta})\,.
\]
This strong $L^2$ convergence implies in particular \eqref{e:distribution-limit}. We henceforth drop the subscript $i$ in the subsequence indexed by $n_i$ for simplicity.  

To get the $L^2$ convergence of $\int_{B_R} \sqrt{f^n} (t,v) \phi (t,v,w)  \dv$, we first prove that it converges a.e. Observe that $\sqrt{f^n}(t,v) \phi (t,v,w)$ converges a.e. in $B_R$ to $\sqrt{f}(t,v) \phi (t,v,w)$ (recall \eqref{e:ae-conv}), while uniform 
integrability of $v \mapsto \sqrt{f^n}(t,v) \phi (t,v,w)$  follows from \eqref{e:fn-mass}. 

We next use \eqref{e:fn-mass} to get
\begin{align*}
\int_{(0,T) \times B_{R+\delta}} \left\{ \int_{B_R} \sqrt{f^n} (t,v) \phi (t,v,w)  \dv \right\}^{2+2\eps} \dt \dw \le\int_{(0,T) \times B_{R+\delta}} \left\{ 2 M_0 \int_{B_R} \phi^2 (t,v,w)  \dv \right\}^{1+\eps} \dt \dw\,.
\end{align*}
This implies that
\[
\lim_{n \to \infty} \left\| \int_{B_R} \sqrt{f^n} (t,v) \phi (t,v,w) \dv \right\|_{L^2((0,T) \times B_{R+\delta})} = \left\| \int_{B_R} \sqrt{f} (t,v) \phi (t,v,w) \dv \right\|_{L^2((0,T) \times B_{R+\delta})}\,.
\]
With a.e. convergence proved above, this implies $L^2$ convergence \cite[Proposition~2.4.6]{willem-book}, and proves \eqref{e:distribution-limit}.

Moreover, the sequence $\{ \sqrt{f^n} (t,v) g^n (t,w)  \}$ is bounded in $L^2(Q_{T,R,\delta})$,
\[ 
\int_{Q_{T,R,\delta}} f^n (t,v) \frac{|\nabla_v \fnkp|^2}{f^n} (t,w) \dt \dv \dw \le 4M_0\left\| \nabla_v\sqrt{\fnkp}\right\|^2_{L^2 ((0,T) \times B_{R+\delta})} \le 4M_0 \mathcal{I}(T,R+\delta)
\]
thanks to \eqref{e:gradfn}, and therefore weakly relatively compact in $L^2 (Q_{T,R,\delta})$. With \eqref{e:distribution-limit} and by uniqueness of the limit, our claim \eqref{e:claim-L2weak} is proved. 

Summarizing, we have proved that
\[
\begin{aligned}
\sqrt{\anin (v-w)}\sqrt{ f^n (t,v) f^n (t,w)}  \left( \nabla_v L^n_+ (t,v) \Psi (t,v)- \nabla_w L^n_+ (t,w) \Psi (t,w) \right)
\\
\rightharpoonup
F^{\kappa,in}_+\Psi(t,v)+\sqrt{\ain(v,w)}\sqrt{f(t,v)f(t,w)}\nabla_wL_+(t,w)(\Psi(t,v)-\Psi(t,w))
\\
=
\sqrt{\ain (v-w) }\sqrt{f (t,v) f (t,w)}  \left( \nabla_v L_+ (t,v) \Psi (t,v)- \nabla_w L_+ (t,w) \Psi (t,w) \right)
\end{aligned}
\]
in $L^2((0,T)\times B_R\times B_{R+\delta})$ for all $T,R>0$ and all $\delta\in(0,\bar\delta]$, and hence
\[
\begin{aligned}
\limsup_{n_i\to\infty}\int_{t_1}^{t_2}\iint_{\R^3 \times \R^3}\!\!-\tfrac12\anin (v\!-\!w) f^{n_i} (t,v) f^{n_i} (t,w) : \left( \nabla_v L^{n_i}_+ (t,v) \Psi (t,v)\!-\! \nabla_w L^{n_i}_+ (t,w) \Psi (t,w) \right)^{\otimes 2}\dv\dw \dt
\\
\le\int_{t_1}^{t_2}\iint_{\R^3 \times \R^3}\!\!-\tfrac12\ain (v\!-\!w) f (t,v) f (t,w) : \left( \nabla_v L_+ (t,v) \Psi (t,v)\!-\! \nabla_w L_+ (t,w) \Psi (t,w) \right)^{\otimes 2}\dv\dw \dt&\,.
\end{aligned}
\]

\paragraph{Short-range interactions (error term).} Integrating by parts in $v$ and in $w$ in \eqref{e:secondterm}, one finds that
\begin{multline*}
\iint_{\R^3 \times \R^3}  \anin (v-w)  (\Psi (t,v) - \Psi (t,w))^2 \nabla_w f^{n,\kappa}_+ (t,w) \cdot \nabla_v f^{n,\kappa} _+ (t,v)   \dv \dw  
\\
= \iint_{\R^3 \times \R^3}    \mathfrak{A}_n (v,w) f^{n,\kappa}_+ (t,w) f^{n,\kappa} _+ (t,v)   \dv \dw\,,
\end{multline*}
with
\begin{align*}
\mathfrak{A}_n (v,w) := & (-\nabla^2_z : \anin) (v-w) \bigg(\Psi (v) - \Psi(w) \bigg)^2 
\\
&- 2 (\nabla_z \cdot \anin)(v-w)(\Psi (v) - \Psi (w)) \cdot \bigg(\nabla_v \Psi (v) +\nabla_w \Psi (w)\bigg) 
\\
&- 2 \anin (v-w) \nabla_v \Psi (v) \cdot \nabla_w \Psi (w)\,.
\end{align*}
By Lemma~\ref{l:d2anin}, where the first and second derivatives of $\anin$ are computed, observing that 
\[
X'\left(\tfrac{|z|}\delta\right)\not=0\implies\tfrac12<\tfrac{|z|}\delta<1\implies n|z|=n\delta\tfrac{|z|}\delta>1\implies X(n|z|)=1\,,
\]
we find that
\[ 
\mathfrak{A}_n (v,w) = 2\ka (|v-w|) n X'(n |v-w|) \frac{(\Psi (v)-\Psi(w))^2}{|v-w|^2}+ X(n|v-w|)\mathfrak{A} (v,w)\,, 
\]
with $\mathfrak{A}$ given by
\begin{align*}
\mathfrak{A} (v,w) = & (-\nabla^2_z : \ain) (v-w) \bigg(\Psi (v) - \Psi(w) \bigg)^2 
\\
&- 2 (\nabla_z \cdot \ain)(v-w)(\Psi (v) - \Psi (w)) \cdot \bigg(\nabla_v \Psi (v) +\nabla_w \Psi (w)\bigg) 
\\
&- 2 \ain (v-w) \nabla_v \Psi (v) \cdot \nabla_w \Psi (w)\,.
\end{align*}
In particular, we have
\begin{align*}
|\mathfrak{A} (v,w)| \le & \phantom{+} 2 \| \nabla_v \Psi\|_\infty^2 \bigg( \ka'(|v-w|)  +3 \ka (|v-w|) \delta^{-1} \bigg) \un_{\{|v-w|\le \delta \}} 
\\
&+ 8 \| \nabla_v \Psi\|_\infty^2 \ka (|v-w|) \frac{\un_{\{|v-w|\le \delta \}}}{|v-w|} 
\\
& + 2 \| \nabla_v \Psi\|_\infty^2 \ka (|v-w|) \frac{\un_{\{|v-w|\le \delta \}}}{|v-w|}\,.
\end{align*}
Using now that $\ka'(r) \le \ka (r)/r$ and the monotonicity of $\ka$  shows that
\[
|\mathfrak{A} (v,w)| \le 18 \| \nabla_v \Psi\|_\infty^2 \ka (\delta) \frac{\un_{\{|v-w|\le \delta \}}}{|v-w|}\un_{v\in B_R}\un_{w\in B_{R+\delta}}\,,
\]
assuming that $\Psi$ has support in $B_R$. Thus
\[
\sup_{n_i}\int_{t_1}^{t_2}\iint_{B_R\times B_{R+\delta}}(|\mathfrak{A} (v,w)|f^{n_i,\kappa}_+ (t,w) f^{n_i,\kappa} _+ (t,v))^{1+\eps}\dv\dw\dt<\infty\,,
\]
by using Lemma~\ref{l:source-nl} with $\rho=1+\eps$, with $\tfrac1p=\tfrac{2-\eps}{3(1+\eps)}$ and $\tfrac1q=\tfrac{1+4\eps}{2(1+\eps)}$ so that $\tfrac1p+\tfrac2{3q}=1$, and with $\beta=1+\eps$ and $\tfrac1\theta=\tfrac{1-2\eps}{2(1+\eps)}$ so 
that $\tfrac1\beta+\tfrac2{3\theta}>1$, together with Lemma \ref{l:interpolation} and \eqref{e:gradfn}-\eqref{e:fn-mass} to obtain an upper bound for the term above. This implies the uniform integrability of 
$\mathfrak{A} (v,w) f^{n_i,\kappa}_+ (t,w) f^{n_i,\kappa} _+ (t,v)$, and since $f^{n_i,\kappa}_+\to f^\kappa_+$ a.e. as $n_i\to\infty$ by \eqref{e:ae-conv}, we conclude that
\[
\begin{aligned}
\lim_{n_i\to \infty}\int_{t_1}^{t_2} \iint_{\R^3 \times \R^3}  X(n|v-w])\mathfrak{A} (v,w) f^{n_i,\kappa}_+ (t,w) f^{n_i,\kappa} _+ (t,v) \dt  \dv \dw 
\\
=   \int_{t_1}^{t_2} \iint_{\R^3 \times \R^3} \mathfrak{A} (v,w) f^{\kappa}_+ (t,w) f^{\kappa} _+ (t,v) \dt  \dv \dw&\,.
\end{aligned}
\]
On the other hand, observing that $n|v-w|X'(n|v-w|)\le 3\cdot\un_{|v-w|<1/n}$ with $1/n<\delta/2$ implies that
\[
0\le\mathfrak{A}_n (v,w)-X(n|v-w|)\mathfrak{A} (v,w)\le 6k(\delta/2)\|\nabla\Psi\|^2_{L^\infty}\frac{\un_{|v-w|<1/n}}{|v-w|}\,.
\]
Then
\[
\int_{t_1}^{t_2} \iint_{\R^3 \times \R^3}  (\mathfrak{A}_n (v,w)-X(n|v-w|)\mathfrak{A} (v,w)) f^{n_i,\kappa}_+ (t,w) f^{n_i,\kappa} _+ (t,v) \dt  \dv \dw \to 0
\]
as $n_i\to\infty$. (Indeed, the inequality above shows that the integrand converges to $0$ a.e. for $v\not=w$, while uniform integrability follows  by the same argument as for the preceding limit, since $\un_{|v-w|<1/n}<\un_{|v-w|<\delta}$.)

\paragraph{Short-range interactions (conclusion).}
We have proved that
\[
\limsup_{n_i \to \infty} \int_{t_1}^{t_2} \mathcal{T}_1^{+,n_i}(t)\dt \le \int_{t_1}^{t_2} \mathcal{T}_1^+ (f)(t) \dt\,,
\]
with
$$
\begin{aligned}
\mathcal{T}_1^+ (f)(t):=\int_{t_1}^{t_2}\iint_{\R^3 \times \R^3}-\tfrac12\ain (v-w) f (t,v) f (t,w) : \left( \Psi\nabla_v L_+ (t,v) - \Psi\nabla_w L_+ (t,w) \right)^{\otimes 2}\dv\dw \dt
\\
+\tfrac12\int_{t_1}^{t_2} \iint_{\R^3 \times \R^3} \mathfrak{A} (v,w) f^{\kappa}_+ (t,w) f^{\kappa} _+ (t,v) \dt  \dv \dw&\,.
\end{aligned}
$$
We can now rewrite the limit $\mathcal{T}_1^+ (f)$ in a different form. Undoing the symmetrization in the limit of the term \eqref{e:firstterm} and the integration by parts in the limit of the term \eqref{e:secondterm}, we find that
\begin{align*}
\mathcal{T}_1^+ (f)(t) =&   - \iint_{\R^3 \times \R^3} \ain (v-w) f (t,v) f (t,w) \left( \nabla_v L_+ (t,v)- \nabla_w L_+ (t,w)\right) \cdot \nabla_v L_+ (t,v) \varphi (t,v)  \dv \dw 
\\
=&  - \tfrac12 \iint_{\R^3 \times \R^3} \left|F^{\kappa,\inn}_+ (t,v,w)\right|^2 \varphi (t,v) \dv \dw 
\\
& + \tfrac12 \iint_{\R^3 \times \R^3} F^{\kappa,\inn}_+ (t,v,w) \cdot \sqrt{\ain (v-w)}\sqrt{ f(t,v) f(t,w)} \nabla_w L_+ (t,w) \bigg(\varphi (t,v)- \varphi (t,w) \bigg)   \dv \dw\,,
\\
\intertext{recalling that $\nabla_w L_+ (t,w)= \frac{\nabla_w \fkp}{f}(t,w)$, we finally get}
\mathcal{T}_1^+ (f)(t) =&   - \tfrac12 \iint_{\R^3 \times \R^3} \left|F^{\kappa,\inn}_+ (t,v,w)\right|^2 \varphi (t,v) \dv \dw 
\\
& + \tfrac12 \iint_{\R^3 \times \R^3} F^{\kappa,\inn}_+ (t,v,w) \cdot \bigg((\varphi (t,v)- \varphi (t,w)) \sqrt{\ain (v-w)}\bigg)\sqrt{f}(t,v)  \frac{\nabla_w \fkp}{\sqrt{f}} (t,w)   \dv \dw\,.
\end{align*}

\paragraph{Long-range interactions.} We next study $\mathcal{T}_1^{\out,n}$. Since $\aout$ is smooth and decays at infinity, we can write
\begin{equation}
\label{e:t1out}
\begin{aligned}
\mathcal{T}_1^{\out,n} = & -  \iint_{\R^3 \times \R^3} A^\out_n (t,v) : \frac{(\nabla_v \fnkp)^{\otimes 2}}{f^n}  (t,v) \varphi (t,v)  \dv 
\\
&+ \iint_{\R^3 \times \R^3} \aout (v-w) \nabla_w f^n (t,w)  \cdot  \nabla_v \fnkp (t,v) \, \varphi (t,v)  \dv \dw\,.
\end{aligned}
\end{equation}
We study the limit of each term separately.

Write the first term as
\[ 
\frac{(\nabla_v \fnkp)^{\otimes 2}}{f^n} = 4   \left( \nabla_v  \sqrt{f^n}   \right)^{\otimes 2} \un_{f^n \ge \kappa}\,.
\]
Since 
\[
\un_{f^{n_i}\ge \kappa} \nabla_v \sqrt{f^{n_i}}=\nabla_v\left(\sqrt{f^{n_i}}-\sqrt{\kappa}\right)_+ \rightharpoonup\nabla_v\left(\sqrt{f}-\sqrt{\kappa}\right)_+=\un_{f \ge \kappa} \nabla_v \sqrt{f}
\]
in $L^2((0,T) \times B_R)$ (see \eqref{e:gradfn-conv}), and since $A^\out_{n_i}\to A^\out$ locally uniformly by Lemma~\ref{l:anout} and definite positive by Lemma~\ref{l:diff-mat}, one has
\[
\sqrt{A^\out_{n_i}} \un_{f^{n_i} \ge \kappa}\nabla_v \sqrt{f^{n_i}}\rightharpoonup\sqrt{A^\out} \un_{f \ge \kappa}  \nabla_v \sqrt{f}
\]
in $L^2 ((0,T) \times B_R)$. Hence
\begin{multline*}
\limsup_{n_i \to \infty}\int_{t_1}^{t_2} \iint_{\R^3 \times \R^3} -A^\out_{n_i} (t,v) : \frac{(\nabla_v f^{n_i,\kappa}_+)^{\otimes 2}}{f^{n_i}}  (t,v) \varphi (t,v)  \dt \dv 
\\
\le \int_{t_1}^{t_2}  \iint_{\R^3 \times \R^3} -A^\out (t,v) : \frac{(\nabla_v f^{\kappa}_+)^{\otimes 2}}{f}  (t,v) \varphi (t,v)  \dt \dv\,.
\end{multline*}

Integrating by parts with respect to $v$ and  $w$ in the second term in \eqref{e:t1out} yields,
\begin{align*}
\int_{t_1}^{t_2} \iint_{\R^3 \times \R^3}& \aout (v-w) \nabla_w f^n (t,w)  \cdot  \nabla_v \fnkp (t,v) \varphi (t,v) \dt \dv \dw 
\\
=& -\int_{t_1}^{t_2} \iint_{\R^3 \times \R^3} (\nabla_z \cdot \aout) (v-w) f^n (t,w)  \fnkp (t,v) \cdot \nabla_v \varphi (t,v) \dt \dv \dw 
\\
& -\int_{t_1}^{t_2} \iint_{\R^3 \times \R^3} \nabla_z^2 : \aout (v-w) f^n (t,w)  \fnkp (t,v) \varphi (t,v) \dt \dv \dw 
\\
= & -\int_{t_1}^{t_2} \int_{\R^3} (\nabla_v \cdot A^\out_n) (t,v)  \fnkp (t,v) \cdot \nabla_v \varphi (t,v) \dt \dv  
\\
& -\int_{t_1}^{t_2} \int_{\R^3} \nabla_v^2 : A^\out_n (t,v) \fnkp (t,v) \varphi (t,v) \dt \dv\,.
\end{align*}
Using Lemma~\ref{l:anout} shows that the second term in \eqref{e:t1out} converges along the subsequence $n_i$ from \eqref{e:ae-conv}. Summarizing, we have proved that
\[ 
\limsup_{n_i \to \infty } \mathcal{T}_1^{\out,n_i}\le \mathcal{T}_1^\out (f)
\] 
with
\begin{equation}\label{T1out}
\begin{aligned}
\mathcal{T}_1^\out (f)(t):=&-\int_{\R^3} A^\out (t,v) : \frac{(\nabla_v f^{\kappa}_+)^{\otimes 2}}{f}  (t,v) \varphi (t,v)  \dv
\\
&+ \iint_{\R^3 \times \R^3} \aout (v-w) \nabla_w f (t,w)  \cdot  \nabla_v \fkp (t,v) \varphi (t,v) \dv \dw\,.
\end{aligned}
\end{equation}

\paragraph{The depleted nonlinearity.} Since $\un_{f^n\ge\kappa}\nabla_vf^n=\nabla_v(f^n -\kappa)_+=\nabla_vf^{n,\kappa}_+$ while $\un_{f^n<\kappa}\nabla_vf^n=\nabla_v(f^n\wedge\kappa)$,
\[ 
\mathcal{T}_1^{\wedge,n}  =   \iint_{\R^3 \times \R^3} \anin (v-w)  \nabla_w (f^n \wedge \kappa)  (t,w) \cdot  \nabla_v f^{n,\kappa}_+ (t,v) \varphi (t,v) \dv \dw\,.
\]
Integrating by parts with respect to $v$ and $w$,
\begin{multline}
\label{e:different-form}
\mathcal{T}_1^{\wedge,n} =  \iint_{\R^3 \times \R^3} \bigg( (-\nabla_z^2 : \anin) (v-w) \varphi (t,v) 
\\
- (\nabla_z \cdot \anin) (v-w) \cdot \nabla_v \varphi (t,v) \bigg)  (f^n \wedge \kappa)  (t,w) f^{n,\kappa}_+ (t,v) \dt \dv \dw\,.
\end{multline}
Recalling the formula for $\nabla_z^2 \colon \anin$ given by \eqref{e:d2anin}, we write $\mathcal{T}_1^{\wedge}= \mathcal{T}_{1,0}^\wedge+\mathcal{T}_{1,1}^\wedge+ \mathcal{T}_{1,2}^\wedge + \mathcal{T}_{1,3}^\wedge$ 
with
\begin{align*}
\mathcal{T}_{1,0}^{\wedge,n } &= 2 \iint_{\R^3 \times \R^3}  \frac{\ka' (|v-w|)X(n|v-w|) (1-X(\delta^{-1}|v-w|))}{|v-w|^2} (f^n \wedge \kappa)  (t,w)  f^{n,\kappa}_+ (t,v) \varphi (t,v) \dt \dv \dw\,,
\\
\mathcal{T}_{1,1}^{\wedge,n } &= 2 \iint_{\R^3 \times \R^3}  \frac{nX'(n|v-w|)\ka (|v-w|)}{|v-w|^2} (f^n \wedge \kappa)  (t,w)  f^{n,\kappa}_+ (t,v) \varphi (t,v) \dt \dv \dw\,,
\\
\mathcal{T}_{1,2}^{\wedge,n } &= - 2 \iint_{\R^3\times \R^3} \frac{\ka (|v-w|)}{|v-w|^2} \delta^{-1}X'(\delta^{-1}|v-w|) (f^n \wedge \kappa)  (t,w) f^{n,\kappa}_+ (t,v) \varphi (t,v) \dt \dv \dw\,,
\\
\mathcal{T}_{1,3}^{\wedge,n } &= 2 \iint_{\R^3\times \R^3} X(n|v-w|) \frac{(1-X(\frac{|v-w|}\delta))\ka (|v-w|) }{|v-w|^2} \tfrac{v-w}{|v-w|}\cdot \nabla_v \varphi (t,v) (f^n \wedge \kappa)  (t,w) f^{n,\kappa}_+ (t,v) \dt \dv \dw\,.
\end{align*}

First
\begin{equation}
\label{e:t10wv}
\begin{aligned}
\lim_{n_i \to \infty}\int_{t_1}^{t_2} \mathcal{T}_{1,0}^{\wedge,n_i} (f^{n_i}) \dt = \int_{t_1}^{t_2} \mathcal{T}_{1,0}^{\wedge} (f) \dt\,,\quad\text{ where}
\\
\mathcal{T}_{1,0}^{\wedge} (f):=2 \iint_{\R^3 \times \R^3}  \frac{\ka' (|v-w|)(1-X(\delta^{-1}|v-w|))}{|v-w|^2} (f\wedge \kappa)  (t,w)  f^{\kappa}_+ (t,v) \varphi (t,v) \dt \dv \dw&\,.
\end{aligned}
\end{equation}
Indeed, \eqref{e:ae-conv} shows that the integrand converges a.e. in $\Omega_{R,\delta} = \{ (t,v,w) : t \in (t_1,t_2), v \in B_R, |v-w|<\delta\}$. To prove uniform integrability, write
\begin{align*}
\int_{\Omega_{R,\delta}} &\left(\frac{\ka'(|v-w|)}{|v-w|^{2}}\right)^{1+\eps}   |\varphi (t,v)|^{1+\eps}(f^n \wedge \kappa)^{1+\eps}  (t,w) (f^n-\kappa)_+^{1+\eps} (t,v) \dt \dv \dw 
\\
& \le\|\varphi \|_\infty^{1+\eps} \kappa^{1+\eps}   \int_{[t_1,t_2]\times B_R} \left\{ \int_{|v-w| < \delta} \left(\frac{\ka'(|v-w|)}{|v-w|^{2}}\right)^{1+\eps} \dw \right\}    (f^n-\kappa)_+^{1+\eps} (t,v) \dt \dv  
\\
& \le\|\varphi \|_\infty^{1+\eps} \kappa^{1+\eps}  \left\{ \int_{B_1} \left(\frac{\ka'(|z|)}{|z|^{2}}\right)^{1+\eps} \dz \right\} \int_{[t_1,t_2]\times B_R} (f^n-\kappa)_+^{1+\eps} (t,v) \dt \dv\,. 
\end{align*}
Recall that $\ka'(|z|)/|z|^2\in L^{1+\iota} (B_1)$. With $\eps = \min (\iota,2/3)$, using \eqref{e:fn-equi} gives the required uniform integrability.

The argument for $\mathcal{T}^{\wedge,n}_{1,3}$ is very similar except that we first use again that $\ka$ is nondecreasing to write
\[
(1-X(\delta^{-1}|v-w|))\ka (|v-w|) \le \ka (\delta) \un_{B_\delta}(v-w)\,.
\]
We then easily get
\begin{equation}
\label{e:t13wv}
\begin{aligned}
\lim_{i \to \infty}  \int_{t_1}^{t_2} \mathcal{T}_{1,3}^{\wedge,n_i} \dt = \int_{t_1}^{t_2} \mathcal{T}_{1,3}^{\wedge} (f) \dt\,,\quad\text{ where}
\\
\mathcal{T}_{1,3}^{\wedge} (f):=2 \iint_{\R^3 \times \R^3}\frac{(1-X(\frac{|v-w|}\delta))\ka (|v-w|) }{|v-w|^2} \tfrac{v-w}{|v-w|}\cdot \nabla_v \varphi (t,v) (f \wedge \kappa)  (t,w) f^{\kappa}_+ (t,v) \dt \dv \dw&\,.
\end{aligned}
\end{equation}

There is no singularity in $\mathcal{T}_{1,2}^{\wedge,n}$, and thanks to \eqref{e:ae-conv} and \eqref{e:fn-equi}, we easily see that
\begin{equation}
\label{e:t12wv}
\begin{aligned}
\lim_{i \to \infty}\int_{t_1}^{t_2} \mathcal{T}_{1,2}^{\wedge,n_i}\dt = \int_{t_1}^{t_2} \mathcal{T}_{1,2}^{\wedge} (f) \dt\,,\quad\text{ where}
\\
\mathcal{T}_{1,2}^{\wedge} (f):=2 \iint_{\R^3 \times \R^3}\frac{\ka (|v-w|)}{|v-w|^2} \delta^{-1}X'(\delta^{-1}|v-w|) (f\wedge \kappa)  (t,w) f^{\kappa}_+ (t,v) \varphi (t,v) \dt \dv \dw&\,.
\end{aligned}
\end{equation}

Finally, in order to study the upper limit in $\mathcal{T}_{1,1}^{\wedge,n}$, we first observe that 
\begin{align*}
\int_{\R^3}\frac{nX'(n|v-w|) \ka (|v-w|)}{|v-w|^2}  (f^n \wedge \kappa)  (t,w) \dw &\le  \kappa \ka (1/n)\int_{\R^3}  \frac{nX'(n|v-w|)}{|v-w|^2} \dw = 4\pi \ka (1/n)\kappa\,.
\end{align*}
(We used that $\ka$ is nonincreasing while $X'\ge 0$ is supported in $[1/2,1]$ and $X(1)-X(1/2)=1$.) Hence
\begin{equation}
\label{e:t11wv}  
\begin{aligned}
\limsup_{n_i\to\infty}\int_{t_1}^{t_2}\mathcal{T}_{1,1}^{\wedge,n }\dt\le&\lim_{n_i\to\infty}8\pi\ka(1/n_i)\kappa\int_{t_1}^{t_2}\int_{\R^3} f^{n_i,\kappa}_+ (t,v) \varphi (t,v)\dv\dt 
\\
=&8\pi\ka(0)\kappa\int_{t_1}^{t_2}\int_{\R^3} f^{\kappa}_+ (t,v) \varphi (t,v) \dv \dt=:\int_{t_1}^{t_2}\mathcal T^\wedge_{1,1}(f)\dt\,.
\end{aligned}
\end{equation}

Observe that
\[
\begin{aligned}
\mathcal{T}_{1,0}^{\wedge } (f) + \mathcal{T}_{1,1}^{\wedge } (f) + \mathcal{T}_{1,2}^{\wedge } (f) = 8 \pi \ka (0) \kappa \int_{\R^3} \fkp (t,v) \varphi (t,v) \dv
\\
+2 \iint_{\R^3 \times \R^3}  \frac{\ka' (|v-w|)(1-X(\frac{|v-w|}\delta))-\ka(|v-w|)\frac1\delta X'(\frac{|v-w|}\delta)}{|v-w|^2} (f\wedge \kappa)  (t,w)  f^{\kappa}_+ (t,v) \varphi (t,v) \dt \dv \dw
\\
=8 \pi \ka (0) \kappa \int_{\R^3} \fkp (t,v) \varphi (t,v) \dv+2 \iint_{\R^3 \times \R^3}  \frac{\ka' (|v-w|)}{|v-w|^2} (f\wedge \kappa)  (t,w)  f^{\kappa}_+ (t,v) \varphi (t,v) \dt \dv \dw
\\
+\iint_{\R^3 \times \R^3}\nabla^2 : \aout (v-w) (f \wedge \kappa)  (t,w) \fkp (t,v)\varphi (t,v)\dv \dw&\,,
\end{aligned}\]
while
\[
\begin{aligned}
\mathcal{T}_{1,3}^{\wedge} (f):=2 \iint_{\R^3 \times \R^3}\frac{(1-X(\frac{|v-w|}\delta))\ka (|v-w|) }{|v-w|^2} \tfrac{v-w}{|v-w|}\cdot \nabla_v \varphi (t,v) (f \wedge \kappa)  (t,w) f^{\kappa}_+ (t,v) \dt \dv \dw
\\
=-\iint_{\R^3 \times \R^3}(\nabla\cdot\ain)(v-w)\cdot \nabla_v \varphi (t,v) (f \wedge \kappa)  (t,w) f^{\kappa}_+ (t,v) \dt \dv \dw&\,.
\end{aligned}
\]
Hence
\begin{equation}\label{T1wedge}
\begin{aligned}
\limsup_{n_i\to\infty}\mathcal T_1^{\wedge,n_i}\le\mathcal{T}_{1,0}^{\wedge } (f) + \mathcal{T}_{1,1}^{\wedge } (f) + \mathcal{T}_{1,2}^{\wedge } (f) + \mathcal{T}_{1,3}^{\wedge} (f)
\\
=-\iint_{\R^3 \times \R^3}(\varphi(t,v)\nabla^2 : \ain (v-w)+\nabla_v \varphi (t,v)\cdot(\nabla\cdot\ain)(v-w))(f \wedge \kappa)  (t,w) \fkp (t,v)\dv \dw&\,.
\end{aligned}
\end{equation}
Elementary transformations put
\[
\mathcal T_1(f):=\mathcal{T}_1^+ (f)+\mathcal{T}_{1,0}^{\wedge } (f) + \mathcal{T}_{1,1}^{\wedge } (f) + \mathcal{T}_{1,2}^{\wedge } (f) + \mathcal{T}_{1,3}^{\wedge} (f)+\mathcal{T}_1^\out (f)
\]
in the form appearing in the statement of Lemma ~\ref{l:t1-lim}.

Observe indeed that the second term in the last right-hand side of \eqref{T1wedge} is the last integral in the definition of $\mathcal E_1(f)$ in Lemma ~\ref{l:t1-lim}. Similarly, it will be convenient to put the first term in the last right-hand side 
of \eqref{T1wedge} in the form
\[
\begin{aligned}
-\iint_{\R^3 \times \R^3}\varphi(t,v)\nabla^2 : \ain (v-w)(f \wedge \kappa)  (t,w) \fkp (t,v)\dv \dw
\\
=-\iint_{\R^3 \times \R^3}\varphi(t,v)\nabla^2 : a(v-w)(f \wedge \kappa)  (t,w) \fkp (t,v)\dv \dw
\\
+\iint_{\R^3 \times \R^3}\varphi(t,v)\nabla^2 : \aout (v-w)(f \wedge \kappa)  (t,w) \fkp (t,v)\dv \dw&\,.
\end{aligned}
\]
The second term in the right-hand side of the equality above is the fourth integral in the definition of $\mathcal E_1(f)$ in Lemma ~\ref{l:t1-lim}, while the first term on the same right-hand corresponds to the third and fourth terms in the last
right-hand side in the definition of $\mathcal T_1(f)$, according to \eqref{e:hess-a}.

Observe also that the last term in the right-hand side of \eqref{T1out} satisfies
\[
\iint_{\R^3 \times \R^3} \aout (v-w) \nabla_w f (t,w)  \cdot  \nabla_v \fkp (t,v) \varphi (t,v) \dv \dw=\int_{\R^3}(\nabla_v\cdot A^{out})(t,v)\cdot\nabla_v \fkp (t,v) \varphi (t,v) \dv
\]
which can be transformed into the second and third integrals on the right-hand side of the definition of $\mathcal E_1(f)$ in Lemma ~\ref{l:t1-lim} after integrating by parts. The remaining identifications needed in order to complete the proof 
of Lemma ~\ref{l:t1-lim} are obvious.
\end{proof}

\medskip
We next turn to the proof of   Lemma~\ref{l:t2-lim}.

\begin{proof}[Proof of Lemma~\ref{l:t2-lim}]
We  first put $\mathcal{T}_2^n$ in the form
\[
\mathcal{T}_2^n= -\iint_{\R^3 \times \R^3}a_n(v-w) f^n (t,v) f^n(t,w) \left(\frac{\nabla_v f^n (t,v)}{f^n(t,v)} - \frac{\nabla_w f^n (t,w)}{f^n (t,w)} \right) \cdot \lnp \left(\frac{f^n}{\kappa}\right) (t,v) \nabla_v \varphi (t,v)\dv \dw\,.
\]
Splitting again $a_n$ into $\anin+\aout$ with $n > 2/\delta$ (so that $a^{\out}_n=\aout$), we write $\mathcal{T}_2^n=\mathcal{T}_2^{\inn,n}+\mathcal{T}_2^{\out,n}$ with
\[
\begin{aligned}
\mathcal{T}_2^{\inn,n}:=& -\iint_{\R^3 \times \R^3}\anin(v-w) f^n (t,v) f^n(t,w) \left(\tfrac{\nabla_v f^n (t,v)}{f^n(t,v)} - \tfrac{\nabla_w f^n (t,w)}{f^n (t,w)} \right) \cdot \lnp \left(\tfrac{f^n}{\kappa}\right) (t,v) \nabla_v \varphi (t,v)\dv \dw\,,
\\
\mathcal{T}_2^{\out,n}:=& -\iint_{\R^3 \times \R^3}\aout(v-w) f^n (t,v) f^n(t,w) \left(\tfrac{\nabla_v f^n (t,v)}{f^n(t,v)} - \tfrac{\nabla_w f^n (t,w)}{f^n (t,w)} \right) \cdot \lnp \left(\tfrac{f^n}{\kappa}\right) (t,v) \nabla_v \varphi (t,v)\dv \dw\,.
\end{aligned}
\]

\paragraph{Short-range interactions.} As far as $\mathcal{T}_2^{\inn,n}$ is concerned, since $\anin(z) = (1-X(\delta^{-1}|z|)) a_n (z)$, we have
\[
\int_{t_1}^{t_2} \mathcal{T}_2^{\inn,n} \dt = -\int_{t_1}^{t_2} \iint_{\R^3 \times \R^3}  F^n (t,v,w)  \cdot G^n (t,v,w) \nabla_v \varphi (t,v) \dt \dv \dw \,, 
\]
with
\begin{align*}
F^n (t,v,w) &\,= \sqrt{a_n (v-w)}\sqrt{ f^n (t,v) f^n (t,w)} \left( \frac{\nabla_v f^n}{f^n} (t,v) - \frac{\nabla_w f^n}{f^n} (t,w) \right)\,,
\\
G^n (t,v,w) &:= (1-X(\delta^{-1}|v-w|))  \sqrt{f^n (t,v) f^n(t,w)} \lnp \left(\frac{f^n(t,v)}{\kappa}\right)  \sqrt{ a_n(v-w)}\,.
\end{align*}

\begin{claim}\label{claim}
Let  $\Omega_{R,\delta} = \{ (t,v,w) : t \in (t_1,t_2), v \in B_R, |v-w|<\delta\}$. The sequence $G^n$ converges (strongly) in $L^2(\Omega_{R,\delta})$ to 
\[ 
G^\inn(t,v,w)=       (1-X(\delta^{-1}|v-w|))  \sqrt{f (t,v) f(t,w)}  \lnp \left(\frac{f(t,v)}{\kappa}\right) \sqrt{ a(v-w)}\,.
\]
\end{claim}

\noindent
Obviously $G^\inn=0$ in $(t_1,t_2)\times B_R\times\R^3\setminus\Omega_{R,\delta}$. Recalling that $\ain (z) = (1-X(\delta^{-1}|z|)) a(z)$, this claim, together with \eqref{e:dissip-conv}, implies in particular that
\begin{equation}\label{e:t2in-conv}  
\begin{aligned} 
\lim_{n_i\to\infty}\int_{t_1}^{t_2}\mathcal{T}_2^{\inn,n_i}\dt = \int_{t_1}^{t_2}\mathcal{T}_2^{\inn} (f) \dt\,,\quad\text{ where}
\\
\mathcal{T}_2^{\inn} (f):=-\iint_{\R^3 \times \R^3}\ain(v-w) f(t,v) f(t,w) \left(\frac{\nabla_v f(t,v)}{f(t,v)} - \frac{\nabla_w f(t,w)}{f(t,w)} \right) \cdot \lnp \left(\frac{f}{\kappa}\right) (t,v) \nabla_v \varphi (t,v)\dv \dw&\,.
\end{aligned}
\end{equation} 

We are thus left with the task of justifying the claim above.

\begin{proof}[Proof of the claim]
First, we use \eqref{e:ae-conv} to see that $G^n \to G$ almost everywhere in $\Omega_{R,\delta}$. Next we prove that the sequence $(G^n)^2$ is uniformly integrable:
\[
\begin{aligned}
\iiint_{\Omega_{R,\delta}} & (G^n)^{2(1+\eps)} (t,v,w) \dt \dv \dw 
\\
&\le\ka(\delta)^{1+\eps}\int_0^T \int_{B_R} \left\{ (f^n)^{1+\eps}(t,\cdot) \ast_v \frac{\un_{B_\delta}}{|\cdot|^{1+\eps}} \right\} (f^n)^{1+\eps} \left( \lnp \left(\frac{f^n}{\kappa}\right)  \right)^{2(1+\eps)}(t,v)  \dt \dv 
\\
&\le\ka(\delta)^{1+\eps} C_{HLS,1+\eps}\|f^n\|^{1+\eps}_{L^{q_\eps}_tL^{p_\eps}_v((0,T) \times (B_{R+\delta}))} \|\Phi^n\|^{1+\eps}_{L^{\theta_\eps}_t L^{1+\eps}_v((0,T) \times B_R)} \,,
\end{aligned}
\]
where we have used Lemma~\ref{l:source-nl} with $F = \Phi^n := f^n \left(\lnp \left(\frac{f^n}{\kappa}\right)\right)^2$, with $p_\eps:=\tfrac{3(1+\eps)}{2-\eps}$ and $q_\eps:=\tfrac{2(1+\eps)}{1+4\eps}$, and with $\beta_\eps=1+\eps$ and 
$\theta_\eps=\tfrac{2(1+\eps)}{1-4\eps}$ to get the last inequality.

Moreover, $|\lnp r| \le C_\eps r^{\eps/2}$ for some absolute constant $C_\eps>0$, so that 
\[
\begin{aligned}
\iiint_{\Omega_{R,\delta}} & (G^n)^{2(1+\eps)} (t,v,w) \dt \dv \dw 
\\
&\le\ka(\delta)^{1+\eps} C_{HLS,1+\eps}C_\eps^{1+\eps}\|f^n\|^{1+\eps}_{L^{q_\eps}_tL^{p_\eps}_v((0,T) \times (B_{R+\delta}))} \|f^n\|^{(1+\eps)^2}_{L^{\theta_\eps(1+\eps)}_t L^{(1+\eps)^2}_v((0,T) \times B_R)}\,.
\end{aligned}
\]
Our choice of $p_\eps,q_\eps$ is such that $\tfrac1{p_\eps}+\tfrac2{3q_\eps}=1$, while $\tfrac1{(1+\eps)^2}+\tfrac2{3\theta_\eps(1+\eps)}=\tfrac{4-4\eps}{3(1+\eps)^2}>1$ for $\eps>0$ small enough. At this point, we use Lemma \ref{l:interpolation}
with \eqref{e:fn-mass} and \eqref{e:gradfn} to conclude that
\[ 
\sup_{n\ge 1}\int_0^T  \iint_{B_R \times B_{R+\delta}} (G^n)^{2(1+\eps)} (t,v,w) \dt \dv \dw \le C[\eps,T,R,\delta]
\]
for some constant $C[\eps,T,R,\delta]>0$. Thus $(G^n)^2$ is uniformly integrable in $(0,T)\times B_R \times B_{R+\delta}$ and since it converges a.e. towards $G^2$ (along a subsequence), strong $L^2$ convergence follows.
\end{proof}

\paragraph{Long-range interactions.}
As far as $\mathcal{T}_2^{\out,n}$ is concerned, we have
\begin{align}
\nonumber  
\int_{t_1}^{t_2} \mathcal{T}_2^{\out,n} \dt =&-\int_{t_1}^{t_2} \iint_{\R^3 \times \R^3} \aout (v-w) f^n(t,w) \nabla_v f^n (t,v) \cdot \lnp \left(\frac{f^n}{\kappa}\right) (t,v) \nabla_v \varphi (t,v) \dt \dv \dw 
\\
\nonumber    
&+\int_{t_1}^{t_2} \iint_{\R^3 \times \R^3} \aout (v-w) f^n(t,v) \nabla_w f^n (t,w) \cdot \lnp \left(\frac{f^n}{\kappa}\right) (t,v) \nabla_v \varphi (t,v) \dt \dv \dw 
\\
\label{e:t2out-conv}  
=& -2 \int_{t_1}^{t_2} \int_{\R^3}  \nabla_v \sqrt{f^n} (t,v) \cdot A^\out_n (t,v)  \sqrt{f^n (t,v)} \lnp \left(\frac{f^n}{\kappa}\right) (t,v) \nabla_v \varphi (t,v) \dt \dv 
\\
\nonumber 
&+ \int_{t_1}^{t_2} \int_{\R^3} (\nabla_v \cdot A^\out_n) (t,v)  \cdot f^n(t,v)  \lnp \left(\frac{f^n}{\kappa}\right) (t,v) \nabla_v \varphi (t,v) \dt \dv \,,
\end{align}
with $A^\out_n (t,v)= \int_{\R^3} \aout (v-w) f^n(t,w) \dw$. 

Recall that $\nabla_v \sqrt{f^{n_i}}$ converges weakly in $L^2 ((0,T) \times B_R)$ towards $\nabla_v \sqrt{f}$ by \eqref{e:gradfn-conv}. We claim that $A^\out_{n_i}   \sqrt{f^{n_i} } \lnp \left(\frac{f^{n_i}}{\kappa}\right)$ converges strongly 
in $L^2 ((0,T) \times B_R)$ towards $A^\out  \sqrt{f}\lnp \left(\frac{f}{\kappa}\right)$. Indeed, thanks to Lemma~\ref{l:anout}, $A^\out_{n_i}$ converges almost everywhere in $(0,T) \times B_R$ to $A^\out$, and therefore
$A^\out_{n_i}  \sqrt{f^{n_i}} \lnp \left(\frac{f^{n_i}}{\kappa}\right)$ converges a.e. in $(0,T) \times B_R$ towards $A^\out \sqrt{f} \lnp \left(\frac{f}{\kappa}\right)$. Moreover, the square of its Euclidian norm  is uniformly integrable since, by
Lemma~\ref{l:anout},
\begin{align*}
\int_{(0,T) \times B_R} \left|A^\out_n  \sqrt{f^n} \lnp \left(\frac{f^n}{\kappa}\right)\right|^3 \dt \dv & \le C_\delta \int_{(0,T) \times B_R} \left|\sqrt{f^n} \lnp \left(\frac{f^n}{\kappa}\right) \right|^3 \dt \dv \,,
\\
\intertext{use now that $|\lnp r| \le C_\eps r^{\eps}$ for some constant $C_\eps$ depending only on $\eps >0$,}
& \le C_\delta C_\eps^3 \kappa^{- 3 \eps} \int_{(0,T) \times B_R} (f^n)^{\frac32 +  3\eps}    \dt \dv\,.
\end{align*}
One conclude thanks to \eqref{e:fn-equi} by choosing $\eps=1/18$. 

The weak convergence of $\nabla_v \sqrt{f^{n_i}}$, and the strong convergence of $A^\out_{n_i}   \sqrt{f^{n_i} } \lnp \left(\frac{f^{n_i}}{\kappa}\right) $ imply that
\begin{multline}\label{e:t2out-conv2}
\lim_{n_i\to\infty}\int_{t_1}^{t_2} \int_{\R^3}  \nabla_v \sqrt{f^{n_i}} (t,v) \cdot A^\out_{n_i} (t,v)  \sqrt{f^{n_i} (t,v)} \lnp \left(\frac{f^{n_i}}{\kappa}\right) (t,v) \nabla_v \varphi (t,v) \dt \dv \dw
\\ 
=\int_{t_1}^{t_2} \int_{\R^3}  \nabla_v \sqrt{f} (t,v) \cdot A^\out (t,v)  \sqrt{f (t,v)} \lnp \left(\frac{f}{\kappa}\right) (t,v) \nabla_v \varphi (t,v) \dt \dv \dw\,.
\end{multline}

We turn to the second term appearing in \eqref{e:t2out-conv}. Thanks to Lemma~\ref{l:anout}, we know that $\nabla_v \cdot A^\out_{n_i}$ converges  towards $\nabla_v A^\out$ a.e. in $(0,T) \times B_R$. Moreover, arguing as above and 
using again Lemma~\ref{l:anout}, we prove that $\left|(\nabla_v \cdot A^\out_n) (t,v)  \cdot f^n(t,v)  \lnp \left(\frac{f^n}{\kappa}\right) (t,v)\right|^{3/2}$ is bounded in $L^1((0,T)\times B_R)$, and conclude that
\begin{multline}\label{e:t2out-conv3}   
\lim_{n_i\to\infty}\int_{t_1}^{t_2} \int_{\R^3} (\nabla_v \cdot A^\out_{n_i}) (t,v)  \cdot f^{n_i}(t,v)  \lnp \left(\frac{f^{n_i}}{\kappa}\right) (t,v) \nabla_v \varphi (t,v) \dt \dv \dw 
\\ 
=\int_{t_1}^{t_2} \int_{\R^3} (\nabla_v \cdot A^\out) (t,v)  \cdot f(t,v)  \lnp \left(\frac{f}{\kappa}\right) (t,v) \nabla_v \varphi (t,v) \dt \dv \dw\,.
\end{multline}

Combining \eqref{e:t2out-conv}, \eqref{e:t2out-conv2} and \eqref{e:t2out-conv3} shows that
\begin{equation}\label{e:t2out-conv4}  
\begin{aligned} 
\lim_{n_i\to\infty}\int_{t_1}^{t_2}\mathcal{T}_2^{\out,n_i} (f^{n_i}) \dt=\int_{t_1}^{t_2}\mathcal{T}_2^{\out} (f) \dt\,,\quad\text{ where}
\\
\mathcal{T}_2^{\out} (f):= \int_{\R^3} \left(-A^\out (t,v)  \nabla_v f (t,v)+f(t,v)\nabla_v \cdot A^\out(t,v)\right)\lnp \left(\frac{f}{\kappa}\right) (t,v) \nabla_v \varphi (t,v)\dv&\,.
\end{aligned}
\end{equation} 
Combining \eqref{e:t2out-conv4} and \eqref{e:t2in-conv} yields the desired result. Indeed, one easily checks that
\[
\mathcal{T}_2^{\inn} (f)=-\iint_{\R^3 \times \R^3}F(t,v,w)\cdot G(t,v,w)\nabla_v\varphi(t,v)\dv \dw\,,
\]
which is the first integral in the last right-hand side of the definition of $\mathcal T_2(f)$ in Lemma~\ref{l:t2-lim}. 

On the other hand, integrating by parts the first summand in the integrand of the right-hand side of \eqref{e:t2out-conv4}, after observing that $\lnp \left(\frac{f}{\kappa}\right)\nabla_v f=\nabla_vh^\kappa_+(f)$, shows that
\[
\begin{aligned}
\int_{\R^3} \left(-A^\out (t,v)  \nabla_v f (t,v)+f(t,v)\nabla_v \cdot A^\out(t,v)\right)\lnp \left(\frac{f}{\kappa}\right) (t,v) \nabla_v \varphi (t,v)\dv
\\
=\int_{\R^3} \nabla_v\cdot A^\out(t,v)\nabla_v\varphi (t,v)\left(h^\kappa_+(f)(t,v)+f\lnp \left(\frac{f}{\kappa}\right) (t,v)\right)\dv+\int_{\R^3} A^\out (t,v):\nabla_v^2\varphi (t,v)h^\kappa_+(f)(t,v)\dv
\end{aligned}
\]
which are respectively the third and the second integrals in the last right-hand side of the definition of $\mathcal T_2(f)$ in Lemma~\ref{l:t2-lim}.
\end{proof}


\subsection{A technical lemma}

\begin{lemma}\label{l:iteration}
Given $\beta > 1$, for  all $n \ge 1$, one has
\[ 
\sum_{i=0}^{n-1} (n-i) \beta^{i} = (\beta-1)^{-2} (\beta^{n+1} - (n+1) (\beta-1)-1)\,.
\]
\end{lemma}

\begin{proof}
Let $\sigma (\beta)$ denote $\sum_{i=0}^{n-1} (n-i) \beta^{i}$. We observe that \(\sigma(\beta ) = \frac{\partial \Sigma}{\partial \alpha} (1,\beta)\)
with 
\[ 
\Sigma (\alpha, \beta ) = \sum_{i=0}^{n-1} \alpha^{n-i} \beta^{i} = \alpha^n \sum_{i=0}^{n-1}(\beta/\alpha)^i = \alpha^n \frac{(\beta/\alpha)^{n}-1}{\beta/\alpha-1} = \alpha \frac{\beta^{n}-\alpha^{n}}{\beta- \alpha}.
\]
We now differentiate and set $\alpha=1$ to get the result. 
\end{proof}


\def\cprime{$'$}

\end{document}